\documentclass[11pt]{article}
\usepackage{amsmath}
\usepackage{bm}
\usepackage{amsfonts}
\usepackage{amssymb,amsmath,amsthm, amsfonts}
\usepackage{mathrsfs}

\def\sqr#1#2{{\vcenter{\vbox{\hrule height.#2pt
              \hbox{\vrule width.#2pt height#1pt \kern#1pt \vrule width.#2pt}
          \hrule height.#2pt}}}}
\def\sqr#1#2{{\vcenter{\vbox{\hrule height.#2pt
              \hbox{\vrule width.#2pt height#1pt \kern#1pt \vrule width.#2pt}
              \hrule height.#2pt}}}}
\def\3n{\negthinspace \negthinspace \negthinspace }
\def\2n{\negthinspace \negthinspace }
\def\1n{\negthinspace }

\def\={\buildrel \triangle \over =}

\def\esssup{\mathop{\rm esssup}}

\def\min{\mathop{\rm min}}
\def\exp{\mathop{\rm exp}}
\def\sup{\mathop{\rm sup}}
\def\inf{\mathop{\rm inf}}

\def\({\Big (}
\def\){\Big )}
\def\[{\Big[}
\def\]{\Big]}
\def\be{\begin{equation}}

\def\ee{\end{equation}}

\def\square#1{\vbox{\hrule\hbox{\vrule height#1%
     \kern#1\vrule}\hrule}}
\def\rectangle#1#2{\vbox{\hrule\hbox{\vrule height#1%
     \kern#2\vrule}\hrule}}

\font\tenbb=msbm10 \font\sevenbb=msbm7 \font\fivebb=msbm5

\newfam\bbfam
\scriptscriptfont\bbfam=\fivebb \textfont\bbfam=\tenbb
\scriptfont\bbfam=\sevenbb

\newtheorem{lemma}{Lemma}[section]
\newtheorem{remark}{Remark}[section]

\newtheorem{theorem}{Theorem}[section]

\newtheorem{definition}{Definition}[section]
\newtheorem{proposition}{Proposition}[section]

\makeatletter
   
   \@addtoreset{equation}{section}
\makeatother

\begin{document}

\title{A Global Stochastic Maximum Principle for Mean-Field Forward-Backward Stochastic Control Systems with Quadratic Generators}
	\author{Rainer Buckdahn$^{1,2}$,\,\, Juan Li$^{3,2,\ast\ast}$,\,\, Yanwei Li$^{3,\ast\ast}$,\,\, Yi Wang$^{3,\ast\ast}$\\
		{\small $^1$Laboratoire de Math\'{e}matiques de Bretagne Atlantique, Univ Brest, }\\
		{\small UMR CNRS 6205, 6 avenue Le Gorgeu, 29200 Brest, France.}\\
		{\small $^2$Research Center for Mathematics and Interdisciplinary Sciences,}\\
		{\small         Shandong University, Qingdao 266237, P.~R.~China.}\\
				{\small $^3$School of Mathematics and Statistics, Shandong University, Weihai, Weihai 264209, P.~R.~China.}\\
		{\small {\it E-mails: rainer.buckdahn@univ-brest.fr; juanli@sdu.edu.cn; 201916496@mail.sdu.edu.cn;}}\\
		{\small {\it wangyi429@mail.sdu.edu.cn.}}\\
	}
\date{April 10, 2024}

\renewcommand{\thefootnote}{\fnsymbol{footnote}}
\footnotetext[1]{Juan Li is supported by the NSF of P.R. China (NOs. 12031009, 11871037), National Key R and D Program of China (NO. 2018YFA0703900), NSFC-RS (No. 11661130148; NA150344).\\
 \mbox{ } \  \ $^{\ast\ast}$Corresponding authors.}

	\maketitle
	
	\textbf{Abstract}. Our paper is devoted to the study of Peng's stochastic maximum principle (SMP) for a stochastic control problem composed of a controlled forward stochastic differential equation (SDE) as dynamics and a controlled backward SDE which defines the cost functional. Our studies combine the difficulties which come, on one hand, from the fact that the coefficients of both the SDE and the backward SDE are of mean-field type (i.e., they do not only depend on the control process and the solution processes but also on their law), and on the other hand, from the fact that the coefficient of the BSDE is of quadratic growth in $Z$. Our SMP is novel, it extends in a by far non trivial way existing results on SMP.

\\

	
\section{Introduction}
	The studies of our paper are devoted to a recursive stochastic mean-field control problem, i.e., more precisely, a stochastic control problem whose dynamics are given by a controlled stochastic differential equation (SDE), and a controlled backward stochastic differential equation (BSDE) which defines the cost functional. The coefficients of both the SDE and the BSDE do not only depend on control process and the solution processes but also on their law, and that of the BSDE has in addition a quadratic growth in $Z$. Supposing the existence of a control process $u^*$ minimizing the cost functional over all control processes taking their values in an arbitrarily given not necessary convex subset of $\mathbb{R}^n$, our objective is to give a necessary condition for this optimality of $u^*$, that is to get the stochastic maximum principle (SMP) for our stochastic mean-field control problem.
	
	The study of the SMP for stochastic control problems has been one of the main objectives in the literature. Peng \cite{P1990} was the first to derive a SMP for general stochastic control problems with a not necessary convex control state space. Since this seminal paper there have been a lot of works by different authors who have extended Peng's SMP to more general stochastic control problems; the reader can refer, for instance, to \cite{YZ1999},  \cite{H2017}, \cite{HJX2022} and \cite{SW2006} as well as the references therein.
	
	In his recent work \cite{H2017} Hu extended Peng's SMP to stochastic control problems whose cost functional is defined through a controlled BSDE. With his work in which he considered two adjoint equations, he succeeded to solve a problem proposed by Peng \cite{P1998}. The studies of Hu were extended by Hu, Ji and Xu \cite{HJX2022} to decoupled forward-backward stochastic control problems with a coefficient driving the BSDE which has a quadratic growth in $Z$. The BSDE with a generator with quadratic growth in $Z$ were studied the first time by Kobylanski \cite{K2000} in 2000; while she supposed the terminal condition to be bounded, Briand and Hu \cite{BH2006,BH2008} proved the existence and the uniqueness for one-dimensional BSDE with unbounded terminal and convex generator. While here bounded mean oscillation martingales (BMO martingales) are not involved yet, Briand and Confortola \cite{BC2008} investigated BSDE with a generator which satisfies a stochastic Lipschitz condition involving BMO martingale. The approach by \cite{HJX2022} also uses BMO martingales as a crucial tool to study the SMP.
	
	Above works as well as recent studies of mean-field forward-backward SDEs have heavily stimulated us in our work. Mean-field BSDEs were first introduced and studied by Buckdahn, Djehiche, Li and Peng \cite{BDLP2009} and Buckdahn, Li and Peng \cite{BLP2009} in 2009. Using Lions' notion of the derivative of a function over the 2-Wasserstein space of probability laws, Buckdahn, Li and Ma \cite{BLM2016} studied the global stochastic maximum principle for a class of controlled equations under Lipschitz conditions for coefficients of mean-field type. Later Buckdahn, Chen and Li \cite{BCL2021} investigated the partial derivative with respect to a probability measure and studied the SMP of an associated stochastic control problem. Their coefficients depend on the joint distribution of state and control process. In the study of our stochastic control problem we are concerned by controlled mean-field BSDEs whose generating coefficient has a quadratic growth. Mean-field BSDEs with quadratic growth have been studied by different authors, see e.g., \cite{HHTW2022} concerning the one-dimensional case and that of a bounded terminal condition.
	
	Our investigations combine the difficulties related with those of the study of a SMP for a stochastic control problem given in form of a controlled mean-field forward-backward SDE and those related with mean-field BSDEs with quadratic growth. For this, we need to have new methods to overcome the difficulties.
	
More precisely, in our paper we consider the following stochastic control problem consisting of a system of mean-field forward and backward SDEs:
\begin{equation*}
		\left\{
		\begin{aligned}
			&dX_t^u=b(t,X_t^u,P_{X_t^u},u_t)dt+\sigma(t,X_t^u,P_{X_t^u},u_t)dB_t,\ 0\leq t\leq T,\\
			&dY_t^u=- f(t,X_t^u,Y_t^u,Z_t^u,P_{(X_t^u, Y_t^u)},u_t)dt+Z_t^udB_t,\ 0\leq t\leq T,\\
			&X_0^u=x_0,\ Y_T^u=\Phi(X_T^u,P_{X_T^u}),
		\end{aligned}
		\right.
\end{equation*}
whose cost functional is defined by the backward SDE
$$J(u(\cdot)):=Y_0^u.$$
Here $u=(u_t)$ runs the admissible control processes taking their values in an arbitrarily given non empty subset of $\mathbb{R}^n.$
Our goal is to derive Peng's SMP which is satisfied by the optimal control process at which this cost functional attains its minimum.
	
	In our approach we first study the existence and the uniqueness for a general linear mean-field BSDE involving an unbounded BMO coefficient for $Z$ (see Proposition 3.1 concerning an a priori estimate and Theorem 3.1; these results extend in a non trivial way Theorem 10 of \cite{BC2008}). This BSDE turns out to be a crucial tool in our studies.
	We first consider the existence and the uniqueness of the linear mean-field BSDEs with unbounded coefficient. This is important for our variational and  adjoint equations.
	In what follows we introduce the spike variations of the optimal control process $u^*$ and give the adjoint equations for the controlled state equations as well as their estimates. Proposition 4.1 presents the estimates of the Taylor expansion of state process controlled by the spike-varied optimal control. For their proof (namely that of (4.8)), Proposition 5.2 in \cite{BCL2021} plays a crucial role. However, for our following computations the assumptions in Proposition 5.2 \cite{BCL2021} turn out to be too restrictive (see Remark 4.1), and we have to study a novel, non trivial extension, which is given by Proposition 4.2. Its subtle and rather technical proof is postponed to the Appendix.
	
	After the Taylor expansion of the solution process of the forward SDE governed by the spike-varied optimal control, that of the corresponding solution $Y$ of the BSDE is essential. It is given by Theorem 5.1 and involves namely the solutions of the variational equations and those of the BSDEs adjoint to them, but also BSDE (5.1). Since the forward and the backward state equations in our control system are of mean-field type, also the adjoint equations and (5.1) are of mean-field type, see (4.12)-(4.15), and the quadratic growth of the BSDE of the control problem implies that the both adjoint BSDEs and (5.1) have stochastic Lipschitz conditions (that of $Z$ is of BMO type). It's here, where Proposition 3.1 and Theorem 3.1 play a crucial role.
	
Notice that for the Taylor expansion of $Y$, among the different crucial estimates that of\\  $E\[\(\displaystyle\int_0^T|Z^*_t||\overline{\overline{Z}}_t|I_{E_\varepsilon}(t)dt\)^p\]$ plays a central role, where $(Y^*,Z^*)$ is the solution of the BSDE of the control problem associated to the optimal control $u^*$, $(\overline{\overline{Y}},\overline{\overline{Z}})$ is the solution of BSDE (5.1) and $E_\varepsilon\subset [0,T]$ is the set on which the optimal control is spike-varied. This estimate has also a central place in the computations in \cite{HJX2022}. However, (3.41) in \cite{HJX2022} is not clearly redacted. A discussion with the authors has given that (3.41) in \cite{HJX2022} can be obtained, but only by using Lebesgue's differentiation theorem, which has the consequence that (3.41) holds for $E_\varepsilon=[t_0,t_0+\varepsilon]$ only $dt_0$-a.e., for all $\varepsilon\in(0,\varepsilon_{t_0})$, for some $\varepsilon_{t_0}>0$ and a constant $C=C_{t_0}$ depending on $t_0$. However, these constraints are not mentioned there nor in their following computations. For more clearness and for simplicity, we give a different method to overcome the difficulties. We introduce the deterministic set $\Gamma_M:=\{t\in[0,T]:E[|Z^*_t|^2]\leq M\},\ M\geq1$, and  redefined $E_\varepsilon:=E_\varepsilon\cap\Gamma_M$. The estimates made in the preceding sections continue to hold true, those for $Y$ are made with the new $E_\varepsilon$. This leads to an SMP which first holds true for $t\in\Gamma_M$, and taking the limit $M\rightarrow +\infty$ gives the general result, see Theorem 6.1. Let us emphasize that our SMP also generalises  that obtained by \cite{HJX2022}, see Remark 6.1.
	
	Our paper is organized as follows. In Section 2, we introduce the different spaces we need. Some properties of BMO martingale and the definition of the differentiability with respect to a probability measure are also recalled. Finally, the framework of the optimal control problem that we study is given. Section 3 is devoted to the existence and the uniqueness for general mean-field BSDEs with stochastic Lipschitz conditions and the corresponding estimates of mean-field BSDE. In Section 4, we study the first- and the second-order variational equations, and we give some estimates with respect to forward controlled equation and its variational equations. In Section 5, we study the expansion of the solution of the backward controlled equation with spike-varied optimal control. The  corresponding maximum principle is obtained in Section 6. In the Appendix, we present the proof of Proposition 4.2.

\section{Preliminaries}
Let $T>0$ be a fixed time horizon and $(\Omega,\mathcal{F},P)$ be a given complete probability space. Let $\{B_t,\ 0\leq t\leq T\}$ be a standard Brownian motion defined on $(\Omega,\mathcal{F},P)$, with valued in $\mathbb{R}^d$. We assume that there is a sub-$\sigma$-field $\mathcal{F}_0\subset\mathcal{F}$, containing all $P$-null subsets of $\mathcal{F}$, such that\\
\indent\quad(i)\ \ the  Brownian motion $B$ is independent of $\mathcal{F}_0$;

\indent\quad(ii)\ $\mathcal{F}_{0}$ is ``rich enough'', i.e., $\mathcal{P}_{2}(\mathbb{R}^{k})=\{P_{\xi},\ \xi\in L^{2}(\mathcal{F}_{0};\mathbb{R}^{k})\},\ k\geq1$.\\
\noindent Here $P_\xi := P\circ [\xi]^{-1}$ denotes the law of the random variable $\xi$ under the probability $P$.

We denote by $\mathbb{F}=(\mathcal{F}_{t})_{t\in [0,T]}$ the filtration generated by $B$
and augmented by $\mathcal{F}_0$.

\indent Let us introduce some notations and concepts, which are used frequently in what follows.
For $k\in \mathbb{N}$\ and\ $x,y\in\mathbb{R}^{k}$, we denote the norm and the inner product, respectively,
by $|x|=\Big(\sum\limits_{i=1}^{k}x_{i}^{2}\Big)^\frac{1}{2}$, and $\langle x,y\rangle=\sum\limits_{i=1}^{k}x_{i}y_{i}.$
We shall introduce the following spaces, for $1<p,\ q<\infty$:\\
\noindent$\bullet \ L^p(\mathcal{F}_T;\mathbb{R}^k)$ is the set of $\mathcal{F}_T$-measurable random variables $ \xi: \Omega\rightarrow\mathbb{R}^k$ such that $\|\xi\|_{L^{p}}:=$\\
\noindent\mbox{ } \ $\left (E\left[|\xi|^{p}\right]\right)^{\frac{1}{p}}<\infty $.\\
\noindent$\bullet \ L^\infty(\mathcal{F}_T;\mathbb{R}^k)$ is the set of $\mathcal{F}_T$-measurable random variables $ \xi: \Omega\rightarrow\mathbb{R}^k$ such that $\|\xi\|_{\infty}:=$\\
\noindent\mbox{ } \ $\esssup\limits_{\omega\in\Omega}|\xi(\omega)|<\infty $.\\
\noindent$\bullet\ \mathcal{H}^{p,q}(0,T;\mathbb{R}^k)$ is the set of $\mathbb{F}$-progressively measurable processes $\eta: \Omega\times[0,T]\rightarrow \mathbb{R}^k$ with $\|\eta\|_{\mathcal{H}^{p,q}}:=$\\
\noindent\mbox{ } \ $\( E\[\(\displaystyle\int_{0}^{T}|\eta(s)|^{p}ds\)^{\frac{q}{p}}\]\)^{\frac{1}{q}}<\infty$.\\
In particular, we write $\mathcal{H}^{p}(0,T;\mathbb{R}^k)$ when $p=q$.\\
\noindent$\bullet\ \mathcal{S}^p(0,T;\mathbb{R}^k)$ is the set of $\mathbb{F}$-adapted continuous processes  $\eta: \Omega\times[0,T]\rightarrow \mathbb{R}^k$ with
$\|\eta\|_{\mathcal{S}^{p}}:=$\\
\noindent\mbox{ } \ $\(E\[\sup\limits_{0\leq s\leq T}|\eta(s)|^{p}\]\)^{\frac{1}{p}}<\infty$.\\
\noindent$\bullet\ \mathcal{S}^\infty(0,T;\mathbb{R}^k)$ is the set of $\mathbb{F}$-adapted continuous processes  $\eta: \Omega\times[0,T]\rightarrow \mathbb{R}^k$ with
$\|\eta\|_{\infty}:=$\\
\noindent\mbox{ } \ $\esssup\limits_{(\omega,s)\in\Omega\times[0,T]}|\eta_s(\omega)|<\infty$.\\
\noindent$\bullet\ \mathcal{P}_{2}(\mathbb{R}^{k})$ is the set of the probability measures $\mu$ on $(\mathbb{R}^{k},\mathcal{B}(\mathbb{R}^{k}))$ with finite second moment, i.e.,\\
\noindent\mbox{ } \ $\displaystyle\int_{\mathbb{R}^{k}}|x|^2\mu (dx)<\infty$.\\
Here let $\mathcal{P}_{2}(\mathbb{R}^{k})$ be endowed with the $2$-Wasserstein metric: For $\mu,\ \nu\in \mathcal{P}_{2}(\mathbb{R}^k)$,
\begin{equation}\label{eq2.1}
W_{2}(\mu,\nu):=\inf\Big\{\Big(\int_{\mathbb{R}^k\times\mathbb{R}^k}|x-y|^{2}\rho(dxdy)\Big)^{\frac{1}{2}}: \rho\in\mathcal{P}_{2}(\mathbb{R}^{2k}),\ \rho(.\times\mathbb{R}^k)=\mu,\  \rho(\mathbb{R}^k\times.)=\nu\Big\}.
\end{equation}
\noindent$\bullet$ For any $p\geq1,\ BMO_p$ is the set of real-valued $\mathbb{F}$-martingales $M$ such that
\begin{equation}\label{eq2.2}
\|M\|_{BMO_p}:=\sup\limits_{\tau}\Bigr\| (E[|M_T-M_\tau|^p|\mathcal{F}_\tau])^{\frac{1}{p}}\Bigr\|_{\infty}<\infty.
\end{equation}
\noindent\mbox{ } \ where the supremum is taken over all stopping times $\tau\in[0,T]$.

For $M\in BMO_2$ it is called a $BMO$ martingale.

In the following lemma we state the properties of $BMO$ martingales, we refer readers to \cite{HWY1992} or \cite{K1994} for more details.
\begin{lemma}\label{lem2.1}
\begin{itemize}
\item[\rm{(i)}] The energy inequality: Given a BMO martingale $M$, for any $n\geq1$,
$$E[\langle M\rangle_T^n]\leq n!\|M\|^{2n}_{BMO_2}.$$
\item[\rm{(ii)}] If $M\in BMO_2$,  the stochastic exponential $\mathcal{E}(M)_{\cdot}:=\exp{\(M_{\cdot}-\frac{1}{2}\langle M\rangle_{\cdot}\)}$ is a uniform integrable martingale.
\item[\rm{(iii)}] The John-Nirenberg inequality: Let $M\in BMO_2$ and $\delta\in\(0,\|M\|_{BMO_2}^{-2}\)$. Then we have
                    \begin{equation*}
                    E\[\exp\{\delta(\langle M\rangle_T-\langle M\rangle_\tau)\}|\mathcal{F}_\tau\]\leq
                    \(1-\delta\|M\|_{BMO_2}^2\)^{-1},
                    \end{equation*}
     for all stopping times $\tau\in[0,T]$.
\item[\rm{(iv)}] The reverse H\"{o}lder inequality: Let $M\in BMO_2$. Let $\Psi:\ (1,+\infty)\rightarrow\mathbb{R}_+$ be defined by
$$\Psi(x)=\(1+\frac{1}{x}\log\frac{2x-1}{2(x-1)}\)^{\frac{1}{2}}-1.$$
The function $\Psi$ is  nonincreasing $\lim\limits_{x\rightarrow+\infty}\Psi(x)=0$ and $\lim\limits_{x\rightarrow1}\Psi(x)=+\infty$. Let $p_M$ be such that $\Psi(p_M)=\|M\|_{BMO_2}$. For all $p\in(1,p_M)$, and for all stopping time $\tau\in[0,T]$, we have
\begin{equation*}
E\[\mathcal{E}(M)_T^p/\mathcal{E}(M)_\tau^p|\mathcal{F}_\tau\]\leq K(p,\|M\|_{BMO_2}),\ \mathrm{P}\text{-}\mathrm{a.s.},
\end{equation*}
where
\begin{equation*}
K(p,\|M\|_{BMO_2})=2\(1-\frac{2p-2}{2p-1}\exp\bigr\{p^2(\|M\|^2_{BMO_2}+2\|M\|_{BMO_2})\bigr\}\)^{-1}.
\end{equation*}
By $p_M^\ast$ we denote the conjugate exponent of $p_M$, that is, $(p_M)^{-1}+(p_M^\ast)^{-1}=1$.
\end{itemize}
\end{lemma}

\subsection{Derivative with respect to a measure over a Banach space}\label{derivative}
Now let us recall the notions of derivative with respect to probability measures.
Let $\mathcal{K}$ be a real separable Banach space endowed with the norm $|\cdot|_{\mathcal{K}}$. We denote by $\mathcal{K}^\prime=\{l:\ \mathcal{K}\rightarrow\mathbb{R}\ |\ l\ \text{is the bounded linear functional}\}$ the dual space of $\mathcal{K}$ and we endow $\mathcal{K}^\prime\times\mathcal{K}$ with the duality product $\langle l,x\rangle_{\mathcal{K}^\prime\times\mathcal{K}}=l(x),\ l\in \mathcal{K}^\prime,\ x\in\mathcal{K}$.
Recall that
$$\mathcal{P}_2(\mathcal{K})=\Bigr\{m\in \mathcal{P}(\mathcal{K})\ \text{is the probability measure over }(\mathcal{K},\mathcal{B}(\mathcal{K})):\ \int_\mathcal{K}|x|^2_\mathcal{K}m(dx)<\infty\Bigr\}.$$
We now recall the notion of differentiability of a function $f:\ \mathcal{P}_2(\mathcal{K})\rightarrow\mathbb{R}$. 
\begin{definition}\label{def1}
We say that $f:\mathcal{P}_2(\mathcal{K})\rightarrow\mathbb{R}$ has the linear functional derivative $\mathcal{D}_mf:\mathcal{P}_2(\mathcal{K})\times\mathcal{K}\rightarrow\mathbb{R}$, if $\mathcal{D}_mf$ is a continuous function over $\mathcal{P}_2(\mathcal{K})\times\mathcal{K}$ with at most quadratic growth such that, for all $m,\ m^\prime\in\mathcal{P}_2(\mathcal{K})$,
\begin{equation}\label{eq2.13}
f(m^\prime)-f(m)=\int_0^1\int_\mathcal{K}\mathcal{D}_mf(tm^\prime+(1-t)m,x)(m^\prime(dx)-m(dx))dt.
\end{equation}
\end{definition}
Let us suppose that for $f:\mathcal{P}_2(\mathcal{K})\rightarrow\mathbb{R}$ the derivative $\mathcal{D}_mf:\mathcal{P}_2(\mathcal{K})\times\mathcal{K}\rightarrow\mathbb{R}$ exists, is continuous and of at most quadratic growth, and that, for all $m\in\mathcal{P}_2(\mathcal{K})$, $\mathcal{D}_mf(m,\cdot):\mathcal{K}\rightarrow\mathbb{R}$ is differentiable, i.e., there exists $\partial_x(\mathcal{D}_mf)(m,\cdot):\mathcal{K}\rightarrow\mathcal{K}^\prime$ such that, for all $x\in\mathcal{K}$,
\begin{equation}\label{eq2.14}
\mathcal{D}_mf(m,y)-\mathcal{D}_mf(m,x)
=\langle\partial_x(\mathcal{D}_mf)(m,x),y-x\rangle_{\mathcal{K}^\prime\times\mathcal{K}}
    +o(|y-x|_{\mathcal{K}}),\ \text{as}\ \mathcal{K}\ni y\rightarrow x.
\end{equation}

\begin{proposition}\label{prop2.1}
We assume that $\partial_x(\mathcal{D}_mf)$ exists and is continuous and of at most linear growth.
Then the following first order Taylor expansion for $f:\ \mathcal{P}_2(\mathcal{K})\rightarrow\mathbb{R}$ at $m\in\mathcal{P}_2(\mathcal{K})$ holds true:
\begin{equation}\label{eq2.15}
f(m^\prime)=f(m)+\int_\mathcal{K}\mathcal{D}_mf(m,x)(m^\prime(dx)-m(dx))
    +o(W_2(m,m^\prime)),\ \text{as}\ W_2(m,m^\prime)\rightarrow0\ (m^\prime\in\mathcal{P}_2(\mathcal{K})).
\end{equation}
\end{proposition}
The proof of this proposition is a slight extension of the corresponding result in Carmona and Delarue \cite{CD2018}. Moreover, also under our the above assumptions the following result concerning Lion's L-derivative of $f: \mathcal{P}_2(\mathcal{K})\rightarrow\mathbb{R}$ extends easily from $\mathcal{K}=\mathbb{R}^n$\ (see Carmona and Delarue \cite{CD2018}) to the general separable Banach space $\mathcal{K}$.
\begin{proposition}\label{eq2.16}
Given any $\xi\in L^2(\Omega,\mathcal{F},P;\mathcal{K})$, we have
\begin{equation}\label{eq2.17}
f(P_{\xi+\eta})=f(P_\xi)
    +E[\langle\partial_x(\mathcal{D}_mf)(P_\xi,\xi),\eta\rangle_{\mathcal{K}\times\mathcal{K}^\prime}]
    +R(\xi,\eta),\ \eta\in L^2(\Omega,\mathcal{F},P;\mathcal{K}),
\end{equation}
where $R(\xi,\eta)=o((E[|\eta|^2_\mathcal{K}])^\frac{1}{2})$, as $(E[|\eta|^2_\mathcal{K}])^\frac{1}{2}\rightarrow0$.
\end{proposition}
\begin{definition}\label{def2}
The L-derivative $(\partial_\mu f(m,x))$, if it exists, is defined by
\begin{equation}\label{eq2.18}
\partial_\mu f(m,x)=\partial_x(\mathcal{D}_mf)(m,x),\ (m,x)\in \mathcal{P}_2(\mathcal{K})\times\mathcal{K}.
\end{equation}
Moreover, when we speak about the differentiability of $f:\ \mathcal{P}_2(\mathcal{K})\rightarrow\mathbb{R}$, we mean the L-differentiability.
\end{definition}

If, for all fixed $x\in\mathcal{K}$, the map $m\rightarrow\mathcal{D}_mf(m,x)$ has the linear functional derivative, then we say that $f$ has the second-order linear functional derivative and denote it by $\mathcal{D}^2_{mm}f$.  By Definition \ref{def1} we have that $\mathcal{D}^2_{mm}f:\mathcal{P}_2(\mathcal{K})\times\mathcal{K}\times\mathcal{K}\rightarrow\mathbb{R}$ with
\begin{equation*}
\mathcal{D}_{m}f(m^\prime,x)-\mathcal{D}_{m}f(m,x)
=\int_0^1\int_\mathcal{K}\mathcal{D}^2_{mm}f(tm^\prime+(1-t)m,x,x^\prime)(m^\prime(dx^\prime)-m(dx^\prime))dt.
\end{equation*}
If $f$ has the second-order linear functional derivative and  $\mathcal{D}^2_{mm}f=(\mathcal{D}^2_{mm}f(m,x,x^\prime))$ is twice differentiable in the variables $(x,x^\prime)$, then we set
$$\partial_{mm}^2 f(m,x,x^\prime):=\partial^2_{xx^\prime}\mathcal{D}^2_{mm}f(m,x,x^\prime).$$

\begin{proposition}\label{prop2.3}
Assume that the mapping $m\mapsto \mathcal{D}^2_{mm}f(m,\cdot,\cdot)$ is continuous with a modulus $\rho_f$. Then the following result holds true:
\begin{equation}\label{SOD}
\begin{split}
\Bigr|& f(m^\prime)-f(m)-\int_{\mathcal{K}}\mathcal{D}_mf(m,x)(m^\prime(dx)-m(dx))\\
&-\frac{1}{2}\int_{\mathcal{K}}\int_{\mathcal{K}}\mathcal{D}^2_{mm}f(m,x,x^\prime)(m^\prime(dx)-m(dx))
(m^\prime(dx^\prime)-m(dx^\prime))\Bigr|\\
&\leq W_2^2(m,m^\prime)\rho_f(W_2(m,m^\prime)),
\end{split}
\end{equation}
where $\rho_f(r)\rightarrow0$ as $r\rightarrow0$.
\end{proposition}
We shall  introduce some notational convention concerning copies of probability spaces and random variables. We denote by $(\widetilde{\Omega},\widetilde{\mathcal{F}},\widetilde{P})$  a copy of the probability space $(\Omega,\mathcal{F},P)$. Given a random variable or a stochastic process $W$ on $(\Omega,\mathcal{F},P)$,  we denote by $\widetilde{W}$ an independent copy of $W$ but defined on $(\widetilde{\Omega},\widetilde{\mathcal{F}},\widetilde{P})$, i.e., $\widetilde{W}$ under $\widetilde{P}$ has the same law as $W$ under $P$. For any pair of random variables $(\vartheta,\eta)$ on $(\Omega,\mathcal{F},P)$, we denote by $(\widetilde{\vartheta},\widetilde{\eta})$ an independent copy defined on $(\widetilde{\Omega},\widetilde{\mathcal{F}},\widetilde{P})$. The expectation $\widetilde{E}$ with respect to $\widetilde{P}$ concerns only the copied random variables defined over $(\widetilde{\Omega},\widetilde{\mathcal{F}},\widetilde{P})$.
\begin{proposition}\label{prop2.4}
Let $f$ be twice differentiable and satisfy (\ref{SOD}). Then, given any $\xi\in L^2(\Omega,\mathcal{F},P;\mathcal{K})$, we have
\begin{equation*}
\begin{split}
f(P_{\xi+\eta})-f(P_\xi)
=&E[\langle\partial_\mu f(P_\xi,\xi), \eta\rangle_{\mathcal{K}\times\mathcal{K}^\prime}]
 +\frac{1}{2}E[\widetilde{E}[tr(\partial_{\mu\mu}^2f(P_\xi,\xi,\widetilde{\xi})\cdot\eta\otimes\widetilde{\eta})]]\\
 &+\frac{1}{2}E[tr(\partial_x\partial_\mu f(P_\xi,\xi)\cdot\eta\otimes\eta)]+R(\xi,\eta),\
 \eta\in L^2(\Omega,\mathcal{F},P;\mathcal{K}),
\end{split}
\end{equation*}
where $R(\xi,\eta)$ satisfies the following inequality
\begin{equation*}
|R(\xi,\eta)|\leq C\( (E[|\eta|^2])^\frac{3}{2}+E[|\eta|^3] \)\leq CE[|\eta|^3].
\end{equation*}
\end{proposition}
The proof of Propositions \ref{prop2.3} and \ref{prop2.4} can refer to Proposition A.2.4 in \cite{CDLL2019} and to Lemma 2.1  in \cite{BLPR2017}, respectively.

As we have already mentioned, for our studies we are mainly interested in the case $\mathcal{K}=\mathbb{R}^k$. And so, if $f : \mathcal{P}_2(\mathbb{R}^k)\mapsto\mathbb{R}^k$ is differentiable on $\mathcal{P}_2(\mathbb{R}^k)$, we can define for the mapping $\partial_\mu f=((\partial_\mu f)_1,\cdots,(\partial_\mu f)_k)^T:\mathcal{P}_2(\mathbb{R}^k)\times\mathbb{R}^k\mapsto\mathbb{R}^k$  the second order derivatives of $f$. Similar to the first order derivative, we define\ $\partial_y\partial_\mu f(P_{\vartheta},y):\mathcal{P}_2(\mathbb{R}^k)\times\mathbb{R}^k\mapsto\mathbb{R}^k\otimes\mathbb{R}^k$\ and\ $\partial^2_{\mu\mu} f(P_{\vartheta},y,z):=\partial_\mu[\partial_\mu f(\cdot,y)](P_{\vartheta},z):\mathcal{P}_2(\mathbb{R}^k)\times\mathbb{R}^k\times\mathbb{R}^k\mapsto\mathbb{R}^k\otimes\mathbb{R}^k.$
The following both definitions are used frequently. For more details the reader may refer to \cite{BLPR2017}.
\begin{definition}\label{def3}
We say that $f$ belongs to $C_b^{1,1}(\mathcal{P}_2(\mathbb{R}^k))$, if $f: \mathcal{P}_2(\mathbb{R}^k)\rightarrow\mathbb{R}$ is differentiable on $\mathcal{P}_2(\mathbb{R}^k)$ and $\partial_\mu f(\cdot,\cdot):\mathcal{P}_2(\mathbb{R}^k)\times\mathbb{R}^k\rightarrow\mathbb{R}^k$ is bounded and Lipschitz continuous, i.e., there is some positive constant $C$ such that
\begin{itemize}
\item[\rm{(i)}] $|\partial_\mu f(\mu,y)|\leq C,\ \mu\in\mathcal{P}_2(\mathbb{R}^k),\ y\in\mathbb{R}^k;$
\item[\rm{(ii)}] $|\partial_\mu f(\mu,y)-\partial_\mu f(\mu^\prime,y^\prime)|\leq C(W_2(\mu,\mu^\prime)+|y-y^\prime|),\ \mu,\ \mu^\prime\in\mathcal{P}_2(\mathbb{R}^k),\ y,\ y^\prime\in\mathbb{R}^k.$
\end{itemize}
\end{definition}
\begin{definition}\label{def4}
We say that $f\in C_b^{2,1}(\mathcal{P}_2(\mathbb{R}^k))$, if $f\in C_b^{1,1}(\mathcal{P}_2(\mathbb{R}^k))$ and
\begin{itemize}
\item[\rm{(i)}] $(\partial_\mu f)_j(\cdot,y)\in C_b^{1,1}(\mathcal{P}_2(\mathbb{R}^k))$, for all $y\in\mathbb{R}^k,\ 1\leq j\leq k$, and $\partial_\mu^2 f:\mathcal{P}_2(\mathbb{R}^k)\times\mathbb{R}^k\times\mathbb{R}^k\rightarrow\mathbb{R}^k\otimes\mathbb{R}^k$ are bounded and Lipschitz-continuous;
\item[\rm{(ii)}] $\partial_\mu f(\mu,\cdot):\mathbb{R}^k\rightarrow\mathbb{R}^k$ is differentiable, for every $\mu\in\mathcal{P}_2(\mathbb{R}^k)$, and its derivative $\partial_y\partial_\mu f:\mathcal{P}_2(\mathbb{R}^k)\times\mathbb{R}^k\rightarrow\mathbb{R}^k\otimes\mathbb{R}^k$ is bounded and Lipschitz-continuous.
\end{itemize}
\end{definition}

\subsection{Problem formulation}
Now we consider the following mean-field forward-backward stochastic control system:
\begin{equation}\label{eq2.5}
\left\{
\begin{aligned}
&dX_t^u=b(t,X_t^u,P_{X_t^u},u_t)dt+\sigma(t,X_t^u,P_{X_t^u},u_t)dB_t,\ 0\leq t\leq T,\\
&dY_t^u=- f(t,X_t^u,Y_t^u,Z_t^u,P_{(X_t^u, Y_t^u)},u_t)dt+Z_t^udB_t,\ 0\leq t\leq T,\\
&X_0^u=x_0\in\mathbb{R}^k,\ Y_T^u=\Phi(X_T^u,P_{X_T^u}),
\end{aligned}
\right.
\end{equation}
where $(b,\sigma):[0,T]\times\mathbb{R}^k\times\mathcal{P}_2(\mathbb{R}^k)\times U\rightarrow\mathbb{R}^k\times \mathbb{R}^{k\times d}$, $f:[0,T]\times\mathbb{R}^k\times\mathbb{R}\times\mathbb{R}^d\times\mathcal{P}_2(\mathbb{R}^{k+1})
\times U\rightarrow\mathbb{R}$, $\Phi:\mathbb{R}^k\times\mathcal{P}_2(\mathbb{R}^k)\rightarrow\mathbb{R}$, and the control domain $U$ is a nonempty subset of $\mathbb{R}^n$. Notice that the control set $U$ is not necessary convex.

An admissible control $u(\cdot)$ is an $\mathbb{F}$-adapted process with values in $U$ such that
\begin{equation}\label{eq2.6}
E\[\int_0^T|u_t|^pdt\]<\infty,\  \text{for all}\ p>1.
\end{equation}
We denote the set of all admissible controls by $\mathcal{U}_{ad}$.

The assumptions on the coefficients of the control system (\ref{eq2.5}) are given below:
\begin{itemize}
\item [\textbf{(A1)}] $\Phi$ and the first-order and the second-order derivatives of $\Phi$ are continuous and bounded.
\item [\textbf{(A2)}] For all $u\in U$, and for $\psi=b$, $\sigma$, $\psi$, the first-order and the second-order derivatives of $\psi$ with respect to $(x,\mu)$ are continuous in $(x,\mu,u)$;  all the first-order and the second-order derivatives of $\psi$ are bounded; there exists a constant $L>0$ such that
$$|b(t,0,\delta_\mathbf{0},u)|+|\sigma(t,0,\delta_\mathbf{0},u)|\leq L(1+|u|).$$
\item [\textbf{(A3)}] For all $u\in U,\ (t,x,y,z)\in [0,T]\times\mathbb{R}^k\times\mathbb{R}\times\mathbb{R}^d$, the mapping $\mathcal{P}_2(\mathbb{R}^{k+1})\ni \mu\mapsto f(t,x,y,z,\mu,u)$ is in $C_b^{2,1}(\mathcal{P}_2(\mathbb{R}^{k+1}))$ with derivatives bounded over $[0,T]\times\mathbb{R}^k\times\mathbb{R}\times\mathbb{R}^d\times\mathcal{P}_2(\mathbb{R}^{k+1})\times U$.
\item [\textbf{(A4)}] The first-order derivatives of $f$ with respect to $(x,y,\mu)$ are Lipschitz continuous and bounded, and also the second-order derivatives of $f$ with respect to $(x,y,z)$ are continuous and bounded. Moreover, there exist positive constants $L_1$, $\gamma$ and $L$ such that, for all $u,\ u_1,\ u_2 \in U,\ (x,y,z),\ (\bar{x},\bar{y},\bar{z})\in \mathbb{R}^k\times\mathbb{R}\times\mathbb{R}^d$ and $\mu,\ \bar{\mu}\in\mathcal{P}_2(\mathbb{R}^{k+1})$,
\begin{equation*}
\begin{split}
&|f(t,x,y,z,\mu,u)-f(t,x,\bar{y},\bar{z},\bar{\mu},u)|\leq L(|y-\bar{y}|+W_2(\mu,\bar{\mu})+(1+|z|+|\bar{z}|)|z-\bar{z}|),\\
&|f(t,x,y,z,\mu,u_1)-f(t,x,y,z,\mu,u_2)|\leq L(1+|x|+|y|+|z|),\\
&|f_z(t,x,y,z,\mu,u)|\leq L+\gamma|z|,\\
&|f(t,x,0,0,\delta_\mathbf{0},u)|\leq L_1.
\end{split}
\end{equation*}
\end{itemize}
From the existence and the uniqueness result of Theorem 3.5 in \cite{HHTW2022} (see also Theorem 4.1 in \cite{JLW2023}), we have the following result.
\begin{proposition}\label{thm2.1}
Let the assumptions \textbf{(A1)}-\textbf{(A4)} hold. Then for any $u(\cdot)\in\mathcal{U}_{ad}$ and $p>1$, the control system (\ref{eq2.5}) admits a unique solution $(X^u,Y^u,Z^u)\in\mathcal{S}^p(0,T;\mathbb{R}^k)\times \mathcal{S}^\infty(0,T;\mathbb{R})\times \mathcal{H}^{2,p}(0,T;\mathbb{R}^d)$ such that $Z^u\cdot B=\int_0^{\cdot}Z_s^udB_s\in BMO_2$. Furthermore, we have the following estimates:
\begin{flalign}\label{eq2.7}
&\ {\rm{(i)}}\ \|X^u\|_{\mathcal{S}^p}^p\leq C\Bigr\{|X_0|^p+E\[\(\int_0^T|b(t,0,\delta_{\mathbf{0}},u_t)|dt\)^p+\(\int_0^T|\sigma(t,0,\delta_{\mathbf{0}},u_t)|^2dt\)^\frac{p}{2}\]\Bigr\},&
\end{flalign}
\begin{flalign}\label{eq2.8}
&\ {\rm{(ii)}}\ \|Y^u\|_{\infty}+\|Z^u\cdot B\|_{BMO_2}\leq C_1,&
\end{flalign}
where $C$ depends on $p,\ T,\ \|b_x\|_\infty,\ \|b_\mu\|_\infty,\ \|\sigma_x\|_\infty,\ \|\sigma_\mu\|_\infty$, and $C_1$ depends on $K,\ \|\Phi\|_\infty,\ T$.
\end{proposition}
\begin{remark}\label{re2.1}
From (\ref{eq2.8}) and Lemma \ref{lem2.1} {\rm (i)}, for all $n\geq1,\ u\in\mathcal{U}_{ad}$,
\begin{equation*}
E\[\(\int_0^T|Z_s^u|^2ds\)^n\]\leq C_n E\[\langle Z^u\cdot B\rangle_T^n\]\leq C^1_n \|Z^u\cdot B\|^{2n}_{BMO_2}\leq C_n^u.
\end{equation*}
\end{remark}
Our recursive optimal control problem is to minimize the cost functional $J(u):=Y_0^u$ , which is the solution of backward stochastic differential equation in the control system (\ref{eq2.5}) at $t=0$:
\begin{equation}\label{eq2.9}
\inf\limits_{u\in \mathcal{U}_{ad}}J(u).
\end{equation}
Note that the control $u^{\ast}\in \mathcal{U}_{ad}$ is called an optimal control, if it satisfies
\begin{equation}\label{eq2.10}
J(u^{\ast})=\inf\limits_{u\in \mathcal{U}_{ad}}J(u).
\end{equation}
\section{Mean-field BSDEs with Stochastic Lipschitz Conditions}\label{MFBSDE}
Before studying the optimal control problem, we first consider the existence and the uniqueness of the mean-field BSDE with stochastic Lipschitz condition and its moment estimates. For simplicity, we consider the one-dimensional case, i.e., $k=d=1$.  In this section, we consider the following linear mean-field BSDE:
\begin{equation}\label{eq4.1}
\left\{
\begin{aligned}
dY_t&=-\(\varrho_t Y_t+\varpi_tZ_t+  \widetilde{E}[\varrho_t^\prime\widetilde{Y}_t]+\varphi_t\)dt+Z_tdB_t \ ,\ \ t\in[0,T],\\       
Y_T&=\xi.
\end{aligned}
\right.
\end{equation}
The assumptions of this section are the following: For $p>1$,
\begin{itemize}
\item [\textbf{(B1)}] The process $\varrho$ and $\varpi$ are $\mathbb{F}$-adapted. $\varrho$ is bounded and  $\varpi\cdot B =\int^{\cdot}_0\varpi_t d B_t$ is a $BMO$ martingale. The process  $\varrho^\prime:\ [0,T]\times\Omega\times\widetilde{\Omega}\rightarrow\mathbb{R}$ is defined over the tensor product between $(\Omega,\mathcal{F},P)$ and a copy $(\widetilde{\Omega},\widetilde{\mathcal{F}},\widetilde{P})$, bounded and adapted to the completed filtration generated by $B$ and its copy $\widetilde{B}$ over $(\widetilde{\Omega},\widetilde{\mathcal{F}},\widetilde{P})$.
\item [\textbf{(B2)}] $\xi\in \cap_{p>1}L^p(\mathcal{F}_T;\mathbb{R})$ and $\varphi=(\varphi_t)\in \cap_{p>1}\mathcal{H}^{1,p}(0,T;\mathbb{R})$.
\end{itemize}
\begin{proposition}\label{prop4.1}
Let the assumptions \textbf{(B1)}-\textbf{(B2)} hold, and let $p_M^*$ be the conjugate exponent of $p_M$ for $M_\cdot=\displaystyle\int_0^\cdot \varpi_tdB_t$. Then for any $\kappa_0>p^*_M$, we have
\begin{equation}\label{eq4.2}
E\[\sup\limits_{t\in[0,T]}|Y_t|^\kappa+\(\int_0^T|Z_t|^2dt\)^{\frac{\kappa}{2}}\]\leq CE\[|\xi|^{\kappa_0}+(\int_0^T|\varphi_t|dt)^{\kappa_0}\]^{\frac{\kappa}{\kappa_0}},
\end{equation}
for all $\kappa\in(1\vee p^*_M,\kappa_0)$
where $C$ depends on $\kappa_0,\ T,\ \|\varrho\|_\infty,\ \|\varrho^\prime\|_\infty,\ \|\varpi\cdot B\|_{BMO_2}$, and $(Y,Z)\in\mathcal{S}^\kappa(0,T;\mathbb{R})\times\mathcal{H}^{2,\kappa}(0,T;\mathbb{R})$ is supposed to be a solution of BSDE (\ref{eq4.1}).
\end{proposition}
\begin{proof}
We first introduce $X_t=\exp{\(\displaystyle\int_0^t\varrho_sds\)}$
and $\Gamma_t=\mathcal{E}\(\displaystyle\int_0^\cdot\varpi_sdB_s\)_t$. Applying It\^{o}'s formula to $\Gamma_tX_tY_t$,
\begin{equation}\label{eq4.4}
\Gamma_tX_tY_t=\Gamma_TX_T\xi+\int_t^T\Gamma_sX_s\(\widetilde{E}[\varrho^\prime_s\widetilde{Y}_s]+\varphi_s\)ds 
-\int_t^T\Gamma_sX_s\(\varpi_sY_s+Z_s\)dB_s,\ 0\leq t\leq T.
\end{equation}
By using the reverse H\"{o}lder inequality, we get $E\[\(\displaystyle\int_0^T|\Gamma_tX_t(\varpi_tY_t+Z_t)|^2dt\)^{\frac{1}{2}}\]<\infty$. Thus, by taking conditional expectation in (\ref{eq4.4}), and setting $(\Gamma_s^t)_{s\in[t,T]}:=(\frac{\Gamma_s}{\Gamma_t})_{s\in[t,T]}$ and $(X_s^t)_{s\in[t,T]}:=(\frac{X_s}{X_t})_{s\in[t,T]}$, we have
\begin{equation}\label{eq4.5}
Y_t=E\[\Gamma^t_TX^t_T\xi+\int_t^T\Gamma^t_sX^t_s(\widetilde{E}[\varrho^\prime_s\widetilde{Y}_s]+\varphi_s)ds|\mathcal{F}_t\]. 
\end{equation}

We use (\ref{eq4.5}) (recall that $\varrho^\prime$ is bounded) and the reverse H\"{o}lder inequality (recall that $1<\kappa_0^*<p_M$). Then, for a constant $C$ which varies from line to line, we get
\begin{equation}\label{eq4.7}
\begin{aligned}
|Y_t|\leq& CE\[\Gamma^t_T\(|\xi|+\int_0^T|\varphi_s|ds\)|\mathcal{F}_t\]+C\int_t^TE[|Y_s|]ds\\
     \leq& C\(E\[(\Gamma^t_T)^{\kappa_0^*}|\mathcal{F}_t\]\)^{\frac{1}{\kappa_0^*}}
            \(E\[\(|\xi|+\int_0^T|\varphi_s|ds\)^{\kappa_0}|\mathcal{F}_t\]\)^{\frac{1}{\kappa_0}}+C\int_t^TE[|Y_s|]ds\\
     \leq& C\(E\[\(|\xi|+\int_0^T|\varphi_s|ds\)^{\kappa_0}|\mathcal{F}_t\]\)^{\frac{1}{\kappa_0}}
            +C\int_t^TE[|Y_s|]ds,\ t\in[0,T].
\end{aligned}
\end{equation}
Hence, from (\ref{eq4.7}), due to H\"{o}lder's inequality and Doob's inequality, we have
\begin{equation}\label{eq4.8}
E\[\sup_{s\in[t,T]}|Y_s|^\kappa\]
\leq C\(E\[|\xi|^{\kappa_0}+(\int_0^T|\varphi_t|dt)^{\kappa_0}\]\)^{\frac{\kappa}{\kappa_0}}
    +C\int_t^TE[|Y_s|^\kappa]ds,\ t\in[0,T],
\end{equation}
and Gronwall's inequality yields
\begin{equation}\label{eq4.6}
E\[\sup_{s\in[0,T]}|Y_s|^\kappa\]
\leq C\(E\[|\xi|^{\kappa_0}+\(\int_0^T|\varphi_t|dt\)^{\kappa_0}\]\)^\frac{\kappa}{\kappa_0}.
\end{equation}
We define the stopping times
\begin{equation*}
\tau_n=\inf\Bigr\{t\in[0,T],\ \int_0^t|Z_s|^2ds\geq n\Bigr\}\wedge T,\ n\geq 1.
\end{equation*}
Applying It\^{o}'s formula to $|Y_t|^2$, we get
\begin{equation*}
|Y_0|^2+\int_0^{\tau_n}|Z_s|^2ds=|Y_{\tau_n}|^2+2\int_0^{\tau_n} Y_s(\varrho_sY_s+\varpi_sZ_s+\widetilde{E}[\varrho_s^\prime \widetilde{Y}_s]+\varphi_s)ds-2\int_0^{\tau_n} Y_sZ_sdB_s.
\end{equation*}
Thus,
\begin{equation*}
\begin{aligned}
\(\int_0^{\tau_n}|Z_s|^2ds\)^{\frac{\kappa}{2}}
\leq & C\(\sup_{t\in[0,T]}|Y_t|^\kappa+\sup_{t\in[0,T]}|Y_t|^\kappa(\int_0^T|\varpi_s|^2ds)^\frac{\kappa}{2}+(\int_0^T |\varphi_s|ds)^\kappa+(\int_0^T E[|Y_s|]ds)^\kappa\)\\
&+C\Bigr|\int_0^{\tau_n} Y_sZ_sdB_s\Bigr|^{\frac{\kappa}{2}}
    +\frac{1}{2}\(\int_0^{\tau_n}|Z_s|^2ds\)^\frac{\kappa}{2}.
\end{aligned}
\end{equation*}
Hence, using the Burkholder-Davis-Gundy inequality, H\"{o}lder inequality and energy inequality, we have
\begin{equation*}
\begin{split}
E\[\(\int_0^{\tau_n}|Z_s|^2ds\)^{\frac{\kappa}{2}}\]
&\leq CE\[\sup_{t\in[0,T]}|Y_t|^\kappa+\sup_{t\in[0,T]}|Y_t|^\kappa(\int_0^T|\varpi_s|^2ds)^\frac{\kappa}{2}+(\int_0^T |\varphi_s|ds)^\kappa+(\int_0^T E[|Y_s|]ds)^\kappa\]\\
&\leq C\(E\[|\xi|^{\kappa_0}+(\int_0^T|\varphi_t|dt)^{\kappa_0}\]\)^\frac{\kappa}{\kappa_0}.
\end{split}
\end{equation*}
Indeed, by choosing $\kappa^\prime>1$ such that $\kappa^\prime\kappa<\kappa_0$, we obtain
\begin{equation*}
\begin{aligned}
E\[\sup_{t\in[0,T]}|Y_t|^\kappa(\int_0^T|\varpi_s|^2ds)^\frac{\kappa}{2}\]
\leq& \(E\[\sup_{t\in[0,T]}|Y_t|^{\kappa^\prime\kappa}\]\)^\frac{1}{\kappa^\prime}
    \(E\[(\int_0^T|\varpi_s|^2ds)^\frac{(\kappa^\prime)^*\kappa}{2}\]\)^\frac{1}{(\kappa^\prime)^*}\\
\leq &C\|\varpi\cdot B\|_{BMO_2}^\kappa\(E\[|\xi|^{\kappa_0}+(\int_0^T|\varphi_t|dt)^{\kappa_0}\]\)^\frac{\kappa}{\kappa_0},
\end{aligned}
\end{equation*}
which we combine with (\ref{eq4.6}). To complete the proof, it suffices to let $n\rightarrow\infty$.
\end{proof}


We recall that the classical BSDE (without mean-field terms) with stochastic Lipschitz condition has a unique solution, see \cite{BC2008}, Theorem 10.
\begin{lemma}[\textbf{\cite{BC2008},\ Theorem 10}]
Let the assumptions \textbf{(B1)-(B2)} hold with $\varrho^\prime\equiv0$. Then BSDE (\ref{eq4.1}) has a unique solution $(Y,Z)$ which belongs to $\mathcal{S}^p(0,T;\mathbb{R})\times\mathcal{H}^{2,p}(0,T;\mathbb{R})$ for all $p<\kappa_0$.
\end{lemma}
Now we extend this result to the mean-field case.
\begin{theorem}\label{thm4.1}
Suppose that \textbf{(B1)}-\textbf{(B2)} hold. Then there exists a unique solution $(Y,Z)\in\mathcal{S}^\kappa(0,T;\mathbb{R})\times\mathcal{H}^{2,\kappa}(0,T;\mathbb{R})$ of the BSDE (\ref{eq4.1}) for all $\kappa\in(1\vee p_M^*,\kappa_0)$, where $\kappa_0>p^*_M$.
\end{theorem}
\begin{proof}
We first prove the uniqueness of the solution for BSDE (\ref{eq4.1}). Let $(Y^1,Z^1)$ and $(Y^2,Z^2)$ be two solutions of BSDE (\ref{eq4.1}). Then we can deduce from Proposition \ref{prop4.1}, that $(Y^i,Z^i)\in \mathcal{S}^\kappa(0,T;\mathbb{R})\times\mathcal{H}^{2,\kappa}(0,T;\mathbb{R}),\ i=1,2$, for all $\kappa\geq p_M^*$, and the processes $\delta Y:=Y^1-Y^2$, $\delta \widetilde{Y}:=\widetilde{Y}^{1}-\widetilde{Y}^{2}$ and $\delta Z:=Z^1-Z^2$ satisfy the following BSDE:
\begin{equation*}
\delta Y_t=
\int_t^T\(\varrho_s \delta Y_s+\varpi_s\delta Z_s+\widetilde{E}[\varrho_s^\prime\delta \widetilde{Y}_s]\)ds+\int_t^T \delta Z_sdB_s,\ t\in[0,T].
\end{equation*}
We notice that the above BSDE satisfies our assumptions \textbf{(B1)}-\textbf{(B2)}, and due to Proposition \ref{prop4.1} we obtain $(\delta Y, \delta Z)\equiv(0,0)\ \text{in }\mathcal{S}^\kappa(0,T;\mathbb{R})\times\mathcal{H}^{2,\kappa}(0,T;\mathbb{R})$.

Now we prove the existence of the solution for BSDE (\ref{eq4.1}). For any $n\geq 1$, we define the stopping time $\tau_n$,
$$\tau_n=\inf\Bigr\{t\in[0,T],\ \int_0^t\(|\varphi_s|+\frac{1}{2}|\varpi_s|^2\)ds\geq n\Bigr\}\wedge T.$$
Let $\xi^n=\xi I_{\{|\xi|\leq n\}}$ and $(Y^n,Z^n)$ be the solution of the following BSDE:
\begin{equation}\label{eq4.9}
Y_t^n=\xi^n+\int_t^T I_{\{s\leq\tau_n\}}\(\varrho_s Y^n_s+\varpi_s Z^n_s+\widetilde{E}[\varrho_s^\prime \widetilde{Y}^{n}_s]+\varphi_s\)ds-\int_t^T Z_s^ndB_s.
\end{equation}
Note that
\begin{equation*}
\begin{split}
\Bigr|I_{\{s\leq\tau_n\}}\(\varrho_s y+\varpi_s z+E[\varrho_s^\prime \eta]+\varphi_s\)\Bigr|
\leq \theta_s+K|y|+K E[|\eta|]+\frac{1}{2}|z|^2,
\end{split}
\end{equation*}
for all $(s,y,z,\eta)\in [0,T]\times\mathbb{R}\times\mathbb{R}\times L^2(\mathcal{F}_T;\mathbb{R})$, where $K\in\mathbb{R}_+$ is a bound of the processes $\varrho$ and $\varrho^\prime$, and
$\theta_s=I_{\{s\leq\tau_n\}}\(|\varphi_s|+\frac{1}{2}|\varpi_s|^2\)$, $s\in[0,T]$. Note that $\displaystyle\int_0^T|\theta_s|ds\leq n$, $P$-a.s.

Similar to the proof of Theorem 2.5 in \cite{M2009}, but only with a slight extension, we see that the BSDE (\ref{eq4.9}) has a unique solution $(Y^n,Z^n)$, and $Y^n$ is a bounded process and $Z\cdot B\in BMO_2$.
Moreover, from Proposition \ref{prop4.1} we have (\ref{eq4.2}) for $(Y^n,Z^n)\in\mathcal{S}^\kappa(0,T;\mathbb{R})\times\mathcal{H}^{2,\kappa}(0,T;\mathbb{R})$.

We still need to show $(Y^n,Z^n)_{n\geq 1}$ is a Cauchy sequence.
Let $\underline{Y}^{m,n}=Y^m-Y^n$ and $\underline{Z}^{m,n}=Z^m-Z^n$. We remark that, for $m\geq n\geq 1$,
\begin{equation*}
\left\{
\begin{aligned}
d\underline{Y}^{m,n}_t&=-\(I_{(0,\tau_n]}(t)\varrho_t \underline{Y}^{m,n}_t
        +I_{(0,\tau_n]}(t)\varpi_t \underline{Z}^{m,n}_t+I_{(0,\tau_n]}(t)\widetilde{E}[\varrho^\prime_t \widetilde{\underline{Y}}^{m,n}_t]+\varphi_t^{m,n}\)dt+\underline{Z}_t^{m,n}dB_t,\\
\underline{Y}^{m,n}_T&=\xi I_{\{n<|\xi|\leq m\}},
\end{aligned}
\right.
\end{equation*}
where
$$\varphi_t^{m,n}=I_{(\tau_n,\tau_m]}(t)\(\varrho_t Y^m_t+\varpi_t Z^m_t+\widetilde{E}[\varrho_t^\prime \widetilde{Y}^{m}_t]+\varphi_t\),\ t\in[0,T].$$
Then, as the BSDE for $(\underline{Y}^{m,n},\underline{Z}^{m,n})$ satisfies the assumptions of Proposition \ref{prop4.1},
\begin{equation}\label{eq4.11}
E\[\sup_{t\in[0,T]}|\underline{Y}^{m,n}_t|^\kappa+\(\int_0^T|\underline{Z}_t^{m,n}|^2dt\)^{\frac{\kappa}{2}}\]
\leq C\(E\[|\xi|^{\kappa^\prime} I_{\{n<|\xi|\leq m\}}+\(\int_0^T|\varphi_t^{m,n}|dt\)^{\kappa^\prime}\]\)^\frac{\kappa}{\kappa^\prime},
\end{equation}
for $1<p_M^*<\kappa<\kappa^\prime<\kappa_0$ and $C\in\mathbb{R}_+$ depending on $\kappa,\ \kappa^\prime$.
As $\varrho,\ \varrho^\prime$ are bounded and $\varpi\cdot B$ is BMO, we have, for $1<p_M^*<\kappa<\kappa^{\prime}<\kappa^{\prime\prime}<\kappa_0$,
\begin{equation}\label{eq4.12}
\begin{aligned}
&E\[\(\int_0^T|\varphi_t^{m,n}|dt\)^{\kappa^{\prime}}\]\\
\leq& CE\[\(\int_{\tau_n}^{\tau_m}(|Y_t^m|+E[|Y_t^{m}|]
      +|\varphi_t|+|\varpi_tZ_t^m|)dt\)^{\kappa^{\prime}}\]\\
\leq&C\(E[\sup_{t\in[0,T]}|Y_t^m|^{\kappa^{\prime\prime}}]\)^{\frac{\kappa^\prime}{\kappa^{\prime\prime}}}
        \(E[(\tau_m-\tau_n)^{\overline{\kappa}}]\)^{\frac{\kappa^{\prime}}{\overline{\kappa}}}
+CE\[\(\int_{\tau_n}^{\tau_m}|\varphi_t|dt\)^{\kappa^{\prime}}\]\\
&+C\(E\[\(\int_{\tau_n}^{\tau_m}|\varpi_t|^2dt\)^\frac{\overline{\kappa}}{2}\]\)^{\frac{\kappa^\prime}{\overline{\kappa}}}
    \(E\[\(\int_{0}^{T}|Z^m_t|^2dt\)^\frac{\kappa^{\prime\prime}}{2}\]\)^{\frac{\kappa^{\prime}}{\kappa^{\prime\prime}}},\ m\geq n\geq 1,
\end{aligned}
\end{equation}
for $\overline{\kappa}=(\kappa^{\prime}\kappa^{\prime\prime})/(\kappa^{\prime\prime}-\kappa^{\prime})$.

Consequently, as due to Proposition \ref{prop4.1}, applied to the BSDE for $(Y^m,Z^m)$,
\begin{equation*}
E\[\sup_{t\in[0,T]}|Y_t^{m}|^{\kappa^{\prime\prime}}+\(\int_0^T|Z_t^m|^2dt\)^{\frac{\kappa^{\prime\prime}}{2}}\]
\leq C\(E\[|\xi|^{\kappa_0}+\(\int_0^T|\varphi_t|dt\)^{\kappa_0}\]\)^\frac{\kappa^{\prime\prime}}{\kappa_0}
=C_{\kappa^{\prime\prime}},
\end{equation*}
it follows from (\ref{eq4.11}) and (\ref{eq4.12}) that, for all $p_M^*<\kappa<\kappa^\prime<\kappa^{\prime\prime}<\kappa_0$, $m\geq n\geq 1$,
\begin{equation}\label{eq4.13}
\begin{aligned}
&E\[\sup_{t\in[0,T]}|\underline{Y}_t^{m,n}|^\kappa+\(\int_0^T|\underline{Z}_t^{m,n}|^2dt\)^{\frac{\kappa}{2}}\]\\
\leq &C\(E\[|\xi|^{\kappa^\prime} I_{\{|\xi|\geq n\}}
        +\(\int_{\tau_n}^T|\varphi_t|dt\)^{\kappa^\prime}\]\)^\frac{\kappa}{\kappa^\prime}
+C_{\kappa^\prime,\kappa^{\prime\prime}}\(\(E[(T-\tau_n)^{\overline{\kappa}}]\)^\frac{\kappa}{\overline{\kappa}}
+\(E\[\(\int_{\tau_n}^T|\varpi_t|^2dt\)^\frac{\overline{\kappa}}{2}\]\)^\frac{\kappa}{\overline{\kappa}}\).
\end{aligned}
\end{equation}
However,
$$E\[\(\int_0^T|\varpi_t|^2dt\)^\frac{\overline{\kappa}}{2}\]\leq
C_{\overline{\kappa}}\|\varpi\cdot B\|^{\overline{\kappa}}_{BMO_2}<\infty.$$
Thus, thanks to the dominated convergence theorem we see that the right-hand side of (\ref{eq4.13}) converges to zero, as $n\rightarrow\infty$. Consequently, $(Y^n,Z^n)$ is a Cauchy sequence in $\mathcal{S}^\kappa(0,T;\mathbb{R})\times\mathcal{H}^{2,\kappa}(0,T;\mathbb{R}),\ \kappa\in(1\vee p_M^*,\kappa_0)$, and passing to the limit in (\ref{eq4.9}) yields the existence result.
\end{proof}


\section{Variational Equations}
In this section we  study the first-order and  the second-order variational equations. Let $u^{\ast}\in \mathcal{U}_{ad}$ be the optimal control and $(X^{\ast},Y^{\ast},Z^{\ast})$ be the corresponding state trajectories of $(\ref{eq2.5})$. Recall that the control set is not necessary convex, hence we shall consider the following ``spike variation" of $u^{\ast}$. For any $\varepsilon>0$, we choose a Borel subset $E_\varepsilon\subset [0,T]$ with Borel measure $|E_\varepsilon|\leq\varepsilon$. Given an arbitrary $u\in\mathcal{U}_{ad}$, we define, for all $t\in[0,T]$,
\begin{equation}\label{eq3.1}
u_t^\varepsilon=
\left\{
\begin{aligned}
u_t^\ast,\ \ \ & t\in[0,T]\setminus E_\varepsilon,\\
u_t,\ \ \  & t\in E_\varepsilon.
\end{aligned}
\right.
\end{equation}
Observe that $u^\varepsilon\in\mathcal{U}_{ad}$, and denote by  $(X^{\varepsilon},Y^{\varepsilon},Z^{\varepsilon})$  the solution of the control system (\ref{eq2.5}) with the admissible control $u^\varepsilon$. For simplicity, we assume $k=d=n=1$. However, our result can be easily extended to the multidimensional case.
For notational convenience, for $\psi=b,\sigma,f,\Phi$ and $w=x,y,z$, we put
\begin{equation*}
\begin{split}
&\psi(t):=\psi(t,X^{\ast}_t,Y^{\ast}_t,Z^{\ast}_t,P_{(X_t^{\ast},Y_t^{\ast})},u_t^{\ast});\ \  \psi_w(t):=\psi_w(t,X_t^{\ast},Y_t^{\ast},Z_t^{\ast},P_{(X_t^{\ast},Y_t^{\ast})},u_t^{\ast});\\
&\delta\psi(t):=\psi(t,X_t^{\ast},Y_t^{\ast},Z_t^{\ast},P_{(X_t^{\ast},Y_t^{\ast})},u_t^{\varepsilon})-\psi(t);\  \delta\psi_w(t):=\psi_w(t,X^{\ast}_t,Y^{\ast}_t,Z^{\ast}_t,P_{(X^{\ast}_t,Y^{\ast}_t)},u^{\varepsilon}_t)-\psi_w(t).
\end{split}
\end{equation*}
We recall the notational convention concerning copies of probability spaces and random variables, we have introduced at the end of Section \ref{derivative}. And so we define
\begin{equation*}
\begin{split}
&\widetilde{\psi}_{\mu}(t):=(\partial_{\mu}\psi)(t,X^{\ast}_t,Y^{\ast}_t,Z^{\ast}_t,P_{(X_t^{\ast},Y_t^{\ast})},u_t^{\ast};(\widetilde{X}^{\ast}_t,\widetilde{Y}^{\ast}_t)),\\ &\widetilde{\psi}_{\mu}^\ast(t):=(\partial_{\mu}\psi)(t,\widetilde{X}^{\ast}_t,\widetilde{Y}^{\ast}_t,\widetilde{Z}^{\ast}_t,P_{(X_t^{\ast},Y_t^{\ast})},\widetilde{u}_t^{\ast};(X^{\ast}_t,Y^{\ast}_t)),\\
\end{split}
\end{equation*}
$\widetilde{\psi}_{\mu,x}(t)$ represents the first component of $\widetilde{\psi}_{\mu}(t)$ and, likewise, $\widetilde{\psi}_{\mu,y}(t)$ represents the second component of $\widetilde{\psi}_{\mu}(t)$.
Next, we introduce our notations for the second-order derivatives of $\psi$ :
\begin{equation*}
\begin{split}
\psi_{ww}(t)&:=\psi_{ww}(t,X_t^{\ast},Y_t^{\ast},Z_t^{\ast},P_{(X_t^{\ast},Y_t^{\ast})},u_t^{\ast}),\\
\widetilde{\psi}_{w\mu}(t)&:=\partial_w(\partial_{\mu}\psi)(t,X^{\ast}_t,Y^{\ast}_t,Z^{\ast}_t,P_{(X_t^{\ast},Y_t^{\ast})},u_t^{\ast};(\widetilde{X}^{\ast}_t,\widetilde{Y}^{\ast}_t)),\\
\widetilde{\psi}_{\widetilde{x}\mu}(t)&:=\partial_{\widetilde{x}}(\partial_{\mu}\psi)(t,X^{\ast}_t,Y^{\ast}_t,Z^{\ast}_t,P_{(X_t^{\ast},Y_t^{\ast})},u_t^{\ast};(\widetilde{X}^{\ast}_t,\widetilde{Y}^{\ast}_t)),\\
\widetilde{\psi}_{\widetilde{y}\mu}(t)&:=\partial_{\widetilde{y}}(\partial_{\mu}\psi)(t,X^{\ast}_t,Y^{\ast}_t,Z^{\ast}_t,P_{(X_t^{\ast},Y_t^{\ast})},u_t^{\ast};(\widetilde{X}^{\ast}_t,\widetilde{Y}^{\ast}_t)),\\
\widetilde{\psi}^\ast_{\widetilde{x}\mu}(t)&:=\partial_{\widetilde{x}}(\partial_{\mu}\psi)(t,\widetilde{X}^{\ast}_t,\widetilde{Y}^{\ast}_t,\widetilde{Z}^{\ast}_t,P_{(X_t^{\ast},Y_t^{\ast})},\widetilde{u}_t^{\ast};(X^{\ast}_t,Y^{\ast}_t)),\\
\widetilde{\psi}^\ast_{\widetilde{y}\mu}(t)&:=\partial_{\widetilde{y}}(\partial_{\mu}\psi)(t,\widetilde{X}^{\ast}_t,\widetilde{Y}^{\ast}_t,\widetilde{Z}^{\ast}_t,P_{(X_t^{\ast},Y_t^{\ast})},\widetilde{u}_t^{\ast};(X^{\ast}_t,Y^{\ast}_t)).
\end{split}
\end{equation*}
The first-order and the second-order variational equations for the stochastic differential equation in (\ref{eq2.5}) are
\begin{equation}\label{eq3.2}
\left\{
\begin{aligned}
dX_t^1&=\(b_x(t)X_t^1+\widetilde{E}[\widetilde{b}_{\mu}(t) \widetilde{X}_t^1]\)dt+\(\sigma_x(t)X_t^1+\widetilde{E}[\widetilde{\sigma}_{\mu}(t) \widetilde{X}_t^1]+\delta\sigma(t)I_{E_\varepsilon}(t)\)dB_t,\\
X_0^1&=0,
\end{aligned}
\right.
\end{equation}
and
\begin{equation}\label{eq3.3}
\left\{
\begin{aligned}
dX_t^2=&\(b_x(t)X_t^2+\widetilde{E}[\widetilde{b}_{\mu}(t) \widetilde{X}_t^2]+\frac{1}{2}b_{xx}(t)(X_t^1)^2+\frac{1}{2}\widetilde{E}[\widetilde{b}_{\widetilde{x}\mu}(t)(\widetilde{X}_t^1)^2]+\delta b(t)I_{E_\varepsilon}(t)\)dt\\
&+\(\sigma_x(t)X_t^2+\widetilde{E}[\widetilde{\sigma}_{\mu}(t) \widetilde{X}_t^2]+\frac{1}{2}\sigma_{xx}(t)(X_t^1)^2+\frac{1}{2}\widetilde{E}[\widetilde{\sigma}_{\widetilde{x}\mu}(t)(\widetilde{X}_t^1)^2]+\delta\sigma_x(t)X_t^1I_{E_\varepsilon}(t)\)dB_t,\\
X_0^2=&0.
\end{aligned}
\right.
\end{equation}
We can easily obtain the following results.
\begin{proposition}\label{prop3.1}
Suppose that \textbf{(A2)} holds true. For every $p>1$, (\ref{eq3.2}) and (\ref{eq3.3}) have a unique solution $X^1$ and $X^2\in\mathcal{S}^p(0,T;\mathbb{R})$, respectively. Moreover,
\begin{equation}\label{eq3.4}
E\[\sup_{t\in[0,T]}|X_t^\varepsilon-X_t^\ast|^p\]= O(\varepsilon^{\frac{p}{2}}),
\end{equation}
\begin{equation}\label{eq3.5}
E\[\sup_{t\in[0,T]}|X_t^1|^p\]= O(\varepsilon^{\frac{p}{2}}),
\end{equation}
\begin{equation}\label{eq3.6}
E\[\sup_{t\in[0,T]}|X_t^2|^p\]= O(\varepsilon^{p}),
\end{equation}
\begin{equation}\label{eq3.7}
E\[\sup_{t\in[0,T]}|X_t^\varepsilon-X_t^\ast-X_t^1|^p\]= O(\varepsilon^{p}),
\end{equation}
\begin{equation}\label{eq3.8}
E\[\sup_{t\in[0,T]}|X_t^\varepsilon-X_t^\ast-X_t^1-X_t^2|^p\]= o(\varepsilon^{p}),
\end{equation}
as $\varepsilon\downarrow0$.
Moreover, we have the following estimate:
\begin{equation}\label{eq3.18}
\begin{aligned}
E\[&\Bigr|\Phi(X_T^\varepsilon,P_{X_T^\varepsilon})-\Phi(X_T^*,P_{X_T^*})-\Phi_x(T)(X_T^1+X_T^2)
-\widetilde{E}\[\widetilde{\Phi}_{\mu}(T)(\widetilde{X}_T^1+\widetilde{X}_T^2)\]-\frac{1}{2}\Phi_{xx}(T)(X_T^1)^2\\
&-\frac{1}{2}\widetilde{E}\[\widetilde{\Phi}_{\widetilde{x}\mu}(T)(\widetilde{X}_T^1)^2\]\Bigr|^p\]
=o(\varepsilon^p),\ \varepsilon\downarrow0.
\end{aligned}
\end{equation}
\end{proposition}
By the boundedness assumptions on the derivatives of $b$ and $\sigma$, we can easily see that for the linear mean-field  SDEs (\ref{eq3.2}) and (\ref{eq3.3}) there exists a unique solution. As concerns the  estimates (\ref{eq3.4})-(\ref{eq3.18}), we can refer to the Propositions 4.2 and 5.2 as well as Lemma 5.2 in \cite{BLM2016}.

We also  need the following  extension of Proposition 5.2 in \cite{BCL2021}.
\begin{proposition}\label{lem4.1}
Let $\theta\in\mathcal{H}^2(0,T)$. Then there exists a function $\rho:\mathbb{R}_+\rightarrow\mathbb{R}_+$ such that $\rho(\varepsilon)\rightarrow 0$ as $\varepsilon\rightarrow0$, and
$$\int_0^T|E[\theta_tX_t^1]|^2dt\leq\rho(\varepsilon)\varepsilon,\ \varepsilon>0.$$
\end{proposition}

\begin{remark}
Unlike Proposition \ref{lem4.1} which only requires that $\theta\in \mathcal{H}^2(0,T)$, Proposition 5.2 in \cite{BCL2021} supposes in addition that, for some $C\in\mathbb{R}_+$ real constant
$$E[|\theta_t|^2]\leq C,\ t\in[0,T].$$
However, this additional assumption does not apply to the estimates we have to do in the sequel.
\end{remark}
For convenience of the reader and a better readability of the manuscript we give the proof, but it will be shifted to the Appendix.

In order to obtain the variational equations associated with the BSDE in (\ref{eq2.5}), we consider the following both adjoint equations:
\begin{equation}\label{eq3.9}
\left\{
\begin{aligned}
dp(t)&=-F(t)dt+q_1(t)dB_t+q_2(t)d\widetilde{B}_t,\ t\in[0,T],\\
p(T)&=(p_1(T),p_2(T))=(\Phi_x(T),\widetilde{\Phi}_{\mu}(T)),
\end{aligned}
\right.
\end{equation}
\begin{equation}\label{eq3.10}
\left\{
\begin{aligned}
dP(t)&=-G(t)dt+Q_1(t)dB_t+Q_2(t)d\widetilde{B}_t,\ t\in[0,T],\\
P(T)&=(P_1(T),P_2(T))=(\Phi_{xx}(T),\widetilde{\Phi}_{\widetilde{x}\mu}(T)),
\end{aligned}
\right.
\end{equation}
where $F(t)=(F_1(t),F_2(t))$ and $G(t)=(G_1(t),G_2(t))$ are processes adapted to the completed filtration $\mathbb{F}^{B,\widetilde{B}}$ generated by $(B,\widetilde{B})$; $F$ and $G$ are defined further down. They guarantee that (\ref{eq3.9}) and (\ref{eq3.10}) have a unique solution in $\mathcal{S}^p([0,T];\mathbb{R}^2)\times\mathcal{H}^{2,p}([0,T];\mathbb{R}^4)$, $p>1$.
The equations (\ref{eq3.9}) and (\ref{eq3.10}) can be split as follows. For $p(t)=(p_1(t),p_2(t))$, we have
\begin{equation}\label{eq3.11}
\left\{
\begin{aligned}
dp_1(t)&=-F_1(t)dt+q_{1,1}(t)dB_t+q_{1,2}(t)d\widetilde{B}_t,\ t\in[0,T],\\
p_1(T)&=\Phi_x(T),
\end{aligned}
\right.
\end{equation}
and
\begin{equation}\label{eq3.12}
\left\{
\begin{aligned}
dp_2(t)&=-F_2(t)dt+q_{2,1}(t)dB_t+q_{2,2}(t)d\widetilde{B}_t,\ t\in[0,T],\\
p_2(T)&=\widetilde{\Phi}_{\mu}(T),
\end{aligned}
\right.
\end{equation}
where
\begin{equation*}
\begin{split}
F_1(t)=&f_x(t)+p_1(t)(b_x(t)+f_y(t)+\sigma_x(t)f_z(t))+q_{1,1}(t)(\sigma_x(t)+f_z(t)),\\
F_2(t)=&\widetilde{f}_{\mu,x}(t)+p_1(t)(\widetilde{b}_\mu(t)+f_z(t)\widetilde{\sigma}_\mu(t))+\widetilde{p}_1(t)\widetilde{f}_{\mu,y}(t)
        +q_{1,1}(t)\widetilde{\sigma}_\mu(t)+p_2(t)(\widetilde{b}_x(t)+f_y(t))\\
       &+\widehat{E}[\widehat{p}_2(t)(\widehat{b}_\mu^*(t)+\widehat{f}_{\mu,y}(t))]+q_{2,1}(t)f_z(t)
       +\widehat{E}[\widehat{q}_{2,2}(t)\widehat{\sigma}_\mu^*(t)]+q_{2,2}(t)\widetilde{\sigma}_x(t),
\end{split}
\end{equation*}
and for $P(t)=(P_1(t),P_2(t))$,
\begin{equation}\label{eq3.13}
\left\{
\begin{aligned}
dP_1(t)&=-G_1(t)dt+Q_{1,1}(t)dB_t+Q_{1,2}(t)d\widetilde{B}_t,\ t\in[0,T],\\
P_1(T)&=\Phi_{xx}(T),
\end{aligned}
\right.
\end{equation}
\begin{equation}\label{eq3.14}
\left\{
\begin{aligned}
dP_2(t)&=-G_2(t)dt+Q_{2,1}(t)dB_t+Q_{2,2}(t)d\widetilde{B}_t,\ t\in[0,T],\\
P_2(T)&=\widetilde{\Phi}_{\widetilde{x}\mu}(T),
\end{aligned}
\right.
\end{equation}
where
\begin{equation*}
\begin{split}
G_1(t)=&p_1(t)b_{xx}(t)+\sigma_{xx}(t)(f_z(t)p_1(t)+q_{1,1}(t))+P_1(t)(2b_x(t)+|\sigma_x(t)|^2+f_y(t)+2f_z(t)\sigma_x(t))\\
        &+\langle\partial^2_{(x,y,z)}f(t)(1,p_1(t),p_1(t)\sigma_x(t)+q_{1,1}(t))^\top,(1,p_1(t),p_1(t)\sigma_x(t)+q_{1,1}(t))^\top\rangle\\
        &+Q_{1,1}(t)(2\sigma_x(t)+f_z(t)),\\
G_2(t)=&p_1(t)\widetilde{b}_{\widetilde{x}\mu}(t)+q_{1,1}(t)\widetilde{\sigma}_{\widetilde{x}\mu}(t)+p_2(t)\widetilde{b}_{xx}(t)
        +\widehat{E}[\widehat{p}_2(t)\widehat{b}^*_{\widetilde{x}\mu}(t)]+q_{2,2}(t)\widetilde{\sigma}_{xx}(t)
        +\widehat{E}[\widehat{q}_{2,2}(t)\widehat{\sigma}^*_{\widetilde{x}\mu}(t)]\\
        &+2P_{2}(t)\widetilde{b}_x(t)+P_{2}(t)|\widetilde{\sigma}_x(t)|^2+2Q_{2,2}(t)\widetilde{\sigma}_x(t)
        +\widetilde{f}_{\mu,y}(t)\widetilde{P}_1(t)+\widehat{E}[\widehat{f}_{\mu,y}(t)\widehat{P}_2(t)]+f_y(t)P_2(t)\\
        &+\langle\partial_{(\widetilde{x},\widetilde{y})}\widetilde{f}_{\mu}(t)(1,p_1(t))^\top,(1,p_1(t))^\top\rangle
        +f_z(t)(p_1(t)\widetilde{\sigma}_{\widetilde{x}\mu}(t)+ Q_{2,1}(t)).
\end{split}
\end{equation*}
Notice that, in (\ref{eq3.11}) and (\ref{eq3.13}), the terminal values $\Phi_{x}(T)$ and $\Phi_{xx}(T)$ are in $L^2(\mathcal{F}_T; \mathbb{R})$, and the driving coefficients are adapted with respect to the filtration generated by $B$. Hence, we see
 $q_{1,2}(t)=0$ and $Q_{1,2}(t)=0$, $t\in[0,T]$.
Applying It\^{o}'s formula to $M_1(t)=p_1(t)(X_t^1+X_t^2)+\frac{1}{2}P_1(t)(X_t^1)^2$ , we get
\begin{equation}\label{eq3.15}
\begin{aligned}
&dM_1(t)
=d\(p_1(t)(X_t^1+X_t^2)+\frac{1}{2}P_1(t)(X_t^1)^2\)\\
=&\Bigr\{\alpha_{1,1}(t)I_{E_\varepsilon}(t)+\alpha_{1,2}(t)(X_t^1+X_t^2)+\widetilde{E}[\alpha_{1,3}(t)(\widetilde{X}_t^1+\widetilde{X}_t^2)]+\frac{1}{2}\alpha_{1,4}(t)(X_t^1)^2
     +\widetilde{E}[\frac{1}{2}\alpha_{1,5}(t)(\widetilde{X}_t^1)^2]\\
 &\ \  +\alpha_{1,6}(t)X_t^1I_{E_\varepsilon}(t)+\widetilde{E}[\alpha_{1,7}(t)X_t^1\widetilde{X}_t^1]+\widetilde{E}[\alpha_{1,8}(t)\widetilde{X}_t^1]I_{E_\varepsilon}(t)
    +\frac{1}{2}\alpha_{1,9}(t)(\widetilde{E}[\widetilde{\sigma}_{\mu}(t)\widetilde{X}_t^1])^2\Bigr\}dt\\
 &+\Bigr\{\beta_{1,1}(t)I_{E_\varepsilon}(t)+\beta_{1,2}(t)(X_t^1+X_t^2)+\widetilde{E}[\beta_{1,3}(t)(\widetilde{X}_t^1+\widetilde{X}_t^2)]+\frac{1}{2}\beta_{1,4}(t)(X_t^1)^2\\
 &\ \ \ \ +\widetilde{E}[\frac{1}{2}\beta_{1,5}(t)(\widetilde{X}_t^1)^2]+\beta_{1,6}(t)X_t^1I_{E_\varepsilon}(t)+\widetilde{E}[\beta_{1,7}(t)X_t^1\widetilde{X}_t^1]\Bigr\}dB_t,\ t\in[0,T],
\end{aligned}
\end{equation}
where
\begin{equation*}
\begin{split}
\alpha_{1,1}(t)&=p_1(t)\delta b(t)+q_{1,1}(t)\delta \sigma(t)+\frac{1}{2}P_1(t)(\delta \sigma(t))^2,\\
\alpha_{1,2}(t)&=p_1(t)b_x(t)+q_{1,1}(t)\sigma_x(t)-F_1(t),\ \ \
\alpha_{1,3}(t)=p_1(t)\widetilde{b}_{\mu}(t)+q_{1,1}(t)\widetilde{\sigma}_{\mu}(t),\\
\alpha_{1,4}(t)&=p_1(t)b_{xx}(t)+q_{1,1}\sigma_{xx}(t)+2P_1(t)b_x(t)+2Q_{1,1}(t)\sigma_x(t)+P_1(t)(\sigma_x(t))^2-G_1(t),\\
\alpha_{1,5}(t)&=p_1(t)\widetilde{b}_{\widetilde{x}\mu}(t)+q_{1,1}\widetilde{\sigma}_{\widetilde{x}\mu}(t),\ \ \ \
\alpha_{1,6}(t)=q_{1,1}(t)\delta\sigma_x(t)+P_1(t)\sigma_x(t)\delta\sigma(t)+Q_{1,1}(t)\delta\sigma(t),\\
\alpha_{1,7}(t)&=P_1(t)\widetilde{b}_{\mu}(t)+P_1(t)\sigma_x(t)\widetilde{\sigma}_{\mu}(t)+Q_{1,1}(t)\widetilde{\sigma}_{\mu}(t),\
\alpha_{1,8}(t)=P_1(t)\delta\sigma(t)\widetilde{\sigma}_{\mu}(t),\
\alpha_{1,9}(t)=P_1(t),\\
\beta_{1,1}(t)&=p_1(t)\delta\sigma(t),\ \ \ \
\beta_{1,2}(t)=q_{1,1}(t)+p_1(t)\sigma_x(t),\ \ \
\beta_{1,3}(t)=p_1(t)\widetilde{\sigma}_{\mu}(t),\\
\beta_{1,4}(t)&=p_{1}(t)\sigma_{xx}(t)+Q_{1,1}(t)+2P_1(t)\sigma_x(t),\ \ \
\beta_{1,5}(t)=p_1(t)\widetilde{\sigma}_{\widetilde{x}\mu}(t),\\
\beta_{1,6}(t)&=p_1(t)\delta\sigma_{x}(t)+P_1(t)\delta\sigma(t),\ \ \
\beta_{1,7}(t)=P_1(t)\widetilde{\sigma}_{\mu}(t).
\end{split}
\end{equation*}
Next, applying It\^{o}'s formula to $M_2(t)=\widetilde{E}\[p_2(t)(\widetilde{X}_t^1+\widetilde{X}_t^2)+\frac{1}{2}P_2(t)(\widetilde{X}_t^1)^2\]$, we have
\begin{equation}\label{eq3.16}
\begin{aligned}
&dM_2(t)
=\widetilde{E}\[(p_2(t)\widetilde{b}_x(t)+q_{2,2}(t)\widetilde{\sigma}_x(t)-F_2(t))(\widetilde{X}_t^1+\widetilde{X}_t^2)
        +p_2(t)E[\widetilde{b}_{\mu}^*(t)(X_t^1+X_t^2)]\]dt\\
&+\widetilde{E}\[\frac{1}{2}\(p_2(t)\widetilde{b}_{xx}(t)+q_{2,2}(t)\widetilde{\sigma}_{xx}(t)+2P_2(t)\widetilde{b}_{x}(t)
        +P_2(t)(\widetilde{\sigma}_x(t))^2+2Q_{2,2}(t)\widetilde{\sigma}_x(t)-G_2(t)\)(\widetilde{X}_t^1)^2\]dt\\
&+\widetilde{E}\[q_{2,2}(t)E[\widetilde{\sigma}_{\mu}^*(t)(X_t^1+X_t^2)]
        +\frac{1}{2}p_2(t)E[\widetilde{b}_{\widetilde{x}\mu}^*(t)(X_t^1)^2]
        +\frac{1}{2}q_{2,2}(t)E[\widetilde{\sigma}_{\widetilde{x}\mu}^*(t)(X_t^1)^2]\]dt\\
&+\widetilde{E}\[\frac{1}{2}P_2(t)(E[\widetilde{\sigma}_{\mu}^*(t)X_t^1])^2
        +\(p_2(t)\delta\widetilde{b}(t)+q_{2,2}(t)\delta\widetilde{\sigma}(t)+\frac{1}{2}P_2(t)(\delta\widetilde{\sigma}(t))^2\)I_{E_\varepsilon}(t)\]dt\\
&+\widetilde{E}\[\(P_2(t)\widetilde{\sigma}_x(t)\delta\widetilde{\sigma}(t)+q_{2,2}(t)\delta\widetilde{\sigma}_x(t)
        +Q_{2,2}(t)\delta\widetilde{\sigma}(t)\)\widetilde{X}_t^1I_{E_\varepsilon}(t)+P_2(t)\delta\widetilde{\sigma}(t)E[\widetilde{\sigma}_{\mu}^*(t)X_t^1]I_{E_\varepsilon}(t)\]dt\\
&+\widetilde{E}\[P_2(t)\widetilde{X}_t^1E[\widetilde{b}_{\mu}^*(t)X_t^1]+P_2(t)\widetilde{\sigma}_x(t)\widetilde{X}_t^1E[\widetilde{\sigma}_{\mu}^*(t)X_t^1]
        +Q_{2,2}(t)\widetilde{X}_t^1E[\widetilde{\sigma}_{\mu}^*(t)X_t^1]\]dt\\
&+\widetilde{E}\[q_{2,1}(t)(\widetilde{X}_t^1+\widetilde{X}_t^2)+\frac{1}{2}Q_{2,1}(t)(\widetilde{X}_t^1)^2\]dB_t,\ t\in[0,T].
\end{aligned}
\end{equation}
We introduce a new copy $(\widehat{\Omega},\widehat{\mathcal{F}},\widehat{P})$ of the probability space $(\Omega,\mathcal{F},P)$, and putting for $\widetilde{b}^*_\mu(t)=(\partial_\mu b)(t,\widetilde{X}_t^*,$ $P_{X_t^*},\widetilde{u}_t^*;X_t^*)$, $b_\mu(t,\omega,\widetilde{\omega}):=(\partial_\mu b)(t,X_t^*(\widetilde{\omega}),P_{X_t^*},u_t^*(\widetilde{\omega});X_t^*(\omega))$, we have the following relation:
\begin{equation}\label{relation}
\begin{split}
&\widetilde{E}\[p_2(t)E[\widetilde{b}_{\mu}^*(t)(X_t^1+X_t^2)]\]\\
=&\widetilde{E}\[p_2(t,\omega,\widetilde{\omega})E[b_{\mu}(t,\omega,\widetilde{\omega})(X_t^1(\omega)+X_t^2(\omega))]\]\\
=&\widetilde{E}\[p_2(t,\omega,\widetilde{\omega})\widehat{E}[b_{\mu}(t,\widehat{\omega},\widetilde{\omega})(X_t^1(\widehat{\omega})+X_t^2(\widehat{\omega}))]\]
=\widetilde{E}\[\widehat{E}[p_2(t,\omega,\widetilde{\omega})b_{\mu}(t,\widehat{\omega},\widetilde{\omega})(X_t^1(\widehat{\omega})+X_t^2(\widehat{\omega}))]\]\\
=&\widehat{E}\[\widetilde{E}[p_2(t,\omega,\widehat{\omega})b_{\mu}(t,\widehat{\omega},\widetilde{\omega})(X_t^1(\widetilde{\omega})+X_t^2(\widetilde{\omega}))]\]
=\widetilde{E}\[\widehat{E}[p_2(t,\omega,\widehat{\omega})b_{\mu}(t,\widehat{\omega},\widetilde{\omega})(X_t^1(\widetilde{\omega})+X_t^2(\widetilde{\omega}))]\]\\
=&\widetilde{E}\[\widehat{E}[p_2(t,\omega,\widehat{\omega})b_{\mu}(t,\widehat{\omega},\widetilde{\omega})](X_t^1(\widetilde{\omega})+X_t^2(\widetilde{\omega}))\]
=\widetilde{E}\[\widehat{E}[\widehat{p}_2(t)\widehat{b}^*_{\mu}(t)](\widetilde{X}_t^1+\widetilde{X}_t^2)\],\ t\in[0,T].
\end{split}
\end{equation}
Here, abusing notation, we have written $E[\zeta(\omega)],\ \widetilde{E}[\zeta(\widetilde{\omega})],\ \cdots $ for $\displaystyle\int_\Omega \zeta(\omega)P(d\omega)$, $\displaystyle\int_{\widetilde{\Omega}} \zeta(\widetilde{\omega})\widetilde{P}(d\widetilde{\omega}),\cdots$, respectively.
Thus, (\ref{eq3.16}) can be rewritten as
\begin{equation}\label{eq3.17}
\begin{aligned}
&dM_2(t)
=d\widetilde{E}\[p_2(t)(\widetilde{X}_t^1+\widetilde{X}_t^2)+\frac{1}{2}P_2(t)(\widetilde{X}_t^1)^2\]\\
=&\Bigr\{\widetilde{E}[\alpha_{2,1}(t)]I_{E_\varepsilon}(t)+\widetilde{E}[\alpha_{2,2}(t)(\widetilde{X}_t^1+\widetilde{X}_t^2)]
    +\frac{1}{2}\widetilde{E}[\alpha_{2,3}(t)(\widetilde{X}_t^1)^2]
    +\frac{1}{2}\widetilde{E}[\alpha_{2,4}(t)(E[\widetilde{\sigma}^*_{\mu}(t)X_t^1])^2]\\
    &+\widetilde{E}[\widehat{E}[\alpha_{2,5}(t)\widehat{X}_t^1\widetilde{X}_t^1]]
    +\widetilde{E}[\alpha_{2,6}(t)\widetilde{X}_t^1]I_{E_\varepsilon}(t)\Bigr\}dt
    +\Bigr\{\widetilde{E}[\beta_{2,2}(t)(\widetilde{X}_t^1+\widetilde{X}_t^2)] +\frac{1}{2}\widetilde{E}[\beta_{2,3}(t)(\widetilde{X}_t^1)^2]\Bigr\}dB_t,
\end{aligned}
\end{equation}
where\vspace{-0.8cm}
\begin{equation*}
\begin{split}
\alpha_{2,1}(t)=&p_2(t)\delta\widetilde{b}(t)+q_{2,2}(t)\delta\widetilde{\sigma}(t)+\frac{1}{2}P_2(t)(\delta\widetilde{\sigma}(t))^2,\\
\alpha_{2,2}(t)=&p_2(t)\widetilde{b}_x(t)+q_{2,2}(t)\widetilde{\sigma}_x(t)+\widehat{E}[\widehat{p}_2(t)\widehat{b}_{\mu}^*(t)]
        +\widehat{E}[\widehat{q}_{2,2}(t)\widehat{\sigma}_{\mu}^*(t)]-F_2(t),\\
\alpha_{2,3}(t)=&p_2(t)\widetilde{b}_{xx}(t)+q_{2,2}(t)\widetilde{\sigma}_{xx}(t)+2P_2(t)\widetilde{b}_{x}(t)
        +P_2(t)(\widetilde{\sigma}_x(t))^2+2Q_{2,2}(t)\widetilde{\sigma}_x(t)\\
        &+\widehat{E}[\widehat{p}_2(t)\widehat{b}_{\widetilde{x}\mu}^*(t)]
        +\widehat{E}[\widehat{q}_{2,2}(t)\widehat{\sigma}_{\widetilde{x}\mu}^*(t)]-G_2(t),\\
\alpha_{2,4}(t)=& P_2(t),\ \
\alpha_{2,5}(t)=\widehat{P}_2(t)\widehat{b}_{\mu}^*(t)+\widehat{P}_2(t)\widehat{\sigma}_{x}(t)\widehat{\sigma}_{\mu}^*(t)
        +\widehat{Q}_{2,2}(t)\widehat{\sigma}_{\mu}^*(t),\\
\alpha_{2,6}(t)=&q_{2,2}(t)\delta\widetilde{\sigma}_x(t)+P_2(t)\widetilde{\sigma}_x(t)\delta\widetilde{\sigma}(t)
        +Q_{2,2}(t)\delta\widetilde{\sigma}(t)+\widehat{E}[\widehat{P}_2(t)\delta\widehat{\sigma}(t)\widehat{\sigma}_{\mu}^*(t)],\\
\beta_{2,2}(t)=&q_{2,1}(t),\ \ \ \beta_{2,3}(t)=Q_{2,1}(t).
\end{split}
\end{equation*}
Let us define $M(t)=M_1(t)+M_2(t),\ t\in[0,T]$.
Combining (\ref{eq3.15}) and (\ref{eq3.17}), we see that
\begin{equation}\label{eq3.19}
\begin{split}
dM(t)=&\Bigr\{\alpha_1(t)I_{E_\varepsilon}(t)+\alpha_{2}(t)(X_t^1+X_t^2)+\widetilde{E}[\alpha_{3}(t)(\widetilde{X}_t^1+\widetilde{X}_t^2)]+\frac{1}{2}\alpha_{4}(t)(X_t^1)^2
    +\frac{1}{2}\widetilde{E}[\alpha_{5}(t)(\widetilde{X}_t^1)^2]\Bigr\}dt\\
 &+R_1^\varepsilon(t)dt+\Bigr\{\beta_{1}(t)I_{E_\varepsilon}(t)+\beta_{2}(t)(X_t^1+X_t^2)+\widetilde{E}[\beta_{3}(t)(\widetilde{X}_t^1+\widetilde{X}_t^2)]+\frac{1}{2}\beta_{4}(t)(X_t^1)^2\\
 &+\frac{1}{2}\widetilde{E}[\beta_{5}(t)(\widetilde{X}_t^1)^2]+\beta_{6}(t)X_t^1I_{E_\varepsilon}(t)+R_2^\varepsilon(t)\Bigr\}dB_t,
\end{split}
\end{equation}
where
\begin{equation*}
\begin{split}
\alpha_{1}(t)=&\alpha_{1,1}(t)+\widetilde{E}[\alpha_{2,1}(t)],\ \
\alpha_{2}(t)=\alpha_{1,2}(t),\ \
\alpha_{3}(t)=\alpha_{1,3}(t)+\alpha_{2,2}(t),\ \
\alpha_{4}(t)=\alpha_{1,4}(t),\\
\alpha_{5}(t)=&\alpha_{1,5}(t)+\alpha_{2,3}(t),\ \
\beta_{1}(t)=\beta_{1,1}(t),\ \
\beta_{2}(t)=\beta_{1,2}(t),\ \
\beta_{3}(t)=\beta_{1,3}(t)+\beta_{2,2}(t),\\
\beta_{4}(t)=&\beta_{1,4}(t),\ \
\beta_{5}(t)=\beta_{1,5}(t)+\beta_{2,3}(t),\ \
\beta_{6}(t)=\beta_{1,6}(t),\\
R_1^\varepsilon(t)=&\alpha_{1,6}(t)X_t^1I_{E_\varepsilon}(t)+\widetilde{E}[\alpha_{1,7}(t)X_t^1\widetilde{X}_t^1]+\widetilde{E}[\alpha_{1,8}(t)\widetilde{X}_t^1]I_{E_\varepsilon}(t)
+\frac{1}{2}\alpha_{1,9}(t)(\widetilde{E}[\widetilde{\sigma}_{\mu}(t)\widetilde{X}_t^1])^2\\
&+\frac{1}{2}\widetilde{E}[\alpha_{2,4}(t)(E[\widetilde{\sigma}^*_{\mu}(t)X_t^1])^2]
    +\widetilde{E}[\widehat{E}[\alpha_{2,5}(t)\widehat{X}_t^1\widetilde{X}_t^1]]
    +\widetilde{E}[\alpha_{2,6}(t)\widetilde{X}_t^1]I_{E_\varepsilon}(t),\\
R_2^\varepsilon(t)=&\widetilde{E}[\beta_{1,7}(t)X_t^1\widetilde{X}_t^1].
\end{split}
\end{equation*}
Now let us estimate $R_1^\varepsilon$ and $R_2^\varepsilon$. We begin with the term $\alpha_{1,6}(t)X_t^1I_{E_\varepsilon}(t)$ in $R_1^\varepsilon(t)$. Thanks to assumption \textbf{(A2)}, we have
\begin{equation*}
|\alpha_{1,6}(t)|\leq C(|q_{1,1}(t)|+|P_1(t)|+|Q_{1,1}(t)|)=:\zeta_1(t),
\end{equation*}
where $\zeta_1$ is independent of $\varepsilon$ and
\begin{equation*}
E\[\(\int_0^T |\zeta_1(t)|^2 dt\)^2\]\leq C\(E\[\sup_{t\in[0,T]}|P_1(t)|^4\]+E\[\(\int_0^T(|q_{1,1}(t)|^2+|Q_{1,1}(t)|^2)dt\)^2\]\)<\infty.
\end{equation*}
Hence,
\begin{equation*}
E\[\(\int_0^T |\alpha_{1,6}(t)|^2 I_{E_\varepsilon}(t)dt\)^2\]
\leq E\[\(\int_0^T|\zeta_1(t)|^2 I_{E_\varepsilon}(t)dt\)^2\]=:\overline{\rho}(\varepsilon),
\end{equation*}
where, thanks to the dominated convergence theorem, $\overline{\rho}(\varepsilon)\downarrow 0$, as $\varepsilon\downarrow 0$.
Consequently, we have
\begin{equation*}
\begin{split}
&E\[\(\int_0^T |\alpha_{1,6}(t)X_t^1|I_{E_\varepsilon}(t) dt\)^2\]
\leq \varepsilon E\[\int_0^T |\alpha_{1,6}(t)X_t^1|^2I_{E_\varepsilon}(t) dt\]\\
&\leq \varepsilon \(E\[\sup_{t\in[0,T]}|X_t^1|^4\]\)^{\frac{1}{2}}
        \(E\[\(\int_0^T |\alpha_{1,6}(t)|^2I_{E_\varepsilon}(t) dt\)^2\]\)^{\frac{1}{2}}
\leq C\varepsilon^2 \sqrt{\overline{\rho}(\varepsilon)},\ \varepsilon>0.
\end{split}
\end{equation*}
Similarly, but with application of Proposition \ref{lem4.1}, we get
$$E\[\(\displaystyle\int_0^T|\widetilde{E}[\alpha_{1,8}(t)\widetilde{X}_t^1]|I_{E_\varepsilon}(t)dt\)^2\]\leq C\varepsilon^2\overline{\rho}(\varepsilon)$$
and   $$E\[\(\displaystyle\int_0^T|\widetilde{E}[\alpha_{2,6}(t)\widetilde{X}_t^1]|I_{E_\varepsilon}(t)dt\)^2\]\leq C\varepsilon^2\overline{\rho}(\varepsilon),\ \varepsilon>0,$$
for some $\overline{\rho}:\mathbb{R}_+\rightarrow\mathbb{R}_+$ with $\overline{\rho}(\varepsilon)\downarrow0$ as $\varepsilon\downarrow0$.\\
Furthermore,
\begin{equation*}
\begin{aligned}
&|\widetilde{E}[\alpha_{1,7}(t)X_t^1\widetilde{X}_t^1]|\\
=&|P_1(t)X_t^1\widetilde{E}[\widetilde{b}_{\mu}(t)\widetilde{X}_t^1]
    +P_1(t)\sigma_x(t)X_t^1\widetilde{E}[\widetilde{\sigma}_{\mu}(t)\widetilde{X}_t^1]
    +Q_{1,1}(t)X_t^1\widetilde{E}[\widetilde{\sigma}_{\mu}(t)\widetilde{X}_t^1]|\\
\leq& C(|P_1(t)|+|Q_{1,1}(t)|)|X_t^1|(|\widetilde{E}[\widetilde{b}_{\mu}(t)\widetilde{X}_t^1]|
        +|\widetilde{E}[\widetilde{\sigma}_{\mu}(t)\widetilde{X}_t^1]|),\ t\in[0,T].
\end{aligned}
\end{equation*}
As $\theta:=\widetilde{b}_{\mu},\ \widetilde{\sigma}_{\mu}\in\mathcal{H}^2(0,T)$ is bounded, we can apply Proposition 5.2 in \cite{BCL2021} and conclude that there exists $\rho_t:\mathbb{R}_+\rightarrow\mathbb{R}_+$ with $\rho_t(\varepsilon)\downarrow0$, as $\varepsilon\downarrow0$, $|\rho_t(\varepsilon)|\leq C,\ t\in[0,T],\ \varepsilon>0$ (for some $C$ depending only on $\theta$), such that
\begin{equation}\label{18-1}
|\widetilde{E}[\theta_t\widetilde{X}_t^1]|\leq \rho_t(\varepsilon)\sqrt{\varepsilon},\ \varepsilon>0,\ t\in[0,T].
\end{equation}
Thus,
$$|\widetilde{E}[\alpha_{1,7}(t)X_t^1\widetilde{X}_t^1]|\leq C\sqrt{\varepsilon}\rho_t(\varepsilon)(|P_1(t)|+|Q_{1,1}(t)|)|X_t^1|$$
and
\begin{equation*}
\begin{aligned}
&E\[\(\int_0^T|\widetilde{E}[\alpha_{1,7}(t)X_t^1\widetilde{X}_t^1]|dt\)^2\]
\leq C\varepsilon E\[\(\int_0^T\rho_t(\varepsilon)(|P_1(t)|+|Q_{1,1}(t)|)dt\cdot\sup_{t\in[0,T]}|X_t^1|\)^2\]\\
&\leq C\varepsilon\int_0^T\rho_t(\varepsilon)^2dt
    \(E\[\(\sup_{t\in[0,T]}|P_1(t)|^2+\int_0^T|Q_{1,1}(t)|^2dt\)^2\]\)^{\frac{1}{2}}
    \(E[\sup_{t\in[0,T]}|X_t^1|^4]\)^{\frac{1}{2}}\\
&\leq C\rho(\varepsilon)\varepsilon^2,\ \varepsilon>0,
\end{aligned}
\end{equation*}
where $\rho(\varepsilon):=\displaystyle\int_0^T\rho_t(\varepsilon)^2dt\downarrow0$, as $\varepsilon\downarrow0$.

A similar argument using Proposition \ref{lem4.1} allows to show that
$$E\[\(\displaystyle\int_0^T|\widetilde{E}[\widehat{E}[\alpha_{2,5}(t)\widehat{X}_t^1\widetilde{X}_t^1]]|dt\)^2\]\leq C\varepsilon^2\rho(\varepsilon), \ \text{for}\ \rho(\varepsilon)\downarrow0.$$
For $\alpha_{1,9}(t)(\widetilde{E}[\widetilde{\sigma}_{\mu}(t)\widetilde{X}_t^1])^2$ and $\widetilde{E}[\alpha_{2,4}(t)(E[\widetilde{\sigma}^*_{\mu}(t)X_t^1])^2]$ in $R_1^\varepsilon$, we just consider the first term here. From the above estimate of $\widetilde{E}[\widetilde{\sigma}_{\mu}(t)\widetilde{X}_t^1]$, we see with the help of Proposition \ref{lem4.1}
\begin{equation*}
\begin{aligned}
&E\[\(\int_0^T|\alpha_{1,9}(t)(\widetilde{E}[\widetilde{\sigma}_{\mu}(t)\widetilde{X}_t^1])^2|dt\)^2\]
\leq CE\[\sup_{t\in[0,T]}|P_1(t)|^2\(\int_0^T|\widetilde{E}[\widetilde{\sigma}_{\mu}(t)\widetilde{X}_t^1]|^2dt\)^2\]\\
&\leq CE[\sup_{t\in[0,T]}|P_1(t)|^2]\(\int_0^T(\varepsilon\rho_t(\varepsilon)^2)dt\)^2
\leq C\varepsilon^2\rho(\varepsilon).
\end{aligned}
\end{equation*}
Using the same method we get $E\[\(\displaystyle\int_0^T|\widetilde{E}[\alpha_{2,4}(t)(E[\widetilde{\sigma}^*_{\mu}X_t^1])^2]|dt\)^2\]\leq C\varepsilon^2\rho(\varepsilon).$\\
Summarizing the above estimates, we obtain
\begin{equation}\label{eq4.14}
E\[\(\displaystyle\int_0^T|R_1^\varepsilon(t)|dt\)^2\]\leq C\varepsilon^2\rho(\varepsilon),
\end{equation}
where $\rho(\varepsilon)\rightarrow0$, as $\varepsilon\rightarrow0$.
Finally, concerning $R_2^\varepsilon$ we observe that, as $\beta_{1,7}(t)=P_1(t)\widetilde{\sigma}_{\mu}(t)$, from (\ref{18-1}),
\begin{equation}\label{eq4.15}
\begin{aligned}
&E\[\(\int_0^TR_2^\varepsilon(t)dB_t\)^2\]
=E\[\(\int_0^T\widetilde{E}[\beta_{1,7}(t)X_t^1\widetilde{X}_t^1]dB_t\)^2\]
= E\[\int_0^T|\widetilde{E}[\beta_{1,7}(t)X_t^1\widetilde{X}_t^1]|^2dt\]\\
&\leq E\[\(\sup_{t\in[0,T]}|P_1(t)|^2\)\(\sup_{t\in[0,T]}|X_t^1|^2\)
        \int_0^T|\widetilde{E}[\widetilde{\sigma}_{\mu}(t)\widetilde{X}_t^1]|^2dt\]\\
&\leq \(E\[\sup_{t\in[0,T]}|P_1(t)|^4\]\)^\frac{1}{2}
    \(E\[\sup_{t\in[0,T]}|X_t^1|^4\]\)^\frac{1}{2}\cdot\varepsilon\rho(\varepsilon)
\leq C\varepsilon^2\rho(\varepsilon).
\end{aligned}
\end{equation}
\section{The Estimates of $Y$}
Let us consider the following BSDE:
\begin{equation}\label{eq5.1}
\left\{
\begin{aligned}
d\overline{\overline{Y}}_t&=-\(\widehat{E}[\widehat{f}_{\mu,y}(t)\widehat{\overline{\overline{Y}}}_t]
    +f_y(t)\overline{\overline{Y}}_t+f_z(t)\overline{\overline{Z}}_t+(\alpha_1(t)+\Delta f(t))I_{E_\varepsilon}(t)\)dt
    +\overline{\overline{Z}}_tdB_t,\ t\in[0,T],\\
\overline{\overline{Y}}_T&=0,
\end{aligned}
\right.
\end{equation}
\ \ \ \ \ \ where\vspace{-0.7cm}
\begin{equation*}
\begin{split}
\Delta f(t)I_{E_\varepsilon}(t)
            =&f(t,X^*_t,Y^*_t,Z^*_t+p_1(t)\delta\sigma(t)I_{E_\varepsilon}(t),P_{(X^*_t,Y^*_t)},u_t^\varepsilon)-f(t)\\
            =&(f(t,X^*_t,Y^*_t,Z^*_t+p_1(t)\delta\sigma(t),P_{(X^*_t,Y^*_t)},u_t)-f(t))I_{E_\varepsilon}(t).
\end{split}
\end{equation*}
\begin{proposition}\label{prop5.1}
Assume that \textbf{(A1)}-\textbf{(A4)} hold true. Then, for all $\kappa\in (1\vee p_K^*,\kappa_0)$, where $K:=\displaystyle\int_0^\cdot f_z(t)dB_t$, we have\vspace{-0.4cm}
\begin{equation}\label{eq5.2}
E\[\sup_{t\in[0,T]}|\overline{\overline{Y}}_t|^\kappa+\(\int_0^T|\overline{\overline{Z}}_t|^2dt\)^\frac{\kappa}{2}\]
\leq C\varepsilon^{\frac{\kappa}{2}}\rho(\varepsilon),
\end{equation}
where $\rho(\cdot)$ is a positive function defined on $(0,\infty)$ such that $\rho(\varepsilon)\rightarrow 0$, as $\varepsilon\rightarrow 0$.
\end{proposition}
\begin{proof}
We observe that the processes $f_{\mu,y}(\cdot)$ and $f_y(\cdot)$ are $\mathbb{F}$-adapted and bounded, while for the $\mathbb{F}$-adapted process $f_z(\cdot)$ we have $f_z(t)\leq L+\gamma|Z_t^*|,\ t\in[0,T]$ (see assumption \textbf{(A4)}). From Proposition \ref{thm2.1} we know that $Z^*\cdot B$ is a $BMO$ martingale. Then, obviously, also $f_z(\cdot)\cdot B\in BMO_2$, and we can apply Proposition \ref{prop4.1}. This yields, for $\kappa\in(1\vee p_K^*,\kappa_0)$,
\begin{equation*}
\begin{split}
E\[\sup_{t\in[0,T]}|\overline{\overline{Y}}_t|^\kappa+\(\int_0^T|\overline{\overline{Z}}_t|^2dt\)^\frac{\kappa}{2}\]
\leq C\(E\[\(\int_0^T|\alpha_1(t)+\Delta f(t)|I_{E_\varepsilon}(t)dt\)^{\kappa_0}\]\)^\frac{\kappa}{\kappa_0}.
\end{split}
\end{equation*}
Note that due to assumption \textbf{(A.4)},
\begin{equation*}
\begin{aligned}
|\Delta f(t)|I_{E_\varepsilon}(t)
&\leq C(1+|X_t^*|+|Y_t^*|+|Z_t^*|(1+|p_1(t)\delta\sigma(t)|)^2)I_{E_\varepsilon}(t)\\
&\leq C(1+|X_t^*|+|Y_t^*|+|Z_t^*|(1+|p_1(t)|^2))I_{E_\varepsilon}(t).
\end{aligned}
\end{equation*}
Thus, by Proposition \ref{prop3.1} and a standard estimate for BSDE (\ref{eq3.11}),
\begin{equation*}
\begin{aligned}
&E\[\(\int_0^T|\Delta f(t)|I_{E_\varepsilon}(t)dt\)^{\kappa_0}\]\\
\leq & CE\[1+\sup_{t\in[0,T]}(|X_t^*|+|Y_t^*|)^{\kappa_0}\]\varepsilon^{\kappa_0}
    +CE\[\(1+\sup_{t\in[0,T]}|p_1(t)|^2\)^{\kappa_0}\(\int_0^T|Z_t^*|I_{E_\varepsilon}(t)dt\)^{\kappa_0}\]\\
\leq & C\varepsilon^{\kappa_0}+C\(E\[\(\int_0^T|Z_t^*| I_{E_\varepsilon}(t)dt\)^{2\kappa_0}\]\)^{\frac{1}{2}},\ \varepsilon>0.
\end{aligned}
\end{equation*}
Note that
\begin{equation*}
\(E\[\(\int_0^T|Z_t^*|I_{E_\varepsilon}(t)dt\)^{2\kappa_0}\]\)^\frac{1}{2}
\leq \(E\[\(\int_0^T|Z_t^*|^2I_{E_\varepsilon}(t)dt\)^{\kappa_0}\]\)^\frac{1}{2}\varepsilon^\frac{\kappa_0}{2},
\end{equation*}
where, as $E\[\(\displaystyle\int_0^T|Z_t^*|^2dt\)^{\kappa_0}\]<\infty$ (see Remark  \ref{re2.1}), it follows from the dominated convergence theorem that
\begin{equation*}
\rho^1(\varepsilon):=C\(E\[\(\int_0^T|Z_t^*|^2I_{E_\varepsilon}(t)dt\)^{\kappa_0}\]\)^\frac{1}{2}\downarrow0,\ \ \text{as}\ \ \varepsilon\downarrow0.
\end{equation*}
Hence,
\begin{equation*}
E\[\(\int_0^T|\Delta f(t)|I_{E_\varepsilon}(t)dt\)^{\kappa_0}\]\leq C\rho^1(\varepsilon)\varepsilon^\frac{\kappa_0}{2},\ \ \varepsilon>0,\ \ \text{with}\ \ \rho^1(\varepsilon)\downarrow0\ \  \text{as}\ \ \varepsilon\downarrow0.
\end{equation*}
With similar arguments using the definition of $\alpha_1(t)$, we get
\begin{equation*}
E\[\(\int_0^T|\alpha_1(t)|I_{E_\varepsilon}(t)dt\)^{\kappa_0}\]
\leq C\rho^2(\varepsilon)\varepsilon^\frac{\kappa_0}{2},\ \varepsilon>0.
\end{equation*}
Hence,
\begin{equation*}
\(E\[\(\int_0^T|\alpha_1(t)+\Delta f(t)|I_{E_\varepsilon}(t)dt\)^{\kappa_0}\]\)^\frac{\kappa}{\kappa_0}
\leq C\rho(\varepsilon)\varepsilon^\frac{\kappa}{2},\ \varepsilon>0,
\end{equation*}
where $\rho(\varepsilon)=\(\rho^1(\varepsilon)+\rho^2)(\varepsilon)\)^\frac{\kappa}{\kappa_0}$.
The proof is complete.
\end{proof}
Let $(\overline{Y}^\varepsilon,\overline{Z}^\varepsilon)$ and $(\check{Y}^\varepsilon,\check{Z}^\varepsilon)$ be the unique solution of the following BSDE, respectively:
\begin{equation}\label{eq5.3}
\left\{
\begin{aligned}
d\overline{Y}^\varepsilon_t=&-\Bigr\{f(t,X_t^\varepsilon,Y_t^\varepsilon,Z_t^\varepsilon,P_{(X_t^\varepsilon,Y_t^\varepsilon)},u_t^\varepsilon)
+\alpha_1(t)I_{E_\varepsilon}(t)+\alpha_{2}(t)(X_t^1+X_t^2)+\widetilde{E}[\alpha_{3}(t)(\widetilde{X}_t^1+\widetilde{X}_t^2)]\\
&\ \ \ \ +\frac{1}{2}\alpha_{4}(t)(X_t^1)^2+\frac{1}{2}\widetilde{E}[\alpha_{5}(t)(\widetilde{X}_t^1)^2]\Bigr\}dt
+\overline{Z}_t^\varepsilon dB_t,\ t\in[0,T],\\
\overline{Y}^\varepsilon_T=&\Phi(X^*_T,P_{X^*_T}),
\end{aligned}
\right.
\end{equation}
and
\begin{equation}\label{eq5.4}
\left\{
\begin{aligned}
d\check{Y}^\varepsilon_t=&-R_1^\varepsilon(t)dt+(\check{Z}_t^\varepsilon-R_2^\varepsilon(t))dB_t,\ t\in[0,T],\\
\check{Y}^\varepsilon_T=&\Phi(X^\varepsilon_T,P_{X^\varepsilon_T})-\Phi(X^*_T,P_{X^*_T})-M(T).
\end{aligned}
\right.
\end{equation}
We notice that the coefficients of (\ref{eq5.3}) and (\ref{eq5.4}) do not depend on the solutions, so by the classical BSDE theory, it allows to show the existence and the uniqueness for the both BSDEs.
Then, since $\overline{Y}^\varepsilon_T+\check{Y}^\varepsilon_T=Y^\varepsilon_T-M(T)$ and the sum of the driving coefficients of $\overline{Y}^\varepsilon$ and $\check{Y}^\varepsilon$ equals to the difference between that of $Y^\varepsilon$ and the drift coefficient of $M$ (see (\ref{eq3.19})), it follows from the uniqueness of the solution of the involved BSDEs that
\begin{equation}\label{eq5.15}
\overline{Y}_t^\varepsilon+\check{Y}_t^\varepsilon=Y^\varepsilon_t-M(t),\ t\in[0,T],
\end{equation}
and
\begin{equation*}
\begin{aligned}
\overline{Z}_t^\varepsilon+(\check{Z}_t^\varepsilon-R_2^\varepsilon(t))
    =&Z_t^\varepsilon-\Bigr\{\beta_1(t)I_{E_\varepsilon}(t)+\beta_2(t)(X_t^1+X_t^2)
        +\widetilde{E}[\beta_3(t)(\widetilde{X}_t^1+\widetilde{X}_t^2)]+\frac{1}{2}\beta_4(t)(X_t^1)^2\\
        &\ \ \ \ \ \ \  +\frac{1}{2}\widetilde{E}[\beta_5(t)(\widetilde{X}_t^1)^2]+\beta_6(t)X_t^1I_{E_\varepsilon}(t)+R_2^\varepsilon(t)\Bigr\},\ t\in[0,T].
\end{aligned}
\end{equation*}
Consequently,
\begin{equation}\label{eq5.19}
\begin{aligned}
\check{Z}_t^\varepsilon=&Z_t^\varepsilon-\overline{Z}_t^\varepsilon-\(\beta_1(t)I_{E_\varepsilon}(t)+\beta_2(t)(X_t^1+X_t^2)
        +\widetilde{E}[\beta_3(t)(\widetilde{X}_t^1+\widetilde{X}_t^2)]+\frac{1}{2}\beta_4(t)(X_t^1)^2\\
        & \quad  \qquad  \qquad +\frac{1}{2}\widetilde{E}[\beta_5(t)(\widetilde{X}_t^1)^2]+\beta_6(t)X_t^1I_{E_\varepsilon}(t)\),\ t\in[0,T].
\end{aligned}
\end{equation}

We give the estimates of (\ref{eq5.4}) in the following proposition.
\begin{proposition}\label{prop5.2}
Under the assumptions \textbf{(A1)}-\textbf{(A4)}, we have
\begin{equation}\label{eq5.5}
E\[\sup_{t\in[0,T]}|\check{Y}_t^\varepsilon|^2+\int_0^T|\check{Z}_t^\varepsilon|^2dt\]
\leq C\varepsilon^2\rho(\varepsilon),
\end{equation}
where $\rho:\mathbb{R}_+\rightarrow\mathbb{R}_+$, and $\rho(\varepsilon)\rightarrow0$ as $\varepsilon\rightarrow0$.
\end{proposition}
\begin{proof}
BSDE (\ref{eq5.4}) is a classical one. Indeed, $R_1^\varepsilon(\cdot)\in \mathcal{H}^2(0,T;\mathbb{R})$, $R_2^\varepsilon(\cdot)\in \mathcal{S}^2(0,T;\mathbb{R})$ and $\check{Y}_T^\varepsilon\in L^2(\mathcal{F}_T;\mathbb{R})$. Consequently, by a BSDE standard estimate as well as the estimates (\ref{eq4.14}) and (\ref{eq4.15}) for  $R_1^\varepsilon$ and $R_2^\varepsilon$, respectively, we have
\begin{equation*}
\begin{split}
E\[\sup_{t\in[0,T]}|\check{Y}_t^\varepsilon|^2+\int_0^T|\check{Z}_t^\varepsilon|^2dt\]
\leq C\(E[|\check{Y}_T^\varepsilon|^2]+E[(\int_0^T|R_1^\varepsilon(t)|dt)^2]+E[(\int_0^T|R_2^\varepsilon(t)|^2dt)]\)
\leq C\varepsilon^2\rho(\varepsilon).
\end{split}
\end{equation*}
Indeed, as
\begin{equation*}
\begin{aligned}
\check{Y}_T^\varepsilon=& \Phi(X_T^\varepsilon,P_{X_T^\varepsilon})-\Phi(X_T^*,P_{X_T^*})-M(T)\\
    =& \Phi(X_T^\varepsilon,P_{X_T^\varepsilon})-\Phi(X_T^*,P_{X_T^*})
        -\Bigr\{\Phi_x(T)(X_T^1+X_T^2)+\frac{1}{2}\Phi_{xx}(T)(X_T^1)^2\\
       & +\widetilde{E}\[\widetilde{\Phi}_\mu(T)(\widetilde{X}_T^1+\widetilde{X}_T^2)
        +\frac{1}{2}\widetilde{\Phi}_{\widetilde{x}\mu}(T)(\widetilde{X}_T^1)^2\]\Bigr\},
\end{aligned}
\end{equation*}
it follows from (\ref{eq3.18}) that $E[|\check{Y}_T^\varepsilon|^2]=o(\varepsilon^2)$. The proof is complete now.
\end{proof}
Let us define $\delta \overline{Y}_t:=\overline{Y}^\varepsilon_t-Y^*_t-\overline{\overline{Y}}_t$ and
$\delta \overline{Z}_t:=\overline{Z}^\varepsilon_t-Z^*_t-\overline{\overline{Z}}_t$, $t\in[0,T]$. Then from (\ref{eq5.1}), (\ref{eq5.3}) and the BSDE in (\ref{eq2.5}), it follows that $(\delta \overline{Y},\delta \overline{Z})$ is the unique solution of the  following BSDE:
\begin{equation}\label{eq5.11}
\left\{
\begin{split}
d(\delta \overline{Y}_t)=&
    -\Bigr\{(f(t,X^\varepsilon_t,Y^\varepsilon_t,Z^\varepsilon_t,P_{(X^\varepsilon_t,Y^\varepsilon_t)},u_t^\varepsilon)
    -f(t,X^*_t,Y^*_t,Z^*_t,P_{(X^*_t,Y^*_t)},u_t^*))
    -\widehat{E}[\widehat{f}_{\mu,y}(t)\widehat{\overline{\overline{Y}}}_t]\\
    &+\alpha_{2}(t)(X_t^1+X_t^2)+\widetilde{E}[\alpha_{3}(t)(\widetilde{X}_t^1+\widetilde{X}_t^2)]+\frac{1}{2}\alpha_{4}(t)(X_t^1)^2
    +\frac{1}{2}\widetilde{E}[\alpha_{5}(t)(\widetilde{X}_t^1)^2]\Bigr\}dt\\
    &+\Bigr\{f_y(t)\overline{\overline{Y}}_t+f_z(t)\overline{\overline{Z}}_t+\Delta f(t)I_{E_\varepsilon}(t)\Bigr\}dt+\delta\overline{Z}_tdB_t,\\
\delta \overline{Y}_T=&0.
\end{split}
\right.
\end{equation}
Remark that, as the driving coefficient of BSDE (\ref{eq5.11}) does not depend on $(\delta \overline{Y},\delta \overline{Z})$, the uniqueness is a consequence of that of the martingale representation.
\begin{proposition}\label{prop5.3}
Let the assumptions \textbf{(A1)}-\textbf{(A4)} hold true. Then, for $1<p<2$,
\begin{equation}\label{eq5.6}
E\[\sup_{t\in[0,T]}|\delta\overline{Y}_t|^p+\(\int_0^T|\delta\overline{Z}_t|^2dt\)^{\frac{p}{2}}\]
\leq C\varepsilon^p\rho(\varepsilon),
\end{equation}
where $\rho:\mathbb{R}_+\rightarrow\mathbb{R}_+$ with $\rho(\varepsilon)\rightarrow0$, as $\varepsilon\rightarrow0$.
\end{proposition}
\begin{proof}
For the sake of shortness, we introduce the following notations.
\begin{equation*}
\begin{split}
\beta(t):=&\beta_{2}(t)(X_t^1+X_t^2)+\widetilde{E}[\beta_{3}(t)(\widetilde{X}_t^1+\widetilde{X}_t^2)]+\frac{1}{2}\beta_{4}(t)(X_t^1)^2
            +\frac{1}{2}\widetilde{E}[\beta_{5}(t)(\widetilde{X}_t^1)^2]+\beta_{6}(t)X_t^1I_{E_\varepsilon}(t),\\
\Theta_t^*:=&(X^*_t,Y_t^*,Z_t^*+p_1(t)\delta\sigma(t)I_{E_\varepsilon}(t)),\ \ \
\Theta_t:=(X^*_t,Y_t^*,Z_t^*),\\
\Theta_t^1:=&(X^*_t+X^1_t+X^2_t,\overline{Y}_t^\varepsilon+M(t),\overline{Z}_t^\varepsilon+\beta(t)+p_1(t)\delta\sigma(t)I_{E_\varepsilon}(t)),\\
\Theta_t^2:=&(X^*_t+X^1_t+X^2_t,Y_t^*+\overline{\overline{Y}}_t+M(t),Z_t^*+\overline{\overline{Z}}_t+\beta(t)+p_1(t)\delta\sigma(t)I_{E_\varepsilon}(t)),\\
\Theta_t^3:=&(X^*_t+X^1_t+X^2_t,Y_t^*+\overline{\overline{Y}}_t+M(t),Z_t^*+\overline{\overline{Z}}_t+\beta(t)),\ t\in[0,T].
\end{split}
\end{equation*}
Recall (\ref{eq5.15}), $Y^\varepsilon_t=\overline{Y}^\varepsilon_t+\check{Y}^\varepsilon_t+M(t)$, $t\in[0,T]$. Hence
\begin{equation*}
\begin{split}
&f(t,X^\varepsilon_t,Y^\varepsilon_t,Z^\varepsilon_t,P_{(X^\varepsilon_t,Y^\varepsilon_t)},u_t^\varepsilon)
    -f(t,X^*_t,Y^*_t,Z^*_t,P_{(X^*_t,Y^*_t)},u_t^*)\\
=&\Delta f(t)I_{E_\varepsilon}(t)+
    f(t,X^\varepsilon_t,Y^\varepsilon_t,Z^\varepsilon_t,P_{(X^\varepsilon_t,Y^\varepsilon_t)},u_t^\varepsilon)
    -f(t,X^*_t,Y^*_t,Z^*_t+p_1(t)\delta\sigma(t)I_{E_\varepsilon}(t),P_{(X^*_t,Y^*_t)},u_t^\varepsilon)\\
    &-f(t,\Theta_t^1,P_{(X^*_t+X^1_t+X^2_t,\overline{Y}_t^\varepsilon+M(t))},u_t^\varepsilon)
    +f(t,\Theta_t^1,P_{(X^*_t+X^1_t+X^2_t,\overline{Y}_t^\varepsilon+M(t))},u_t^\varepsilon)\\
    &-f(t,\Theta_t^2,P_{(X^*_t+X^1_t+X^2_t,Y_t^*+\overline{\overline{Y}}_t+M(t))},u_t^\varepsilon)
    +f(t,\Theta_t^2,P_{(X^*_t+X^1_t+X^2_t,Y_t^*+\overline{\overline{Y}}_t+M(t))},u_t^\varepsilon)\\
=&\Delta f(t)I_{E_\varepsilon}(t)+I_1(t)+I_2(t)+I_3(t),
\end{split}
\end{equation*}
where
\begin{equation}\label{eq5.18}
\begin{aligned}
I_1(t):=&f(t,\Theta_t^1,P_{(X^*_t+X^1_t+X^2_t,\overline{Y}_t^\varepsilon+M(t))},u_t^\varepsilon)
        -f(t,\Theta_t^2,P_{(X^*_t+X^1_t+X^2_t,Y_t^*+\overline{\overline{Y}}_t+M(t))},u_t^\varepsilon),\\
I_2(t):=&f(t,\Theta_t^2,P_{(X^*_t+X^1_t+X^2_t,Y_t^*+\overline{\overline{Y}}_t+M(t))},u_t^\varepsilon)
        -f(t,X^*_t,Y^*_t,Z^*_t+p_1(t)\delta\sigma(t)I_{E_\varepsilon}(t),P_{(X^*_t,Y^*_t)},u_t^\varepsilon),\\
I_3(t):=&f(t,X^\varepsilon_t,Y^\varepsilon_t,Z^\varepsilon_t,P_{(X^\varepsilon_t,Y^\varepsilon_t)},u_t^\varepsilon)
        -f(t,\Theta_t^1,P_{(X^*_t+X^1_t+X^2_t,\overline{Y}_t^\varepsilon+M(t))},u_t^\varepsilon),\ t\in[0,T].
\end{aligned}
\end{equation}
For $I_1$, we have
\begin{equation*}
\begin{split}
I_1(t)=&f(t,\Theta_t^1,P_{(X^*_t+X^1_t+X^2_t,\overline{Y}_t^\varepsilon+M(t))},u_t^\varepsilon)
        -f(t,\Theta_t^2,P_{(X^*_t+X^1_t+X^2_t,Y_t^*+\overline{\overline{Y}}_t+M(t))},u_t^\varepsilon)\\
      =&\Bigr\langle\int_0^1\partial_{(x,y,z)}f(t,\Theta_t^2+\lambda(\Theta_t^1-\Theta_t^2),
            P_{(X^*_t+X^1_t+X^2_t,Y_t^*+\overline{\overline{Y}}_t+M(t)+\lambda\delta\overline{Y}_t)},u_t^\varepsilon)d\lambda,\Theta_t^1-\Theta_t^2\Bigr\rangle\\
       &+\widehat{E}\[\int_0^1\partial_{\mu, y}f(t,\Theta_t^2+\lambda(\Theta_t^1-\Theta_t^2),
            P_{(X^*_t+X^1_t+X^2_t,Y_t^*+\overline{\overline{Y}}_t+M(t)+\lambda\delta\overline{Y}_t)},u_t^\varepsilon;\\
            &\ \ \ \ \ \ \ \ \ \ \ \ \ \ \ \ \ \ \ \ \ \ (\widehat{X}^*_t+\widehat{X}^1_t+\widehat{X}^2_t,\widehat{Y}_t^*+\widehat{\overline{\overline{Y}}}_t+\widehat{M}(t)+\lambda\widehat{\delta\overline{Y}}_t))d\lambda\cdot\widehat{\delta\overline{Y}}_t\],\ t\in[0,T].
\end{split}
\end{equation*}
To shorten the notations, we put
\begin{equation*}
\begin{aligned}
A_1(t):=&\int_0^1\partial_{(x,y,z)}f(t,\Theta_t^2+\lambda(\Theta_t^1-\Theta_t^2),
            P_{(X^*_t+X^1_t+X^2_t,Y_t^*+\overline{\overline{Y}}_t+M(t)+\lambda\delta\overline{Y}_t)},u_t^\varepsilon)d\lambda,\\
\widehat{A}_1(t):=&\int_0^1\partial_{\mu, y}f(t,\Theta_t^2+\lambda(\Theta_t^1-\Theta_t^2),
            P_{(X^*_t+X^1_t+X^2_t,Y_t^*+\overline{\overline{Y}}_t+M(t)+\lambda\delta\overline{Y}_t)},u_t^\varepsilon;\\
            &\ \ \ \ \ \ \ \ \ \ \ \ \ \ \ (\widehat{X}^*_t+\widehat{X}^1_t+\widehat{X}^2_t,\widehat{Y}_t^*+\widehat{\overline{\overline{Y}}}_t+\widehat{M}(t)+\lambda\delta\widehat{\overline{Y}}_t))d\lambda.
\end{aligned}
\end{equation*}
$I_1(t)$ can be rewritten as
\begin{equation}\label{eq5.7}
I_1(t)=\langle A_1(t),\Theta_t^1-\Theta_t^2\rangle+\widehat{E}[\widehat{A}_1(t)\cdot\widehat{\delta\overline{Y}}_t].
\end{equation}
For $I_2(t)$, we have
\begin{equation}\label{eq5.8}
\begin{split}
I_2(t)=&f(t,\Theta_t^2,P_{(X^*_t+X^1_t+X^2_t,Y_t^*+\overline{\overline{Y}}_t+M(t))},u_t^\varepsilon)
        -f(t,\Theta_t^*,P_{(X^*_t,Y^*_t)},u_t^\varepsilon)\\
      =&\Bigr\langle\int_0^1\partial_{(x,y,z)}f(t,\Theta_t^*+\lambda(\Theta_t^2-\Theta_t^*),
            P_{(X^*_t+\lambda(X^1_t+X^2_t),Y_t^*+\lambda(\overline{\overline{Y}}_t+M(t)))},u_t^\varepsilon)d\lambda,\Theta_t^2-\Theta_t^*\Bigr\rangle\\
       &+\widehat{E}\[\Bigr\langle\int_0^1\partial_{\mu}f(t,\Theta_t^*+\lambda(\Theta_t^2-\Theta_t^*),
            P_{(X^*_t+\lambda(X^1_t+X^2_t),Y_t^*+\lambda(\overline{\overline{Y}}_t+M(t)))},u_t^\varepsilon;\\
            &\ \ \ \ \ (\widehat{X}^*_t+\lambda(\widehat{X}^1_t+\widehat{X}^2_t),\widehat{Y}_t^*+\lambda(\widehat{\overline{\overline{Y}}}_t+\widehat{M}(t))))d\lambda,
            (\widehat{X}^1_t+\widehat{X}^2_t,\widehat{\overline{\overline{Y}}}_t+\widehat{M}(t))\Bigr\rangle\]\\
      =&\langle f_{(x,y,z)}(t),\Theta_t^2-\Theta_t^*\rangle
        +\widehat{E}\[\langle\widehat{f}_{\mu}(t),(\widehat{X}^1_t+\widehat{X}^2_t,\widehat{\overline{\overline{Y}}}_t+\widehat{M}(t))\rangle\]+I_{2,1}(t)+I_{2,2}(t),
\end{split}
\end{equation}
where
\begin{equation*}
\begin{split}
I_{2,1}(t):=&\Bigr\langle\int_0^1\partial_{(x,y,z)}\(f(t,\Theta_t^*+\lambda(\Theta_t^2-\Theta_t^*),
            P_{(X^*_t+\lambda(X^1_t+X^2_t),Y_t^*+\lambda(\overline{\overline{Y}}_t+M(t)))},u_t^\varepsilon)-f_{(x,y,z)}(t)\)d\lambda,\\
            &\ \ \ \ \ \Theta_t^2-\Theta_t^*\Bigr\rangle,\\
I_{2,2}(t):=&\widehat{E}\[\Bigr\langle\int_0^1\(\partial_{\mu}f(t,\Theta_t^*+\lambda(\Theta_t^2-\Theta_t^*),
            P_{(X^*_t+\lambda(X^1_t+X^2_t),Y_t^*+\lambda(\overline{\overline{Y}}_t+M(t)))},u_t^\varepsilon;\\
            &\ \ \ \ \ \ (\widehat{X}^*_t+\lambda(\widehat{X}^1_t+\widehat{X}^2_t),\widehat{Y}_t^*+\lambda(\widehat{\overline{\overline{Y}}}_t+\widehat{M}(t)))-\widehat{f}_{\mu}(t)\)d\lambda,
            (\widehat{X}^1_t+\widehat{X}^2_t,\widehat{\overline{\overline{Y}}}_t+\widehat{M}(t))\Bigr\rangle\].
\end{split}
\end{equation*}
Now let us consider $I_{2,1}(t)$,
\begin{equation*}
\begin{split}
&I_{2,1}(t)\\
=&\Bigr\langle\!\int_0^1\!\!\(f_{(x,y,z)}(t,\Theta_t^*\!+\!\lambda(\Theta_t^2\!-\!\Theta_t^*),
            P_{(X^*_t+\lambda(X^1_t+X^2_t),Y_t^*+\lambda(\overline{\overline{Y}}_t+M(t)))},u_t)
            \!-\!f_{(x,y,z)}(t)\!\)I_{E_\varepsilon}(t)d\lambda,\Theta_t^2\!-\!\Theta_t^*\Bigr\rangle\\
            +&\Bigr\langle\int_0^1\(f_{(x,y,z)}(t,\Theta_t+\lambda(\Theta_t^3-\Theta_t),
            P_{(X^*_t+\lambda(X^1_t+X^2_t),Y_t^*+\lambda(\overline{\overline{Y}}_t+M(t)))},u_t^*)-f_{(x,y,z)}(t)\)d\lambda,
            \Theta_t^3-\Theta_t\Bigr\rangle\\
            -&\Bigr\langle\!\int_0^1\!\!\(f_{(x,y,z)}(t,\Theta_t\!+\!\lambda(\Theta_t^3\!-\!\Theta_t),
            P_{(X^*_t+\lambda(X^1_t+X^2_t),Y_t^*+\lambda(\overline{\overline{Y}}_t+M(t)))},u_t^*)
            \!-\!f_{(x,y,z)}(t)\!\)I_{E_\varepsilon}(t)d\lambda,
            \Theta_t^3\!-\!\Theta_t\!\Bigr\rangle\\
 =:&I^1_{2,1}(t)+I^2_{2,1}(t)+I^3_{2,1}(t).
\end{split}
\end{equation*}
We begin with the estimate of $I^1_{2,1}(t)$. For this we note that $f_{(x,y,z)}=(f_x,f_y,f_z)$ and that, due to assumption \textbf{(A.4)} and Proposition \ref{prop3.1}, for all $p> 1$,
\begin{equation*}
\begin{split}
&E\[\(\!\!\int_0^T\!\Bigr|\! \int_0^1\!\!(f_{x}(t,\Theta_t^*\!+\!\lambda(\Theta_t^2\!-\!\Theta_t^*),
            P_{(X^*_t+\lambda(X^1_t+X^2_t),Y_t^*+\lambda(\overline{\overline{Y}}_t+M(t)))},u_t)
            \!-\!f_x(t))I_{E_\varepsilon}(t)
            (X^1_t\!+\!X^2_t)d\lambda\Bigr|dt\)^p\]\\
&\leq C\varepsilon^pE\[\sup_{t\in[0,T]}|X^1_t|^p+\sup_{t\in[0,T]}|X^2_t|^p\]
    \leq C(\varepsilon^\frac{3p}{2}+\varepsilon^{2p}).
\end{split}
\end{equation*}
Similarily, from assumption \textbf{(A.4)} and the Proposition \ref{prop3.1} and \ref{prop5.1},
\begin{equation*}
\begin{split}
&E\[\(\!\!\int_0^T\!\Bigr|\! \int_0^1\!\!(f_{y}(t,\Theta_t^*\!+\!\lambda(\Theta_t^2\!-\!\Theta_t^*),
            P_{(X^*_t+\lambda(X^1_t+X^2_t),Y_t^*+\lambda(\overline{\overline{Y}}_t+M(t)))},u_t)
            \!-\!f_y(t))I_{E_\varepsilon}(t)
            (\overline{\overline{Y}}_t\!+\!M(t))d\lambda\Bigr|dt\)^p\]\\
&\leq C\varepsilon^pE\[\sup_{t\in[0,T]}|\overline{\overline{Y}}_t|^p+\sup_{t\in[0,T]}|M(t)|^p\]
 \leq C\varepsilon^\frac{3p}{2},\ \varepsilon>0.
\end{split}
\end{equation*}
Since $|f_z(t,x,y,z,\mu)|\leq L+\gamma|z|$, it holds
\begin{equation}\label{eq5.16}
\begin{split}
\alpha_\varepsilon:=&E\[\(\int_0^T\Bigr| \int_0^1(f_{z}(t,\Theta_t^*+\lambda(\Theta_t^2-\Theta^*_t),
            P_{(X^*_t+\lambda(X^1_t+X^2_t),Y_t^*+\lambda(\overline{\overline{Y}}_t+M(t)))},u_t)-f_z(t))I_{E_\varepsilon}(t)\\
            &\ \ \ \ \ \ \ \ (\overline{\overline{Z}}_t+\beta(t))d\lambda\Bigr|dt\)^p\]\\
    \leq& C_p E\[\(\int_0^T(1+\gamma(|Z_t^*|+|\overline{\overline{Z}}_t+\beta(t)+p_1(t)\delta\sigma(t)I_{E_\varepsilon}(t)|))
        I_{E_\varepsilon}(t)|\overline{\overline{Z}}_t+\beta(t)|dt\)^p\].
\end{split}
\end{equation}
The hardest term in (\ref{eq5.16}) to estimate is $E\[\(\displaystyle\int_0^T|Z_t^*||\overline{\overline{Z}}_t|I_{E_\varepsilon}(t)dt\)^p\]$.
A straight-forward argument gives only
\begin{equation*}
\begin{aligned}
&E\[\(\int_0^T|Z_t^*||\overline{\overline{Z}}_t|I_{E_\varepsilon}(t)dt\)^p\]
\leq E\[\(\int_0^T|Z_t^*|^2I_{E_\varepsilon}(t)dt\)^\frac{p}{2}
        \(\int_0^T|\overline{\overline{Z}}_t|^2dt\)^\frac{p}{2}\]\\
&\leq \(E\[\(\int_0^T|Z_t^*|^2I_{E_\varepsilon}(t)dt\)^p\]\)^\frac{1}{2}
            \(E\[\(\int_0^T|\overline{\overline{Z}}_t|^2dt\)^p\]\)^\frac{1}{2}
\leq \rho^1(\varepsilon)\cdot\rho^2(\varepsilon)\varepsilon^\frac{p}{2}
=:\rho(\varepsilon)\varepsilon^\frac{p}{2},
\end{aligned}
\end{equation*}
but we need an estimate which leads to a convergence speed quicker than $\varepsilon^p$, as $\varepsilon\downarrow0$.

In order to overcome the constraints related with (3.41) in \cite{HJX2022} (see our explanations in the introduction), we give a simpler and more clear approach.
We introduce the deterministic set  $\Gamma_M:=\{t\in[0,T]:\ E[|Z_t^*|^2]\leq M\},\ M\geq 1.$ Note that, as $(X_t^*,Y_t^*,Z_t^*)$ does not depend on $\varepsilon>0$ nor on $E_\varepsilon$, neither $\Gamma_M$ depends on $\varepsilon>0$. Moreover, as
$$\int_0^TE[|Z_t^*|^2]dt=E\[\int_0^T|Z_t^*|^2dt\]<\infty,$$
we have $I_{\bigcup_{\substack{M\geq 1}}\Gamma_M}(t)=1$, dt-a.e.

Let $M\geq 1$ be arbitrarily fixed, $E_\varepsilon:= E_\varepsilon\cap\Gamma_M\subset\Gamma_M$. Observe that also for the redefined $E_\varepsilon$ we have $|E_\varepsilon|\leq\varepsilon$. For all $1< p<2$, from H\"{o}lder's inequality,
\begin{equation*}
\begin{aligned}
E\[\(\int_0^T|Z_t^*||\overline{\overline{Z}}_t|I_{E_\varepsilon}(t)dt\)^p\]
\leq &E\[\(\int_{E_\varepsilon}|Z_t^*|^2dt\)^\frac{p}{2}
        \(\int_0^T|\overline{\overline{Z}}_t|^2dt\)^\frac{p}{2}\]\\
\leq &\(E\[\int_{E_\varepsilon}|Z_t^*|^2dt\]\)^\frac{p}{2}
        \(E\[\(\int_0^T|\overline{\overline{Z}}_t|^2dt\)^\frac{p}{2-p}\]\)^\frac{2-p}{2}\\
= & \(\int_{E_\varepsilon}E[|Z_t^*|^2]dt\)^\frac{p}{2}
        \(E\[\(\int_0^T|\overline{\overline{Z}}_t|^2dt\)^\frac{p}{2-p}\]\)^\frac{2-p}{2}.
\end{aligned}
\end{equation*}
Here $\(\displaystyle\int_{E_\varepsilon}E[|Z_t^*|^2]dt\)^\frac{p}{2}\leq M^\frac{p}{2}\varepsilon^\frac{p}{2}$. Put $\beta:=\frac{4p}{2-p}$. For $1<p<2$, we have $\beta>1$. We choose $\kappa,\ \kappa_0$ such that $\beta\vee p_{f_z}^*<\kappa<\kappa_0$. Then due to Proposition \ref{prop5.1},
\begin{equation*}
\begin{split}
&\(E\[\(\int_0^T|\overline{\overline{Z}}_t|^2dt\)^\frac{p}{2-p}\]\)^\frac{2-p}{2}
\leq \(E\[\(\int_0^T|\overline{\overline{Z}}_t|^2dt\)^\frac{\beta}{2}\]\)^\frac{p}{\beta}
\leq \(E\[\(\int_0^T|\overline{\overline{Z}}_t|^2dt\)^\frac{\kappa}{2}\]\)^\frac{p}{\kappa}\\
&\leq \(\widetilde{\rho}(\varepsilon)\varepsilon^\frac{\kappa}{2}\)^\frac{p}{\kappa}
=:\rho(\varepsilon)\varepsilon^\frac{p}{2},
\end{split}
\end{equation*}
for some $\rho:\mathbb{R}_+\rightarrow\mathbb{R}_+$ with $\rho(\varepsilon)\downarrow0$, as $\varepsilon\downarrow0$. Hence, it follows that for all $1<p<2$,
\begin{equation}\label{eq5.17}
E\[\(\int_0^T|Z_t^*||\overline{\overline{Z}}_t|I_{E_\varepsilon}(t)dt\)^p\]
\leq M^\frac{p}{2}\rho(\varepsilon)\varepsilon^p,\ \varepsilon>0.
\end{equation}
The remaining terms in (\ref{eq5.16}) can be estimated in a more classical way. Indeed,
\begin{equation*}
\begin{aligned}
\alpha_\varepsilon
\leq &C_p E\[\(\int_{E_\varepsilon}|Z_t^*||\overline{\overline{Z}}_t|dt\)^p\]
        +C_p E\[\(\int_{E_\varepsilon}(1+|\overline{\overline{Z}}_t+\beta(t)+p_1(t)\delta\sigma(t)I_{E_\varepsilon}(t)|)
            |\overline{\overline{Z}}_t+\beta(t)|dt\)^p\]\\
        &+C_pE\[\(\int_{E_\varepsilon}|Z_t^*||\beta(t)|dt\)^p\].
\end{aligned}
\end{equation*}
Thanks to Proposition \ref{prop3.1}, for all $p>1$,
$$E\[\(\int_{E_\varepsilon}|\beta(t)|^2dt\)^p\]\leq \rho_p(\varepsilon)\varepsilon^p,\ \varepsilon>0,$$
with $\rho_p(\varepsilon)\downarrow0$, as $\varepsilon\downarrow0$. Combining this with Proposition \ref{prop5.1} and the argument for (\ref{eq5.17}), we get for all $1< p<2$,
$$\alpha_\varepsilon\leq C_M\rho(\varepsilon)\varepsilon^p,$$
for some function $\rho:\mathbb{R}_+\rightarrow\mathbb{R}_+$ with $\rho(\varepsilon)\downarrow0$, as $\varepsilon\downarrow0$.
Hence, summarising the above estimates, we have
$$E\[\(\int_0^T|I^1_{2,1}(t)|dt\)^p\]\leq C\varepsilon^p\rho(\varepsilon).$$
With the same arguments we also get
$$E\[\(\int_0^T|I^3_{2,1}(t)|dt\)^p\]\leq C\varepsilon^p\rho(\varepsilon).$$
We notice that
\begin{equation*}
\begin{split}
&f_{(x,y,z)}(t,\Theta_t+\lambda(\Theta_t^3-\Theta_t),
            P_{(X^*_t+\lambda(X^1_t+X^2_t),Y_t^*+\lambda(\overline{\overline{Y}}_t+M(t)))},u_t^*)-f_{(x,y,z)}(t)\\
=&\lambda\int_0^1\partial^2_{(x,y,z)}f(t,\Theta_t+\theta\lambda(\Theta_t^3-\Theta_t),
            P_{(X^*_t+\theta\lambda(X^1_t+X^2_t),Y_t^*+\theta\lambda(\overline{\overline{Y}}_t+M(t)))},u_t^*)(\Theta_t^3-\Theta_t)d\theta\\
            &+\widehat{E}\[\Bigr\langle\lambda\int_0^1\partial_{\mu}f_{(x,y,z)}(t,\Theta_t+\theta\lambda(\Theta_t^3-\Theta_t),
            P_{(X^*_t+\theta\lambda(X^1_t+X^2_t),Y_t^*+\theta\lambda(\overline{\overline{Y}}_t+M(t)))},u_t^*;\\
                & \ \ \ \ \ \ \ \ \ (\widehat{X}_t^*+\theta\lambda(\widehat{X}_t^1+\widehat{X}_t^2),\widehat{Y}_t^*+\theta\lambda(\widehat{\overline{\overline{Y}}}_t+\widehat{M}(t))))d\theta,
                (\widehat{X}_t^1+\widehat{X}_t^2,\widehat{\overline{\overline{Y}}}_t+\widehat{M}(t))\Bigr\rangle\]\\
=:&\mathbb{I}_1(t)+\mathbb{I}_2(t),\ t\in[0,T].
\end{split}
\end{equation*}
We have
\begin{equation*}
\begin{split}
&\Bigr\langle\!\int_0^1\!\lambda\!\int_0^1\!\partial^2_{(x,y,z)}f(t,\Theta_t\!+\!\theta\lambda(\Theta_t^3\!-\!\Theta_t),
            P_{(X^*_t+\theta\lambda(X^1_t+X^2_t),Y_t^*+\theta\lambda(\overline{\overline{Y}}_t+M(t)))},u_t^*)
            (\Theta_t^3\!-\!\Theta_t)d\theta d\lambda, \Theta_t^3\!-\!\Theta_t\!\Bigr\rangle\\
=&\Bigr\langle\frac{1}{2}\partial^2_{(x,y,z)}f(t)(\Theta_t^3-\Theta_t),\Theta_t^3-\Theta_t\Bigr\rangle
        +\mathbb{I}_1^1(t),\ t\in[0,T],
\end{split}
\end{equation*}
where\  \ $\mathbb{I}_1^1(t)\!:=\Bigr\langle\!\displaystyle\int_0^1\!\!\lambda\!\!\displaystyle\int_0^1\!\!\(\partial^2_{(x,y,z)}f(t,\Theta_t+\theta\lambda(\Theta_t^3-\Theta_t),
            P_{(X^*_t+\theta\lambda(X^1_t+X^2_t),Y_t^*+\theta\lambda(\overline{\overline{Y}}_t+M(t)))},u_t^*)-\partial^2_{(x,y,z)}f(t)\!\)$
            $(\Theta_t^3-\Theta_t)d\theta d\lambda, \Theta_t^3-\Theta_t\Bigr\rangle$.

Let us remark that thanks to the Propositions \ref{prop3.1}, \ref{lem4.1} and \ref{prop5.1},
$$E\[\(\int_0^T|\Theta_t^3-\Theta_t|^2dt\)^p\]\leq C\varepsilon^p.$$
As $p_2\in \mathcal{S}^p([0,T];\mathbb{R})$ and $\beta_2\in\mathcal{H}^{2,p}(0,T;\mathbb{R})$, $p> 1$, Proposition \ref{lem4.1} allows to conclude that
$$E\[\(\int_0^T|\widehat{E}[p_2(t)\widehat{X}_t^1]+\widehat{E}[\beta_3(t)\widehat{X}_t^1]|^2dt\)^p\]\leq C\varepsilon^p\rho(\varepsilon).$$
Hence, recalling that $\Theta_t^3-\Theta_t=(X_t^1+X_t^2,\overline{\overline{Y}}_t+M(t),\overline{\overline{Z}}_t+\beta(t))$, we see that $(\Theta_t^3-\Theta_t)-(X_t^1,p_1(t)X_t^1,\beta_2(t)X_t^1)
=(X_t^2,\overline{\overline{Y}}_t+M(t)-p_1(t)X_t^1,\overline{\overline{Z}}_t+\beta(t)-\beta_2(t)X_t^1)$,  and from the Propositions \ref{prop3.1}, \ref{lem4.1} and \ref{prop5.1} we get
\begin{equation*}
\begin{split}
E\[\(\int_0^T|(\Theta_t^3-\Theta_t)-(X_t^1,p_1(t)X_t^1,\beta_2(t)X_t^1)|^2dt\)^p\]\leq C\varepsilon^p\rho(\varepsilon).
\end{split}
\end{equation*}
Thus, using the boundedness and the Lipschitz property of $\partial^2_{(x,y,z)}f$, we have
\begin{equation*}
\begin{split}
E\[\(\int_0^T|\mathbb{I}_1^1(t)|dt\)^p\]
\leq CE\[\(\int_0^T(|\Theta_t^3-\Theta_t|^3\wedge|\Theta_t^3-\Theta_t|^2)dt\)^p\]
\leq C\varepsilon^p\rho(\varepsilon),
\end{split}
\end{equation*}
where $a\wedge b=\min(a,b),\ a, b\in\mathbb{R}$.
Indeed,
\begin{equation*}
\begin{aligned}
|\mathbb{I}_1^1(t)|
\leq & |\Theta_t^3-\Theta_t|^2\wedge|\Theta_t^3-\Theta_t|^3\\
\leq &C|(\Theta_t^3-\Theta_t)-(X_t^1,p_1(t)X_t^1,\beta_2(t)X_t^1)|^2+C(|X_t^1|^3+|p_1(t)X_t^1|^3)\\
        &+|\beta_2(t)X_t^1|^2\wedge|\beta_2(t)X_t^1|^3.
\end{aligned}
\end{equation*}
The estimate of $E\[\(\displaystyle\int_0^T(|X_t^1|^3+|p_1(t)X_t^1|^3)dt\)^p\]\leq C \varepsilon^\frac{3p}{2},\ \varepsilon>0$, is a direct consequence of Proposition \ref{prop3.1}. Concerning the last term in the above inequality we note that we only have $\beta_2\in\mathcal{H}^{2,q}(0,T;\mathbb{R})$, $q>1$. So we proceed as follows:
$$|\beta_2(t)X_t^1|^2\wedge|\beta_2(t)X_t^1|^3
\leq \sup_{t\in[0,T]}|X_t^1|^2(|\beta_2(t)|^2(1\wedge|\beta_2(t)X_t^1|)),\ t\in[0,T],$$
and so,
\begin{equation*}
\begin{aligned}
&E\[\(\int_0^T(|\beta_2(t)X_t^1|^2\wedge|\beta_2(t)X_t^1|^3)dt\)^p\]\\
\leq & \(E\[\sup_{t\in[0,T]}|X_t^1|^{4p}\]\)^\frac{1}{2}\(E\[\(\int_0^T|\beta_2(t)|^2(1\wedge|\beta_2(t)X_t^1|)dt\)^{2p}\]\)^\frac{1}{2}
\leq \rho(\varepsilon)\varepsilon^p,\ \varepsilon>0,
\end{aligned}
\end{equation*}
where, due to Proposition \ref{prop3.1} and the dominated convergence theorem,
$$\rho(\varepsilon)
=\(E\[\(\int_0^T|\beta_2(t)|^2(1\wedge|\beta_2(t)X_t^1|)dt\)^{2p}\]\)^\frac{1}{2}
\rightarrow 0,\ \text{as}\ \varepsilon\downarrow0.$$
To simplify the notation, we put $\Theta_t^4:=(X_t^1,p_1(t)X_t^1,\beta_2(t)X_t^1)$ and $\Theta_t^{\prime4}:=(X_t^1,p_1(t)X_t^1)$.
Then, from the above estimates we obtain
\begin{equation*}
\begin{split}
&\Bigr\langle\!\int_0^1\!\lambda\!\int_0^1\!\partial^2_{(x,y,z)}f(t,\Theta_t\!+\!\theta\lambda(\Theta_t^3\!-\!\Theta_t),
            P_{(X^*_t+\theta\lambda(X^1_t+X^2_t),Y_t^*+\theta\lambda(\overline{\overline{Y}}_t+M(t)))},u_t^*)
            (\Theta_t^3\!-\!\Theta_t)d\theta d\lambda, \Theta_t^3\!-\!\Theta_t\Bigr\rangle\\
=&\Bigr\langle\frac{1}{2}\partial^2_{(x,y,z)}f(t)\Theta_t^4,\Theta_t^4\Bigr\rangle+\mathbb{I}^1_1(t)+\mathbb{I}_1^2(t),
\end{split}
\end{equation*}
where
$$\mathbb{I}_1^2(t):=\Bigr\langle\frac{1}{2}\partial^2_{(x,y,z)}f(t)(\Theta_t^3-\Theta_t),\Theta_t^3-\Theta_t\Bigr\rangle
-\Bigr\langle\frac{1}{2}\partial^2_{(x,y,z)}f(t)\Theta_t^4,\Theta_t^4\Bigr\rangle,$$
and
$$E\[\(\int_0^T|\mathbb{I}_1^1(t)+\mathbb{I}_1^2(t)|dt\)^p\]\leq  C\varepsilon^p\rho(\varepsilon).$$
Using Proposition \ref{lem4.1} again, we show with estimates similar to those above that
$$E\[\(\int_0^T\Bigr|\int_0^1\mathbb{I}_2(t)d\lambda\Bigr||\Theta_t^3-\Theta_t|dt\)^p\]\leq  C\varepsilon^p\rho(\varepsilon).$$
Combining the above estimates we see that
\begin{equation*}
I_{2,1}^2(t)=\Bigr\langle\frac{1}{2}\partial^2_{(x,y,z)}f(t)\Theta_t^4,\Theta_t^4\Bigr\rangle+\mathbb{I}_3(t),
\end{equation*}
where
$\mathbb{I}_3(t):=\mathbb{I}^1_1(t)+\mathbb{I}_1^2(t)+\Bigr\langle\displaystyle\int_0^1\mathbb{I}_2(t)d\lambda,\Theta_t^3-\Theta_t\Bigr\rangle$
 and $E\[\(\displaystyle\int_0^T|\mathbb{I}_3(t)|dt\)^p\]\leq  C\varepsilon^p\rho(\varepsilon)$.

Next we estimate $I_{2,2}(t)$ with computations similar to those for  $I_{2,1}(t)$. So we obtain
\begin{equation*}
I_{2,2}(t)=\widehat{E}\[\Bigr\langle\frac{1}{2}\partial_{(\widehat{x},\widehat{y})}\partial_\mu f(t)\widehat{\Theta}_t^{\prime4},
    \widehat{\Theta}_t^{\prime4}\Bigr\rangle\]+\mathbb{I}_4(t),\ t\in[0,T],
\end{equation*}
where $\mathbb{I}_4(t):=I_{2,2}(t)-\widehat{E}\[\Bigr\langle\frac{1}{2}\partial_{(\widehat{x},\widehat{y})}\partial_\mu f(t)\widehat{\Theta}_t^{\prime4},
    \widehat{\Theta}_t^{\prime4}\Bigr\rangle\]$ with
$E\[\(\displaystyle\int_0^T|\mathbb{I}_4(t)|dt\)^p\]\leq  C\varepsilon^p\rho(\varepsilon).$\\
Combining the above results, we have
\begin{equation}\label{eq5.9}
\begin{split}
I_{2}(t)=&\langle f_{(x,y,z)}(t),\Theta_t^2-\Theta_t^*\rangle
        +\widehat{E}\[\Bigr\langle\widehat{f}_{\mu}(t),(\widehat{X}_t^1+\widehat{X}_t^2,\widehat{\overline{\overline{Y}}}_t+\widehat{M}(t))\Bigr\rangle\]
        +\Bigr\langle\frac{1}{2}\partial^2_{(x,y,z)}f(t)\Theta_t^4,\Theta_t^4\Bigr\rangle\\
        &+\widehat{E}\[\Bigr\langle\frac{1}{2}\partial_{(\widehat{x},\widehat{y})}\partial_\mu f(t)\widehat{\Theta}_t^{\prime4},\widehat{\Theta}_t^{\prime4}\Bigr\rangle\]+\mathbb{I}_5(t),
\end{split}
\end{equation}
where $\mathbb{I}_5(t):=I_{2,1}^1(t)+I_{2,1}^3(t)+\mathbb{I}_3(t)+\mathbb{I}_4(t)$ and
$E\[\(\displaystyle\int_0^T|\mathbb{I}_5(t)|dt\)^p\]\leq  C\varepsilon^p\rho(\varepsilon).$\\
Finally, we consider $I_3(t)$. Recall its definition in (\ref{eq5.18}) and also (\ref{eq5.19}) for $\check{Z}^\varepsilon$. Obviously,
\begin{equation*}
\begin{split}
|I_3(t)|
\leq &C\(|X^\varepsilon_t\!-\!X^*_t\!-\!X^1_t\!-\!X^2_t|
    +|Y^\varepsilon_t\!-\!\overline{Y}_t^\varepsilon\!-\!M(t)|
    +(1\!+\!|Z^\varepsilon_t|\!+\!|\overline{Z}_t^\varepsilon\!+\!\beta(t)\!+\!p_1(t)\delta\sigma(t)I_{E_\varepsilon}(t)|)|\check{Z}^\varepsilon_t|\\
        &+(E[|X^\varepsilon_t-X^*_t-X^1_t-X^2_t|^2])^\frac{1}{2}+(E[|Y^\varepsilon_t-\overline{Y}_t^\varepsilon-M(t)|^2])^\frac{1}{2}\)\\
        \leq &C\(|X^\varepsilon_t-X^*_t-X^1_t-X^2_t|+|\check{Y}^\varepsilon_t|
    +(1+|Z^\varepsilon_t|+|\check{Z}^\varepsilon_t|)|\check{Z}^\varepsilon_t|+(E[|X^\varepsilon_t-X^*_t-X^1_t-X^2_t|^2])^\frac{1}{2}\\
    &+(E[|\check{Y}^\varepsilon_t|^2])^\frac{1}{2}\).
\end{split}
\end{equation*}
We can easily see that, thanks to Proposition \ref{prop3.1}, for all $p>1$,
$$E\[\sup_{t\in[0,T]}|X^\varepsilon_t-X^*_t-X^1_t-X^2_t|^p\]=o(\varepsilon^p),$$
due to Proposition \ref{prop5.2}, for all $1< p<2$,
$$E[\sup_{t\in[0,T]}|\check{Y}^\varepsilon_t|^p]
    \leq \(E[\sup_{t\in[0,T]}|\check{Y}^\varepsilon_t|^2]\)^\frac{p}{2}\leq C\varepsilon^p\rho^\frac{p}{2}(\varepsilon),$$
and
$$E\[\(\int_0^T|Z^\varepsilon_t||\check{Z}^\varepsilon_t|dt\)^p\]
\leq \(E\[\(\int_0^T|Z^\varepsilon_t|^2dt\)^\frac{p}{2-p}\]\)^\frac{2-p}{2}
\(E\[\int_0^T|\check{Z}^\varepsilon_t|^2dt\]\)^\frac{p}{2}
\leq C\varepsilon^p\rho^\frac{p}{2}(\varepsilon).$$
Therefore, we have
\begin{equation}\label{eq5.10}
\begin{split}
&E\[\(\int_0^T|I_{3}(t)|dt\)^p\]
 \leq C\varepsilon^p\rho(\varepsilon).
\end{split}
\end{equation}
Combined with (\ref{eq5.18}), (\ref{eq5.7}), (\ref{eq5.9}) and (\ref{eq5.10}), this yields
\begin{equation*}
\begin{split}
&f(t,X^\varepsilon_t,Y^\varepsilon_t,Z^\varepsilon_t,P_{(X^\varepsilon_t,Y^\varepsilon_t)},u_t^\varepsilon)
    -f(t,X^*_t,Y^*_t,Z^*_t,P_{(X^*_t,Y^*_t)},u_t^*)\\
=&\Delta f(t)I_{E_\varepsilon}(t)
+\langle A_1(t),\Theta_t^1-\Theta_t^2\rangle+\widehat{E}[\widehat{A}_1(t)\cdot\widehat{\delta\overline{Y}}_t]
    +\langle f_{(x,y,z)}(t),\Theta_t^2-\Theta_t^*\rangle\\
        &\!+\!\widehat{E}\[\langle\widehat{f}_{\mu}(t),(\widehat{X}_t^1\!+\!\widehat{X}_t^2,\widehat{\overline{\overline{Y}}}_t\!+\!\widehat{M}(t))\rangle\]
        \!+\!\frac{1}{2}\Bigr\langle\partial^2_{(x,y,z)}f(t)\Theta_t^4,\Theta_t^4\Bigr\rangle
        \!+\!\frac{1}{2}\widehat{E}\[\Bigr\langle\partial_{(\widehat{x},\widehat{y})}\partial_\mu f(t)\widehat{\Theta}_t^{\prime4},\widehat{\Theta}^{\prime4}_t\Bigr\rangle\]\!+\!R_\varepsilon(t),
\end{split}
\end{equation*}
with $R_\varepsilon(t):=\mathbb{I}_5(t)+I_3(t)$ and
$E\[\(\displaystyle\int_0^T|R_\varepsilon(t)|dt\)^p\]\leq  C\varepsilon^p\rho(\varepsilon)$.\\
Thus, (\ref{eq5.11}) can be rewritten as
\begin{equation}\label{eq5.12}
\left\{
\begin{split}
d(\delta \overline{Y}_t)
=&-\Bigr\{\langle A_1(t),\Theta_t^1-\Theta_t^2\rangle+\widehat{E}[\widehat{A}_1(t)\cdot\widehat{\delta\overline{Y}}_t]
    +\langle f_{(x,y,z)}(t),(X_t^1+X_t^2,M(t),\beta(t))\rangle\\
        &+\widehat{E}\[\langle\widehat{f}_{\mu}(t),(\widehat{X}_t^1+\widehat{X}_t^2,\widehat{M}(t))\rangle\]
        +\frac{1}{2}\Bigr\langle\partial^2_{(x,y,z)}f(t)\Theta_t^4,\Theta_t^4\Bigr\rangle
        +\frac{1}{2}\widehat{E}\[\Bigr\langle\partial_{(\widehat{x},\widehat{y})}\partial_\mu f(t)\widehat{\Theta}_t^{\prime4},\widehat{\Theta}_t^{\prime4}\Bigr\rangle\]\\
    &+\alpha_{2}(t)(X_t^1+X_t^2)+\widetilde{E}[\alpha_{3}(t)(\widetilde{X}_t^1+\widetilde{X}_t^2)]+\frac{1}{2}\alpha_{4}(t)(X_t^1)^2
    +\frac{1}{2}\widetilde{E}[\alpha_{5}(t)(\widetilde{X}_t^1)^2]+R_\varepsilon(t)\Bigr\}dt\\
    &+\delta \overline{Z}_t dB_t,\ t\in[0,T],\\
\delta \overline{Y}_T=&0.
\end{split}
\right.
\end{equation}
Using the adjoint equations (\ref{eq3.11})-(\ref{eq3.14}), we have
\begin{equation*}
\begin{split}
\alpha_2(t)=&-f_x(t)-p_1(t)f_y(t)-p_1(t)\sigma_x(t)f_z(t)-q_{1,1}(t)f_z(t),\\
\alpha_3(t)=&-\widetilde{f}_{\mu,x}(t)-p_1(t)f_z(t)\widetilde{\sigma}_\mu(t)-\widetilde{p}_1(t)\widetilde{f}_{\mu,y}(t)
            -p_{2}(t)f_y(t)-\widehat{E}[\widehat{p}_{2}(t)\widehat{f}_{\mu,y}(t)]-q_{2,1}f_z(t),\\
\alpha_4(t)=&-p_1(t)f_z(t)\sigma_{xx}(t)-P_1(t)f_y(t)-2P_1(t)f_z(t)\sigma_x(t)-Q_{1,1}(t)f_z(t)\\
            &-\langle\partial^2_{(x,y,z)}f(t)(1,p_1(t),p_1(t)\sigma_x(t)+q_{1,1}(t))^\top,(1,p_1(t),p_1(t)\sigma_x(t)+q_{1,1}(t))^\top\rangle,\\
\alpha_5(t)=&-\widetilde{P}_1(t)\widetilde{f}_{\mu,y}(t)-f_y(t)P_2(t)-\widehat{E}[\widehat{P}_2(t)\widehat{f}_{\mu,y}(t)]
            -p_1(t)f_z(t)\widetilde{\sigma}_{\widetilde{x}\mu}(t)-f_z(t)Q_{2,1}(t)\\
            &-\langle\partial_{(\widetilde{x},\widetilde{y})}\widetilde{f}_{\mu}(t)(1,p_1(t))^\top,(1,p_1(t))^\top\rangle,\\
\beta(t)=&q_{1,1}(t)(X_t^1+X_t^2)+p_1(t)\sigma_x(t)(X_t^1+X_t^2)+\widetilde{E}[p_1(t)\widetilde{\sigma}_\mu(t)(\widetilde{X}_t^1+\widetilde{X}_t^2)]
            +\widetilde{E}[q_{2,1}(t)(\widetilde{X}_t^1+\widetilde{X}_t^2)]\\
            &+\frac{1}{2}p_1(t)\sigma_{xx}(t)(X_t^1)^2+\frac{1}{2}Q_{1,1}(t)(X_t^1)^2+P_1(t)\sigma_{x}(t)(X_t^1)^2
            +\frac{1}{2}\widetilde{E}[p_1(t)\widetilde{\sigma}_{\widetilde{x}\mu}(t)(\widetilde{X}_t^1)^2]\\
            &+\frac{1}{2}\widetilde{E}[Q_{2,1}(t)(\widetilde{X}_t^1)^2]+p_1(t)\delta\sigma_x(t)X_t^1I_{E_\varepsilon}(t)
            +P_1(t)\delta\sigma(t)X_t^1I_{E_\varepsilon}(t).
\end{split}
\end{equation*}
We should notice that both $(\widetilde{\Omega},\widetilde{\mathcal{F}},\widetilde{P})$ and $(\widehat{\Omega},\widehat{\mathcal{F}},\widehat{P})$ are independent copy spaces of $(\Omega,\mathcal{F},P)$, and so
$$\widetilde{E}[\widetilde{f}_{\mu,x}(t)(\widetilde{X}_t^1+\widetilde{X}_t^2)]
=\widehat{E}[\widehat{f}_{\mu,x}(t)(\widehat{X}_t^1+\widehat{X}_t^2)]$$
and
$$\widetilde{E}[\widetilde{p}_1\widetilde{f}_{\mu,y}(t)(\widetilde{X}_t^1+\widetilde{X}_t^2)]
=\widehat{E}[\widehat{p}_1\widehat{f}_{\mu,y}(t)(\widehat{X}_t^1+\widehat{X}_t^2)].$$
Moreover, in analogy to (\ref{relation}) , we get
\begin{equation*}
\begin{split}
&\widehat{E}\[\widehat{f}_{\mu,y}(t)\widehat{\widetilde{E}[p_2(t)(\widetilde{X}_t^1+\widetilde{X}_t^2)]}\]
=\widetilde{E}\[\widehat{E}[\widehat{f}_{\mu,y}(t)\widehat{p}_2(t)](\widetilde{X}_t^1+\widetilde{X}_t^2)]\]
\end{split}
\end{equation*}
and
$$\widehat{E}\[\widehat{f}_{\mu,y}(t)\widehat{\widetilde{E}[P_2(t)(\widetilde{X}_t^1)^2]}\]
=\widetilde{E}\[\widehat{E}[\widehat{f}_{\mu,y}(t)\widehat{P}_2(t)](\widetilde{X}_t^1)^2]\].$$
Hence, substituting the above relations for $\alpha_j(t),\ j=2,3,4,5$ and $\beta(t)$, (\ref{eq5.12}) reduces to
\begin{equation}\label{eq5.13}
\left\{
\begin{aligned}
d(\delta \overline{Y}_t)
=&-\Bigr\{\langle A_1(t),(0,\delta \overline{Y}_t,\delta\overline{Z}_t)\rangle
        +\widehat{E}[\widehat{A}_1(t)\cdot\widehat{\delta\overline{Y}}_t]
   +f_z(t)p_1(t)\delta\sigma_x(t)X_t^1I_{E_\varepsilon}(t)\\
    &\ \ \ \ +f_z(t)P_1(t)\delta\sigma(t)X_t^1I_{E_\varepsilon}(t)+R_\varepsilon(t)\Bigr\}dt+\delta \overline{Z}_t dB_t,\ t\in[0,T],\\
\delta \overline{Y}_T=&0.
\end{aligned}
\right.
\end{equation}
Recall that $|f_z(t)|\leq L+\gamma|Z_t^*|,\ t\in[0,T]$. Hence, as $E_\varepsilon\subset\Gamma_M$, for $1< p<2$,
\begin{equation*}
\begin{split}
&E\[\(\int_0^T|f_z(t)||(p_1(t)\delta\sigma_x(t)+P_1(t)\delta\sigma(t))X_t^1|I_{E_\varepsilon}(t)dt\)^p\]\\
\leq &E\[\(\int_0^T|f_z(t)|I_{E_\varepsilon}(t)dt
        \sup_{t\in[0,T]}|p_1(t)\delta\sigma_x(t)+P_1(t)\delta\sigma(t)|\sup_{t\in[0,T]}|X_t^1|\)^p\]\\
\leq &\varepsilon^\frac{p}{2}\(\int_0^TE[|f_z(t)|^2] I_{E_\varepsilon}(t)dt\)^\frac{p}{2}
        \(E\[\sup_{t\in[0,T]}|X_t^1|^\frac{4p}{2-p}\]\)^\frac{2-p}{4}
        \(E\[\sup_{t\in[0,T]}|p_1(t)\delta\sigma_x(t)+P_1(t)\delta\sigma(t)|^\frac{4p}{2-p}\]\)^\frac{2-p}{4}\\
\leq &C_M\varepsilon^p\cdot C\varepsilon^\frac{p}{2}=C_M\varepsilon^{\frac{3p}{2}}.
\end{split}
\end{equation*}
For the sake of simplicity, we put $R^\prime_\varepsilon(t):=f_z(t)p_1(t)\delta\sigma_x(t)X_t^1I_{E_\varepsilon}(t)+f_z(t)P_1(t)\delta\sigma(t)X_t^1I_{E_\varepsilon}(t)
+R_\varepsilon(t)$, $t\in[0,T]$.

This latter estimate and the fact that $A_1(\cdot)$ and $\widehat{A}_1(\cdot)$ are bounded (and so, in particular, $A_1(\cdot)\cdot B\in BMO_2$ ) allow to use Proposition \ref{prop4.1} and so we can conclude:
\begin{equation*}
E\[\sup_{t\in[0,T]}|\delta\overline{Y}_t|^p+\(\int_0^T|\delta\overline{Z}_t|^2dt\)^{\frac{p}{2}}\]
\leq C\varepsilon^p\rho(\varepsilon).
\end{equation*}
The proof is complete now.
\end{proof}
Now the main result of this section is given below.
\begin{theorem}\label{thm5.1}
Suppose the assumptions \textbf{(A1)}-\textbf{(A4)} hold true. Then we have, for all $1<p<2$,
\begin{equation}\label{eq5.14}
E\[\sup_{t\in[0,T]}|Y_t^\varepsilon-Y_t^*-Y_t^1-Y_t^2-\overline{\overline{Y}}_t|^p\]\leq C_p\varepsilon^p\rho(\varepsilon),
\end{equation}
where
$Y_t^1:=p_1(t)X_t^1+\widetilde{E}[p_2(t)\widetilde{X}_t^1]$,
$Y_t^2:=p_1(t)X_t^2+\widetilde{E}[p_2(t)\widetilde{X}_t^2]+\frac{1}{2}P_1(t)(X_t^1)^2+\frac{1}{2}\widetilde{E}[P_2(t)(\widetilde{X}_t^1)^2]$,
and $(\overline{\overline{Y}},\overline{\overline{Z}})$ is the solution of (\ref{eq5.1}).
\end{theorem}
\begin{proof}
We notice that $M(t)=Y_t^1+Y_t^2,\ t\in[0,T]$, and that $Y_t^\varepsilon-Y_t^*-Y_t^1-Y_t^2-\overline{\overline{Y}}_t=\delta\overline{Y}_t+\check{Y}^\varepsilon_t$;
see (\ref{eq5.15}) and the definition of $\delta\overline{Y}$. Thus, according to the Propositions \ref{prop5.2} and \ref{prop5.3}, we have
\begin{equation*}
\begin{split}
&E\[\sup_{t\in[0,T]}|Y_t^\varepsilon-Y_t^*-Y_t^1-Y_t^2-\overline{\overline{Y}}_t|^p\]\\
&\leq CE\[\sup_{t\in[0,T]}|\delta\overline{Y}_t|^p\]+C\(E\[\sup_{t\in[0,T]}|\check{Y}^\varepsilon_t|^2\]\)^\frac{p}{2}
\leq C\varepsilon^p\rho(\varepsilon).
\end{split}
\end{equation*}
\end{proof}

\section{Stochastic Maximum Principle}
In this section, we give the maximum principle. Combining the initial values of (\ref{eq3.2}) and (\ref{eq3.3}), i.e., $X_0^1=0$ and $X_0^2=0$,  we conclude that $Y_0^1=0$ and $Y_0^2=0$. Again, by the result of Theorem \ref{thm5.1}, we deduce that
\begin{equation}\label{eq6.1}
Y_0^\varepsilon
=Y_0^*+Y_0^1+Y_0^2+\overline{\overline{Y}}_0+o(\varepsilon)
=Y_0^*+\overline{\overline{Y}}_0+o(\varepsilon).
\end{equation}
According to (\ref{eq6.1}), we have
\begin{equation}\label{eq6.2}
0\leq J(u^\varepsilon)-J(u^*)=Y_0^\varepsilon-Y_0^*=\overline{\overline{Y}}_0+o(\varepsilon).
\end{equation}
We shall introduce the following linear mean-field SDE:
\begin{equation}\label{eq6.3}
\left\{
\begin{aligned}
d\Gamma_t=&\(f_y(t)\Gamma_t+\widehat{E}[\widehat{f}^*_{\mu,y}(t)\widehat{\Gamma}_t]\)dt+f_z(t)\Gamma_tdB_t,\\
\Gamma_0=&1.
\end{aligned}
\right.
\end{equation}
We recall the definition of $\alpha_1(t)$,
$$\alpha_1(t)=p_1(t)\delta b(t)+q_{1,1}(t)\delta \sigma(t)+\frac{1}{2}P_1(t)(\delta \sigma(t))^2
+\widetilde{E}\[p_2(t)\delta\widetilde{b}(t)+q_{2,2}\delta\widetilde{\sigma}(t)+\frac{1}{2}P_2(t)(\delta\widetilde{\sigma}(t))^2\],$$ and BSDE (\ref{eq5.1}) satisfied by $(\overline{\overline{Y}},\overline{\overline{Z}})$.
By applying It\^{o}'s formula to $\Gamma_t\overline{\overline{Y}}_t$ and taking expectation, we deduce
\begin{equation}\label{eq6.4}
\overline{\overline{Y}}_0=\overline{\overline{Y}}_0\Gamma_0=E\[\int_0^T\Gamma_t(\alpha_1(t)+\Delta f(t))I_{E_\varepsilon}(t)dt\].
\end{equation}
Combining  (\ref{eq6.4}) with (\ref{eq6.2}), we have
\begin{equation}\label{eq6.5}
E\[\int_0^T\Gamma_t(\alpha_1(t)+\Delta f(t))I_{E_\varepsilon}(t)dt\]+o(\varepsilon)\geq0.
\end{equation}
Let us define
$H_1(t,x,\nu,u,p,q,P):=pb(t,x,\nu,u)+q\sigma(t,x,\nu,u)+\frac{1}{2}P(\sigma(t,x,\nu,u)-\sigma(t,x,\nu,u_t^*))^2$, $(t,x,\nu,u,p,q,P)\in[0,T]\times\mathbb{R}\times\mathcal{P}_2(\mathbb{R})\times U\times\mathbb{R}^3$, and introduce the Hamiltonian
\begin{equation}\label{eq6.6}
\begin{aligned}
&\mathcal{H}(t,X_t^*,Y_t^*,Z_t^*,u_t,(p_1(t),p_2(t)),(q_{1,1}(t),q_{2,2}(t)),(P_1(t),P_2(t)))\\
&:=H_1(t,X_t^*,P_{X_t^*},u_t,p_1(t),q_{1,1}(t),P_1(t))
    +\widetilde{E}[H_1(t,\widetilde{X}_t^*,P_{X_t^*},\widetilde{u}_t,p_2(t),q_{2,2}(t),P_2(t))]\\
&\ \ \ \ +f(t,X_t^*,Y_t^*,Z_t^*+p_1(t)(\sigma(t,X_t^*,P_{X_t^*},u_t)-\sigma(t,X_t^*,P_{X_t^*},u_t^*)),P_{(X_t^*,Y_t^*)},u_t).
\end{aligned}
\end{equation}
The Hamiltonian is defined as an operator acting over spaces of square integrable random variables. Here, abusing notation we have written the Hamiltonian directly, evaluated at $((X_t^*,Y_t^*,Z_t^*),$ $(p_1(t),p_2(t)),(q_{1,1}(t),q_{2,2}(t)),(P_1(t),P_2(t))).$
Then, as
\begin{equation*}
\begin{aligned}
&(\alpha_1(t)+\Delta f(t))I_{E_\varepsilon}(t)\\
=&\(\mathcal{H}(t,X_t^*,Y_t^*,Z_t^*,u_t,(p_1(t),p_2(t)),(q_{1,1}(t),q_{2,2}(t)),(P_1(t),P_2(t)))\\
&-\mathcal{H}(t,X_t^*,Y_t^*,Z_t^*,u_t^*,(p_1(t),p_2(t)),(q_{1,1}(t),q_{2,2}(t)),(P_1(t),P_2(t)))\)I_{E_\varepsilon}(t),\ t\in[0,T],
\end{aligned}
\end{equation*}
and recalling that we have supposed $E_\varepsilon\subset\Gamma_M=\{t\in[0,T]:\ E[|Z_t^*|^2]\leq M\}$, we can choose $E_\varepsilon=[t_0,t_0+\varepsilon]\cap\Gamma_M,\ t_0\in[0,T)$, and the arbitrariness of $t_0\in[0,T)$ and of $\varepsilon>0$ allows to conclude from (\ref{eq6.5}) that, for all $u\in\mathcal{U}_{ad}$, $dtdP$-a.e.,
\begin{equation*}
\begin{aligned}
&\Gamma_t\(\mathcal{H}(t,X_t^*,Y_t^*,Z_t^*,u_t,(p_1(t),p_2(t)),(q_{1,1}(t),q_{2,2}(t)),(P_1(t),P_2(t)))\\
&-\mathcal{H}(t,X_t^*,Y_t^*,Z_t^*,u_t^*,(p_1(t),p_2(t)),(q_{1,1}(t),q_{2,2}(t)),(P_1(t),P_2(t)))\)I_{\Gamma_M}(t)\geq0,
\end{aligned}
\end{equation*}
for all $M\geq 1$. However, $I_{\cup_{M\geq1}\Gamma_M}(t)=1$, dt-a.e., and so we have
\begin{equation*}
\begin{aligned}
&\Gamma_t\(\mathcal{H}(t,X_t^*,Y_t^*,Z_t^*,u_t,(p_1(t),p_2(t)),(q_{1,1}(t),q_{2,2}(t)),(P_1(t),P_2(t)))\\
&-\mathcal{H}(t,X_t^*,Y_t^*,Z_t^*,u_t^*,(p_1(t),p_2(t)),(q_{1,1}(t),q_{2,2}(t)),(P_1(t),P_2(t)))\)\geq0,
\ dtdP\text{-a.e.},
\end{aligned}
\end{equation*}
for all $t\in[0,T]$.

This allows to formulate now the main result of this paper.
\begin{theorem}\label{thm6.1}
Let the assumptions \textbf{(A1)}-\textbf{(A4)} hold true. Let $u^*\in\mathcal{U}_{ad}$ be the optimal control and $(X^*,Y^*,Z^*)$ be the corresponding state processes of (\ref{eq2.5}). Then we have, for all $u\in\mathcal{U}_{ad}$, $dtdP$-a.e.,
\begin{equation}\label{eq6.7}
\begin{split}
&\mathcal{H}(t,X_t^*,Y_t^*,Z_t^*,u_t,(p_1(t),p_2(t)),(q_{1,1}(t),q_{2,2}(t)),(P_1(t),P_2(t)))\Gamma_t\\
&\geq \mathcal{H}(t,X_t^*,Y_t^*,Z_t^*,u_t^*,(p_1(t),p_2(t)),(q_{1,1}(t),q_{2,2}(t)),(P_1(t),P_2(t)))\Gamma_t,
\end{split}
\end{equation}
where $\mathcal{H}$ is defined by (\ref{eq6.6}), and $\Gamma$, $(p_1,(q_{1,1},q_{1,2}))$, $(p_2,(q_{2,1},q_{2,2}))$, $(P_1,(Q_{1,1},Q_{1,2}))$ and $(P_2,$ $(Q_{2,1},Q_{2,2}))$ are the solutions of (\ref{eq6.3}), (\ref{eq3.11})-(\ref{eq3.14}), respectively.
\end{theorem}
\begin{remark}
Our results extend those in \cite{HJX2022} to the case of mean-field coefficients which depend also on the law of the controlled processes $(X_t^u,Y_t^u)_{t\in[0,T]}$. In the case where $f$ does not depend on the marginal law of $Y_t^u$, i.e., $\widehat{f}^*_{\mu,y}(t)=0$, or even if we only have $\widehat{f}^*_{\mu,y}(t)\geq0$, then it follows from (\ref{eq6.3}) that $\Gamma_t>0,\ t\in[0,T]$, and (\ref{eq6.7}) writes
\begin{equation}\label{eq6.8}
\begin{aligned}
&\mathcal{H}(t,X_t^*,Y_t^*,Z_t^*,u_t,(p_1(t),p_2(t)),(q_{1,1}(t),q_{2,2}(t)),(P_1(t),P_2(t)))\\
&\geq \mathcal{H}(t,X_t^*,Y_t^*,Z_t^*,u_t^*,(p_1(t),p_2(t)),(q_{1,1}(t),q_{2,2}(t)),(P_1(t),P_2(t))),\ \
dtdP\text{-a.e.},
\end{aligned}
\end{equation}
for all $u\in\mathcal{U}_{ad}$. If, finally, none of the coefficients $b,\ \sigma,\ f$ and $\Phi$ depends on the law of the controlled state process $X^u$ nor of that of $Y^u$, then the solution of the adjoint BSDE (\ref{eq3.12}) annuls, $(p_2,(q_{2,1},q_{2,2}))=0$, and so does the solution $(P_2,(Q_{2,1},Q_{2,2}))$ of BSDE (\ref{eq3.14}). This has as consequence that in the definition of the Hamiltonian $\mathcal{H}$ it holds
$$\widetilde{E}\[H_1(t,\widetilde{X}_t^*,P_{X_t^*},\widetilde{u}_t,p_2(t),q_{2,2}(t),P_2(t))\]=0,$$
which makes that (\ref{eq6.8}) takes the equivalent form
\begin{equation*}
\begin{aligned}
&H_1(t,X_t^*,u,p_1(t),q_{1,1}(t),P_1(t))+f(t,X_t^*,Y_t^*,Z_t^*+p_1(t)(\sigma(t,X_t^*,u)-\sigma(t,X_t^*,u_t^*)),u)\\
&\geq H_1(t,X_t^*,u_t^*,p_1(t),q_{1,1}(t),P_1(t))+f(t,X_t^*,Y_t^*,Z_t^*,u_t^*),\ \text{for\ all}\  u\in U,\ dtdP\text{-a.s.}
\end{aligned}
\end{equation*}
This is exactly the optimality condition established in \cite{HJX2022} (Theorem 3.16).
\end{remark}

\section{Appendix}
\subsection{Proof of Proposition \ref{lem4.1}.}
\begin{proof} Without loss of generality we assume $b\equiv0$; the extension to the case with $b$ is straight-forward.

For $\phi\in \mathcal{H}^2(0,T)$, put
$$D[\phi](t):=\int_0^t\(\sigma_x(s)\phi_s+\widehat{E}[\widehat{\sigma}_{\mu}(s)\widehat{\phi}_s]\)dB_s.\ t\in[0,T].$$
Then $D:\mathcal{H}^2(0,T)\rightarrow \mathcal{H}^2(0,T)$ is a bounded linear operator. In fact, for all $\phi\in \mathcal{H}^2(0,T)$, it holds
\begin{equation*}
\begin{split}
|D[\phi]|^2_{\mathcal{H}^2(0,T)}
=&E\[\int_0^T|D[\phi](t)|^2dt\]
=\int_0^TE[|D[\phi](t)|^2]dt\\
\leq &C\int_0^T\int_0^tE[|\phi_s|^2]dsdt
\leq CT|\phi|^2_{\mathcal{H}^2(0,T)}.
\end{split}
\end{equation*}
Denoting by $D^j=D\circ\cdots \circ D$ the composition of $D$ with itself up to the order $j$, we have
\begin{equation*}
\begin{aligned}
&|D^j[\phi]|^2_{\mathcal{H}^2(0,T)}
=E\[\int_0^T|D^j[\phi](t)|dt\]
\leq C\int_0^T\int_0^tE[|D^{j-1}[\phi](s)|^2]dsdt\\
&\leq \cdots \leq C^j\int_0^T\int_0^{t}\int_0^{t_1}\cdots\int_0^{t_{j-1}}E[|\phi_{t_j}|^2]dt_j\cdots dt_2dt_1dt\\
&\leq C^j\int_0^T\int_0^{t}\int_0^{t_1}\cdots\int_0^{t_{j-2}}dt_{j-1}\cdots dt_2dt_1dt
  \int_0^TE[|\phi_{t}|^2]dt
=\frac{(CT)^j}{j!}|\phi|^2_{\mathcal{H}^2(0,T)},\ \phi\in \mathcal{H}^2(0,T).
\end{aligned}
\end{equation*}
Thus the operator norm: $\|D^j\|^2\leq \frac{(CT)^j}{j!},\ j\geq 1$. Let $A_n=\sum\limits^n_{j=1}D^j$, for $m>n$, we have $\|A_m-A_n\|\leq \sum\limits_{j=n}^m\|D^j\|\leq\sum\limits_{j=n}^m\frac{(CT)^j}{j!}\rightarrow0$, as $n,m\rightarrow\infty$, so there exists $A:\mathcal{H}^2(0,T)\rightarrow\mathcal{H}^2(0,T)$ such that $A_n$ converges to $A$ in the space $L(\mathcal{H}^2(0,T))$. Hence, $A\in L(\mathcal{H}^2(0,T))$, let us determine its adjoint operator $A^*\in L(\mathcal{H}^2(0,T))$.

Given $\phi,\psi\in\mathcal{H}^2(0,T)$, let $\vartheta^\psi=(\vartheta^\psi_{t,s})_{0\leq s\leq t}$ be the jointly measurable two-parameter process with $\vartheta_{t,\cdot}^\psi:=\vartheta^\psi_{t,\cdot}I_{[0,t]}\in \mathcal{H}^2(0,T),\ t\in[0,T]$, such that
$$\psi_t=E[\psi_t]+\int_0^t\vartheta^\psi_{t,s}dB_s,\ P\text{-a.s},\ t\in[0,T].$$
Note that $\vartheta_{\cdot,s}=(\vartheta_{t,s}^\psi)_{t\in[s,T]}$ is $\mathcal{B}([s,T])\otimes\mathcal{F}_s$-measurable, $s\in[0,T]$. Thus,
\begin{equation}\label{eqA.2}
\ E\[\int_0^T\int_0^t|\vartheta^\psi_{t,s}|^2dsdt\]
\leq  \int_0^TE\[\int_0^t|\vartheta_{t,s}^\psi|^2ds\]dt
\leq \int_0^TE[\psi_t^2]dt
\leq |\psi|^2_{\mathcal{H}^2(0,T)}.
\end{equation}
Therefore, for the adjoint operator $D^*,\ A^*\in L(\mathcal{H}^2(0,T))$, we have
\begin{equation}\label{eqA.1}
\begin{aligned}
&(D^*[\psi],\phi)_{\mathcal{H}^2(0,T)}=(\psi,D[\phi])_{\mathcal{H}^2(0,T)}
=\int_0^TE\[\psi_t\int_0^t\(\sigma_x(s)\phi_s+\widehat{E}[\widehat{\sigma}_{\mu}(s)\widehat{\phi}_s]\)dB_s\]dt\\
&=\int_0^TE\[\int_0^t\vartheta_{t,s}^\psi\(\sigma_x(s)\phi_s+\widehat{E}[\widehat{\sigma}_{\mu}(s)\widehat{\phi}_s]\)ds\]dt
=\int_0^T\int_0^tE\[\vartheta_{t,s}^\psi\(\sigma_x(s)\phi_s+\widehat{E}[\widehat{\sigma}_{\mu}(s)\widehat{\phi}_s]\)\]dsdt\\
&=\int_0^T\int_0^tE\[\(\vartheta_{t,s}^\psi\sigma_x(s)+\widehat{E}[\widehat{\sigma}^*_{\mu}(s)\widehat{\vartheta}_{t,s}^\psi]\)\phi_s\]dsdt
=E\[\int_0^T\int_s^T\(\vartheta_{t,s}^\psi\sigma_x(s)+\widehat{E}[\widehat{\sigma}^*_{\mu}(s)\widehat{\vartheta}_{t,s}^\psi]\)dt\cdot \phi_sds\].
\end{aligned}
\end{equation}
Thus, defining the jointly measurable process $\eta^\psi$:
$$\eta_s^\psi:=\int_s^T(\vartheta^\psi_{t,s}\sigma_x(s)+\widehat{E}[\widehat{\sigma}_\mu^*(s)\widehat{\vartheta}_{t,s}^\psi])dt,$$
we deduce from (\ref{eqA.1}) that
$$D^*[\psi](s)=\eta_s^\psi=
\int_s^T(\vartheta_{t,s}^\psi\sigma_x(s)+\widehat{E}[\widehat{\sigma}^*_\mu(s)\widehat{\vartheta}_{t,s}^\psi])dt,
\ dsdP\text{-a.e.}$$
Then, from (\ref{eqA.2}),
$$|D^*[\psi]|^2_{\mathcal{H}^2(0,T)}\leq C|\psi|^2_{\mathcal{H}^2(0,T)}.$$
Now we consider the adjoint operator $A^*$ of $A$,
\begin{equation*}
\begin{aligned}
(A^*[\psi],\phi)_{\mathcal{H}^2(0,T)}&=(\psi,A[\phi])_{\mathcal{H}^2(0,T)}
=\lim_{N\rightarrow\infty}\sum_{j=0}^N(\psi,D^j[\phi])_{\mathcal{H}^2(0,T)}\\
&=\lim_{N\rightarrow\infty}\sum_{j=0}^N((D^*)^j[\psi],\phi)_{\mathcal{H}^2(0,T)}
=\lim_{N\rightarrow\infty}(\sum_{j=0}^N(D^*)^j[\psi],\phi)_{\mathcal{H}^2(0,T)}.
\end{aligned}
\end{equation*}
Hence, $A^*[\psi]=\lim_{N\rightarrow\infty}\sum_{j=0}^N(D^*)^j[\psi]$ with respect to the weak topology $\sigma(\mathcal{H}^2(0,T),\mathcal{H}^2(0,T))$. As $(D^N)^*=(D^*)^N$, and $\|(D^N)^*\|=\|D^N\|\leq \frac{(CT)^N}{N!}$, this convergence holds to also in $\mathcal{H}^2(0,T)$.

Now we come back to (\ref{eq3.2}). First of all, let us consider the stochastic integral process $G_t^\varepsilon:=\displaystyle\int_0^t\delta\sigma(s)I_{E_\varepsilon}(s)dB_s,\ t\in[0,T]$. It is obvious that $G^\varepsilon\in \mathcal{H}^2(0,T)$. Moreover, from equation (\ref{eq3.2}) for $X^1$ and the definition of $D$, we have
\begin{equation*}
\begin{aligned}
X^1=&D[X^1]+G^\varepsilon=G^\varepsilon+D[G^\varepsilon+D[X^1]]=G^\varepsilon+D[G^\varepsilon]+D^2[X^1]\\
=&\cdots=\sum_{j=0}^ND^j(G^\varepsilon)+D^{N+1}[X^1],\ \ N\geq1.
\end{aligned}
\end{equation*}
We notice that  $|D^{N+1}[X^1]|_{\mathcal{H}^2(0,T)}\leq\|D^{N+1}\|\cdot |X^1|_{\mathcal{H}^2(0,T)}\rightarrow0$ and $\sum\limits^N_{j=0}D^j\rightarrow A$ in $L(\mathcal{H}^2(0,T))$, as $N\rightarrow\infty$. Hence we have $X^1=A(G^\varepsilon)$. Next, we introduce the identical mapping $I(\phi):=\phi,\ \phi\in \mathcal{H}^2(0,T)$, and observe that $A=\sum_{j\geq0}D^j=I+D\circ A$. So we get $X^1=A(G^\varepsilon)=G^\varepsilon+D[A[G^\varepsilon]]$. Let $\theta\in \mathcal{H}^2(0,T)$ and $\Gamma_\theta=\{t\in[0,T]:\ E[\theta_t^2]<\infty\}$. Note that $\Gamma_\theta$ is of full Borel measure. Similar to the previous notation and representation, we get
$$\theta_t=E[\theta_t]+\int_0^t\vartheta_{t,s}^\theta dB_s,\ P\text{-a.s.},\ t\in\Gamma_\theta.$$
We also remark that
\begin{equation}\label{eqA.3}
E\[\int_0^t|\vartheta_{t,s}^\theta|^2ds\]\leq E[\theta_t^2],\ t\in\Gamma_\theta.
\end{equation}
Then, for all $t\in\Gamma_\theta$,
\begin{equation}\label{eqA.4}
\begin{aligned}
&E[\theta_t X_t^1]=E[\theta_t(G_t^\varepsilon+D[A[G^\varepsilon]](t))]\\
&=E\[\theta_t\int_0^t\(\delta\sigma(s)I_{E_\varepsilon}(t)
        +\{\sigma_x(s)A[G^\varepsilon](s)+\widehat{E}[\widehat{\sigma}_\mu(s)\widehat{A[G^\varepsilon](s)}]\}\)dB_s\]\\
&=E\[\int_0^t\vartheta^\theta_{t,s}\delta\sigma(s)I_{E_\varepsilon}(s) ds\]
    +\int_0^tE\[\(\vartheta_{t,s}^\theta\sigma_x(s)+\widehat{E}[\widehat{\vartheta}_{t,s}^\theta\widehat{\sigma}_\mu^*(s)]\)A[G^\varepsilon](s)\]ds.
\end{aligned}
\end{equation}
Setting, for $t\in\Gamma_\theta$,
\begin{equation*}
\eta_{t,s}^\theta:=
\left\{
\begin{split}
&\vartheta_{t,s}^\theta\sigma_x(s)+\widehat{E}[\widehat{\vartheta}_{t,s}^\theta\widehat{\sigma}_\mu^*(s)],\ &0\leq s\leq t\leq T,\\
&0, &s>t,
\end{split}
\right.
\end{equation*}
we get from (\ref{eqA.3}) that $E\[\displaystyle\int_0^t|\eta_{t,s}^\theta|^2ds\]\leq CE\[\displaystyle\int_0^t|\vartheta_{t,s}^\theta|^2ds\]\leq CE[|\theta_t|^2],\ t\in\Gamma_\theta,$ and (\ref{eqA.4}) can be rewritten as
\begin{equation}\label{eqA.5}
E[\theta_tX_t^1]=E\[\int_0^t \vartheta_{t,s}^\theta\delta\sigma(s)I_{E_\varepsilon}(s)ds\]
    +(\eta_{t,\cdot}^\theta,A[G^\varepsilon])_{\mathcal{H}^2(0,T)},\ t\in\Gamma_\theta.
\end{equation}
Notice that
$$(\eta^\theta_{t,\cdot},A[G^\varepsilon])_{\mathcal{H}^2(0,T)}
=(A^*[\eta_{t,\cdot}],G^\varepsilon)_{\mathcal{H}^2(0,T)}
=E\[\int_0^T A^*[\eta_{t,\cdot}^\theta](s)G_s^\varepsilon ds\].$$
Consequently, for $t\in\Gamma_\theta$, with the notation $u_{t,s}^\theta:=A^*[\eta_{t,\cdot}^\theta](s)$, we get for $\vartheta_{t,s,\cdot}^{u^\theta}=\vartheta_{t,s,\cdot}^{u^\theta}I_{[0,s]}\in {\mathcal{H}^2(0,T)}$ with
$$u_{t,s}^\theta=E[u_{t,s}^\theta]+\int_0^s\vartheta_{t,s,r}^{u^\theta}dB_r,\ P\text{-a.s.},$$
that
\begin{equation}\label{eqA.6}
(\eta_{t,\cdot}^\theta,A[G^\varepsilon])_{\mathcal{H}^2(0,T)}
=\int_0^tE[u_{t,s}^\theta G_s^\varepsilon]ds
=\int_0^tE\[\int_0^s\vartheta_{t,s,r}^{u^\theta}\cdot \delta\sigma(r)I_{E_\varepsilon}(r)dr\]ds
\end{equation}
and, thus,
\begin{equation}\label{eqA.7}
|(\eta_{t,\cdot}^\theta, A[G^\varepsilon])_{\mathcal{H}^2(0,T)}|
\leq C\(E\[\int_0^t\int_{E_\varepsilon\cap[0,s]}|\vartheta_{t,s,r}^{u^\theta}|^2drds\]\)^\frac{1}{2}
\cdot\sqrt{\varepsilon},\ \ \varepsilon>0.
\end{equation}
We observe that
\begin{equation*}
\begin{aligned}
&E\[\int_0^T\int_0^s|\vartheta_{t,s,r}^{u^\theta}|^2drds]
\leq E\[\int_0^T|u_{t,s}^\theta|^2ds\]
=\Bigr|A^*[\eta^\theta_{t,\cdot}]\Bigr|^2_{\mathcal{H}^2(0,T)}
\leq \|A^*\|^2\cdot |\eta_{t,\cdot}^\theta|^2_{\mathcal{H}^2(0,T)}\\
&\leq C\|A\|^2\cdot E\[\int_0^t|\vartheta_{t,s}^\vartheta|^2ds\]
\leq C\|A\|^2\cdot E[|\theta_t|^2],\ t\in\Gamma_\theta,
\end{aligned}
\end{equation*}
and thus,
$$E\[\int_0^T\int_0^T\int_0^s|\vartheta_{t,s,r}^{u^\theta}|^2drdsdt\]
\leq C\|A\|^2E\[\int_0^T|\theta_t|^2dt\].$$
Consequently, from the dominated convergence theorem,
$$\rho(\varepsilon):= CE\[\int_0^T \int_0^T\int_{E_\varepsilon\cap[0,s]}|\vartheta_{t,s,r}^{u^\theta}|^2drdsdt\]\rightarrow0,
\ \text{as}\ \varepsilon\downarrow0,$$
and combining this with (\ref{eqA.7}), we have
$$\int_0^T|(\eta_{t,\cdot}^\theta,A[G^\varepsilon])_{\mathcal{H}^2(0,T)}|^2dt\leq \rho^1(\varepsilon)\varepsilon,\ \varepsilon>0.$$
Similarly,
$$\int_0^T\Bigr|E\[\int_0^t\vartheta_{t,s}^\theta\cdot\delta\sigma(s)I_{E_\varepsilon}(s)ds\]\Bigr|^2dt
\leq C E\[\int_0^T\int_{E_\varepsilon\cap[0,t]}|\vartheta^\theta_{t,s}|^2dsdt\]\cdot \varepsilon
=\rho^2(\varepsilon)\varepsilon,\ \varepsilon>0,$$
where $\rho^2(\varepsilon)
:=CE\[\displaystyle\int_0^T\int_{E_\varepsilon\cap[0,t]}|\vartheta^\theta_{t,s}|^2dsdt\]\downarrow0$, as $\varepsilon\downarrow0$, as also here the dominated convergence theorem can be applied. Indeed, $E\[\displaystyle\int_0^T\int_0^t|\vartheta_{t,s}^\theta|^2dsdt\]\leq E\[\displaystyle\int_0^T|\theta_t|^2dt\]<+\infty$.

We define now $\rho(\varepsilon):=\rho^1(\varepsilon)+\rho^2(\varepsilon)$. Then, from (\ref{eqA.5}) and the estimates above, we have
$$\int_0^T|E[\theta_tX_t^1]|^2dt\leq \rho(\varepsilon)\varepsilon,\ \varepsilon>0,$$
where $\rho(\varepsilon)\rightarrow0$, as $\varepsilon\downarrow0$. \end{proof}

\textbf{Acknowledgments.} \ The authors would like to thank Professors Shaolin Ji and Mingshang Hu for the fruitful discussions.

\end{document}